%!TEX TS-program = xelatex
% BEGIN ARXIV XELATEX GUARD
\ExplSyntaxOn
\msg_new:nnn { tex-studies-arxiv } { xelatex-required }
{ This~arXiv~wrapper~requires~XeLaTeX.~Run~latexmk~-xelatex. }
\sys_if_engine_xetex:F
{ \msg_fatal:nn { tex-studies-arxiv } { xelatex-required } }
\ExplSyntaxOff
% END ARXIV XELATEX GUARD
% Generated by export-arxiv.ps1. Edit content.tex or the yj-* stack, then rerun.
% The yj-* class/package code below is materialized from the current local sources.
% Materialized from yj-paper.cls
\documentclass[12pt, reqno]{amsart}
\makeatletter
\ExplSyntaxOn

% Fix \part: center (amsart uses \raggedright), Roman numerals, bold number
\ExplSyntaxOff

\def
\part{\@startsection{part}{0}%
  \z@{\linespacing\@plus\linespacing}{.5\linespacing}%
{\normalfont\bfseries\centering}}
\def\yj@part{part}
\def\@seccntformat#1{%
  \expandafter\ifx\csname yj@#1\endcsname\yj@part
  \csname the#1\endcsname\@secnumpunct
  \else
  \protect\textup{\protect\@secnumfont
    \csname the#1\endcsname
    \protect\@secnumpunct
  }%
  \fi
}
\ExplSyntaxOn

% --- Begin yj-bib.sty ---
% Materialized from yj-bib.sty
\bool_new:N \g__yjbib_chinese_bool

\keys_define:nn { yj-bib }
{
  chinese .bool_gset:N = \g__yjbib_chinese_bool,
  chinese .default:n   = true,
}

\usepackage[
  backend=biber,
  datamodel=yj-amsalpha,
  backref=true,
  bibstyle=alphabetic,
  citestyle=yj-amsalpha,
  isbn=false,
  doi=false,
  url=false,
  giveninits=true,
  maxnames=4,
  minnames=4,
  date=year,
  labeldate=year
]{biblatex}
\usepackage{biblatex-shortfields}
% Load the AMS modifier after biblatex-shortfields so the policy layer selects
% short fields while the visual style retains its series punctuation, for
% example "GTM, no. 52". Keeping this as a package also lets export-arxiv.ps1
% inline the local formatting without shipping a separate .bbx file.

% --- Begin yj-amsalpha.sty ---
% This file may be distributed and/or modified under the conditions of the
% LaTeX Project Public License, version 1.3c or later.
% Materialized from yj-amsalpha.sty
\newtoggle{bbx:mr}
\DeclareBibliographyOption[boolean]{mr}[true]{\settoggle{bbx:mr}{#1}}

\ExecuteBibliographyOptions{
  abbreviate=true,
  doi=false,
  eprint=true,
  isbn=false,
  maxalphanames=4,
  maxbibnames=999,
  minalphanames=3,
  mr=true,
  related=false,
  url=false
}

% AMS alphabetic labels retain all four names, but abbreviate longer lists to
% three initials followed by a raised plus. AMS purifies the plus out of its
% sorting key, so it must not participate in Biber's corresponding key.

% Keep a fixed two-character date component in Biber's internal label even
% without a year. The printable label below removes the "nd" sentinel. This
% lets sorting discard the date component without truncating undated names,
% while keeping undated labels distinct from labels for years ending in 00.
\DeclareLabelalphaTemplate{
  \labelelement{
    \field[final]{shorthand}
    \field{label}
    \field[strwidth=3,strside=left,ifnames=1]{labelname}
    \field[strwidth=1,strside=left]{labelname}
  }
  \labelelement{
    \field[strwidth=2,strside=right]{year}
    \literal{nd}
  }
}
\tl_new:N \l__yjamsalpha_label_tl
\DeclareFieldFormat{labelalpha}{%
  \ifboolexpr{
    test {\iffieldundef{year}}
    and
    test {\iffieldundef{shorthand}}
  }
  {\savefield{labelalpha}{\l__yjamsalpha_label_tl}%
  \exp_args:NV \tl_range:nnn \l__yjamsalpha_label_tl {1}{-3}}
{#1}}
% BibLaTeX also applies the labelalpha formatter when measuring labelalphawidth,
% so the internal sentinel never adds space to either citations or labels.

% AMS sorts by the name abbreviation and four-digit year before the full
% name/year/title tie-breakers. Sorting anyt instead compares the printed
% two-digit year first, which reverses 1999 and 2001. Preserve explicit modern
% sorting fields; a shorthand is a complete label and has no date to discard.
\DeclareSortingTemplate{yjamsalpha}{
  \sort{\field{presort}}
  \sort[final]{\field{sortkey}}
  \sort{
    \field{shorthand}
    \field[strwidth=-2,strside=left]{labelalpha}
  }
  \sort{
    \field{sortyear}
    \field{year}
    \literal{0000}
  }
  \sort{
    \field{sortname}
    \field{author}
    \field{editor}
    \field{translator}
    \field{sorttitle}
    \field{title}
  }
  \sort{
    \field{sortyear}
    \field{year}
  }
  \sort{
    \field{sorttitle}
    \field{title}
  }
  \sort{
    \field{volume}
    \literal{0}
  }
}
\ExecuteBibliographyOptions{sorting=yjamsalpha}

\renewcommand*{\subtitlepunct}{\addcolon\space}
\setlength{\biblabelsep}{\labelsep}

\DeclareFieldFormat{yjams:titlecase}{\MakeSentenceCase*{#1}}
\DeclareFieldFormat{journaltitle}{#1\isdot}
\DeclareFieldFormat{booktitle}{#1}
\DeclareFieldFormat[article,periodical]{volume}{\mkbibbold{#1}}
\DeclareFieldFormat[article,periodical]{pages}{#1}

% Let the URL formatter supply breakpoints within long arXiv identifiers and
% retain a natural breakpoint before the optional subject classification. Use
% an explicit URI action because \href's file heuristic misclassifies legacy
% identifiers such as hep-th/0012136 when the linked text spans multiple lines.
\DeclareFieldFormat{eprint:arxiv}{%
  arXiv\addcolon\space
  \ifhyperref
  {\hyper@linkurl{%
      \nolinkurl{#1}%
      \iffieldundef{eprintclass}
      {}
    {\addspace\texttt{\mkbibbrackets{\thefield{eprintclass}}}}}%
  {https://arxiv.org/\abx@arxivpath/#1}}
  {\nolinkurl{#1}%
    \iffieldundef{eprintclass}
    {}
{\addspace\texttt{\mkbibbrackets{\thefield{eprintclass}}}}}}
\DeclareFieldAlias{eprint:arXiv}{eprint:arxiv}

% AMS Style Guide §10.1 puts the arXiv identifier before the date. Apply the
% same order to Hexagon preprints. Include legacy exports as misc and articles
% with no journal metadata, independently of identifier display options.
\newcommand*{\yjams@ifstandalonepreprint}{%
  \ifboolexpr{
    not test {\iffieldundef{eprint}}
    and
    (test {\iffieldequalstr{eprinttype}{arxiv}}
      or test {\iffieldequalstr{eprinttype}{arXiv}}
      or test {\iffieldequalstr{eprinttype}{hexagon}}
    or test {\iffieldequalstr{eprinttype}{Hexagon}})
    and
    (test {\ifentrytype{online}}
      or test {\ifentrytype{misc}}
      or (test {\ifentrytype{article}}
        and test {\iffieldundef{journaltitle}}
        and test {\iffieldundef{journalsubtitle}}
        and test {\iffieldundef{journaltitleaddon}}
        and test {\iffieldundef{shortjournal}}
        and test {\iffieldundef{issn}}
        and test {\iffieldundef{volume}}
        and test {\iffieldundef{number}}
and test {\iffieldundef{eid}}))}}

\renewbibmacro*{date}{%
  \yjams@ifstandalonepreprint
  {\iftoggle{bbx:eprint}
    {\usebibmacro{eprint}%
    \newunit}
  {}}
  {}%
\printdate}

% The date macro has already handled the preprint eprint. Keep the native
% terminal DOI and URL handling, and leave the eprint field available.
\letbibmacro{yjams:doi+eprint+url}{doi+eprint+url}
\renewbibmacro*{doi+eprint+url}{%
  \yjams@ifstandalonepreprint
  {\begingroup
    \togglefalse{bbx:eprint}%
    \usebibmacro{yjams:doi+eprint+url}%
  \endgroup}
{\usebibmacro{yjams:doi+eprint+url}}}

% Hexagon resolves both unversioned and versioned identifiers at the site root.
\DeclareFieldFormat{eprint:hexagon}{%
  hexagon\addcolon\space
  \ifhyperref
  {\hyper@linkurl{%
      \nolinkurl{#1}%
      \iffieldundef{eprintclass}
      {}
    {\addspace\texttt{\mkbibbrackets{\thefield{eprintclass}}}}}%
  {https://hexagonmath.org/#1}}
  {\nolinkurl{#1}%
    \iffieldundef{eprintclass}
    {}
{\addspace\texttt{\mkbibbrackets{\thefield{eprintclass}}}}}}
\DeclareFieldAlias{eprint:Hexagon}{eprint:hexagon}

% Apply the same identifier break policy to direct and wrapper loading paths.
% Canonical identifiers remain unchanged; only their legal breakpoints vary.
\setcounter{biburlnumpenalty}{100}
\biburlnumskip=0mu plus 1mu\relax

\DeclareFieldFormat
[article,periodical,book,mvbook,inbook,bookinbook,suppbook,
  collection,mvcollection,incollection,suppcollection,
  proceedings,mvproceedings,inproceedings,
reference,mvreference,inreference]
{number}{\bibstring{number}\addnbspace #1}

\DeclareFieldFormat
[article,inbook,incollection,inproceedings,patent,thesis,unpublished]
{title}{\mkbibemph{#1}}

% AMS and article bibliographies relax line breaking inside the list. Apply
% the same fully justified policy to direct style users and wrapper users.
\AtBeginBibliography{\sloppy}

% Keep BibLaTeX's name hooks, prefix handling and giveninits selection, adding
% only the AMS comma before a suffix in the ordinary given-family name order.
\renewbibmacro*{name:given-family}[4]{%
  \usebibmacro{name:delim}{#2#3#1}%
  \usebibmacro{name:hook}{#2#3#1}%
  \mkbibcompletenamegivenfamily{%
    \ifdefvoid{#2}
    {}
    {\mkbibnamegiven{#2}\isdot\bibnamedelimd}%
    \ifdefvoid{#3}
    {}
    {\mkbibnameprefix{#3}\isdot
      \ifprefchar
      {}
    {\ifuseprefix{\bibnamedelimc}{\bibnamedelimd}}}%
    \mkbibnamefamily{#1}\isdot
\ifdefvoid{#4}{}{\addcomma\bibnamedelimd\mkbibnamesuffix{#4}\isdot}}}

% Main editors take a parenthesized role. Container editors use editorstrg
% directly, retaining their distinct (Name, ed.) punctuation.
\newbibmacro*{yjams:editor+others}{%
  \ifboolexpr{
    test \ifuseeditor
    and
    not test {\ifnameundef{editor}}
  }
  {\printnames{editor}%
    \setunit{\addspace}%
    \printtext[parens]{\usebibmacro{editor+othersstrg}}%
  \clearname{editor}}
{}}
\renewbibmacro*{editor+others}{\usebibmacro{yjams:editor+others}}

% The language list has already undergone BibLaTeX's clearlang processing;
% preserve its localized list formatter and avoid empty parentheses.
\newbibmacro*{yjams:language}{%
  \iflistundef{language}
  {}
  {\setunit{\addspace}%
\printtext[parens]{\printlist{language}}}}

\renewbibmacro*{in:}{}

\renewbibmacro*{title}{%
  \ifboolexpr{
    test {\iffieldundef{title}}
    and
    test {\iffieldundef{subtitle}}
  }
  {}
  {\printtext[title]{%
      \printfield[yjams:titlecase]{title}%
      \setunit{\subtitlepunct}%
    \printfield[yjams:titlecase]{subtitle}}%
  \newunit}%
\printfield{titleaddon}}

\renewbibmacro*{maintitle}{%
  \ifboolexpr{
    test {\iffieldundef{maintitle}}
    and
    test {\iffieldundef{mainsubtitle}}
  }
  {}
  {\printtext[maintitle]{%
      \printfield[yjams:titlecase]{maintitle}%
      \setunit{\subtitlepunct}%
    \printfield[yjams:titlecase]{mainsubtitle}}%
  \newunit}%
\printfield{maintitleaddon}}

% BibLaTeX's current data model distinguishes a work's volume from its
% number within a publication series. Prefer the abbreviated series title
% when supplied and never reinterpret volume as a series number.
\renewbibmacro*{series+number}{%
  \iffieldundef{shortseries}
  {\printfield[noformat]{series}}
  {\printfield[noformat]{shortseries}}%
  \setunit*{\addcomma\space}%
  \printfield{number}%
\newunit}

\renewbibmacro*{maintitle+title}{%
  \iffieldsequal{maintitle}{title}
  {\clearfield{maintitle}%
    \clearfield{mainsubtitle}%
  \clearfield{maintitleaddon}}
  {\iffieldundef{maintitle}
    {}
    {\usebibmacro{maintitle}%
      \newunit\newblock
      \iffieldundef{volume}
      {}
      {\printfield{volume}%
        \printfield{part}%
  \setunit{\subtitlepunct}}}}%
  \usebibmacro{title}%
\newunit}

\renewbibmacro*{maintitle+booktitle}{%
  \iffieldundef{maintitle}
  {}
  {\usebibmacro{maintitle}%
    \newunit\newblock
    \iffieldundef{volume}
    {}
    {\printfield{volume}%
      \printfield{part}%
  \setunit{\subtitlepunct}}}%
  \usebibmacro{booktitle}%
\newunit}

\renewbibmacro*{journal+issuetitle}{%
  \usebibmacro{journal}%
  \setunit*{\addspace}%
  \printfield{volume}%
  \yjams@ifstandalonepreprint
  {\newunit\usebibmacro{date}}
  {\setunit{\addspace}%
    \iffieldundef{year}
    {}
  {\printtext[parens]{\printfield{year}}}}%
  \newunit
  \printfield{number}%
  \setunit{\bibeidpunct}%
  \printfield{eid}%
  \ifboolexpr{
    test {\ifentrytype{article}}
    and
    not test {\iffieldundef{eid}}
    and
    test {\iffieldundef{pages}}
  }
  {\setunit*{\addcomma\space}%
  \printfield{pagetotal}}
  {}%
\newunit}

% Original amsalpha puts article pagination before the note. Keep the
% standard note-before-pages order for every other entry type.
\renewbibmacro*{note+pages}{%
  \ifentrytype{article}
  {\printfield{pages}%
    \usebibmacro{yjams:language}%
    \newunit
  \printfield{note}}
  {\printfield{note}%
    \setunit{\bibpagespunct}%
  \printfield{pages}}%
\newunit}

\newbibmacro*{yjams:publisher+location}{%
  \printlist{publisher}%
  \setunit*{\addcomma\space}%
  \printlist{location}%
\newunit}

\renewbibmacro*{publisher+location+date}{%
  \usebibmacro{yjams:publisher+location}%
  \usebibmacro{date}%
\newunit}

\renewbibmacro*{institution+location+date}{%
  \printlist{institution}%
  \setunit*{\addcomma\space}%
  \printlist{location}%
  \setunit*{\addcomma\space}%
  \usebibmacro{date}%
\newunit}

\newbibmacro*{yjams:container}{%
  \usebibmacro{maintitle+booktitle}%
  \ifboolexpr{
    test {\iffieldundef{venue}}
    and
    test {\iffieldundef{eventyear}}
  }
  {}
  {\setunit{\addspace}%
    \printtext[parens]{%
      \printfield{venue}%
      \setunit*{\addcomma\space}%
  \printeventdate}}%
  \ifnameundef{editor}
  {}
  {\setunit{\addspace}%
    \printtext[parens]{%
      \printnames{editor}%
      \setunit{\printdelim{editortypedelim}}%
  \usebibmacro{editorstrg}}}%
\newunit}

\DeclareFieldFormat{mrnumber}{%
  \ifhyperref
  {\hyper@linkurl{MR\addnbspace #1}{https://mathscinet.ams.org/mathscinet-getitem?mr=#1}}
{MR\addnbspace #1}}
% A MathSciNet number identifies the reviewed work itself.
\DeclareDataInheritance{*}{*}{\noinherit{mrnumber}}

% Standard and custom drivers share this hook before backreferences.
\renewbibmacro*{addendum+pubstate}{%
  \printfield{addendum}%
  \newunit\newblock
  \printfield{pubstate}%
  \ifboolexpr{
    not togl {bbx:mr}
    or test {\iffieldundef{mrnumber}}
  }
  {}
  {\yjams@ifstandalonepreprint
    {}
    {\setunit{\addperiod\space}%
      \printfield{mrnumber}%
      \ifboolexpr{
        test {\iflistundef{pageref}}
        and (not togl {bbx:related} or test {\iffieldundef{related}})
      }
      {\nopunct}
{}}}}

% Keep the standard identifier, publication-state and linking facilities at
% the end of each entry while restoring the original AMS field order.
\newbibmacro*{yjams:finentry}{%
  \usebibmacro{doi+eprint+url}%
  \newunit\newblock
  \usebibmacro{addendum+pubstate}%
  \setunit{\bibpagerefpunct}\newblock
  \usebibmacro{pageref}%
  \newunit\newblock
  \iftoggle{bbx:related}
  {\usebibmacro{related:init}%
  \usebibmacro{related}}
  {}%
\usebibmacro{finentry}}

\DeclareBibliographyDriver{article}{%
  \usebibmacro{bibindex}%
  \usebibmacro{begentry}%
  \usebibmacro{author/translator+others}%
  \setunit{\printdelim{nametitledelim}}\newblock
  \usebibmacro{title}%
  \newunit
  \ifentrytype{article}
  {}
  {\printlist{language}}%
  \newunit\newblock
  \usebibmacro{byauthor}%
  \newunit\newblock
  \usebibmacro{bytranslator+others}%
  \newunit\newblock
  \printfield{version}%
  \newunit\newblock
  \usebibmacro{in:}%
  \usebibmacro{journal+issuetitle}%
  \newunit
  \usebibmacro{byeditor+others}%
  \newunit
  \usebibmacro{note+pages}%
  \newunit\newblock
  \iftoggle{bbx:isbn}
  {\printfield{issn}}
  {}%
  \newunit\newblock
\usebibmacro{yjams:finentry}}

\DeclareBibliographyDriver{book}{%
  \usebibmacro{bibindex}%
  \usebibmacro{begentry}%
  \usebibmacro{author/editor+others/translator+others}%
  \setunit{\printdelim{nametitledelim}}\newblock
  \usebibmacro{maintitle+title}%
  \newunit
  \usebibmacro{byauthor}%
  \newunit\newblock
  \usebibmacro{byeditor+others}%
  \newunit\newblock
  \printfield{edition}%
  \newunit
  \iffieldundef{maintitle}
  {\printfield{volume}%
  \printfield{part}}
  {}%
  \newunit
  \printfield{volumes}%
  \newunit\newblock
  \usebibmacro{series+number}%
  \newunit\newblock
  \usebibmacro{chapter+pages}%
  \newunit\newblock
  \usebibmacro{publisher+location+date}%
  \usebibmacro{yjams:language}%
  \newunit
  \printfield{note}%
  \newunit\newblock
\usebibmacro{yjams:finentry}}

% An edited legacy book maps to collection in BibLaTeX. Keep the same AMS
% ordering while retaining the native collection extent and identifier fields.
\DeclareBibliographyDriver{collection}{%
  \usebibmacro{bibindex}%
  \usebibmacro{begentry}%
  \usebibmacro{editor+others}%
  \setunit{\printdelim{nametitledelim}}\newblock
  \usebibmacro{maintitle+title}%
  \newunit
  \usebibmacro{byauthor}%
  \newunit\newblock
  \usebibmacro{byeditor+others}%
  \newunit\newblock
  \printfield{edition}%
  \newunit
  \iffieldundef{maintitle}
  {\printfield{volume}%
  \printfield{part}}
  {}%
  \newunit
  \printfield{volumes}%
  \newunit\newblock
  \usebibmacro{series+number}%
  \newunit\newblock
  \usebibmacro{chapter+pages}%
  \newunit\newblock
  \usebibmacro{publisher+location+date}%
  \usebibmacro{yjams:language}%
  \newunit
  \printfield{note}%
  \newunit
  \printfield{pagetotal}%
  \newunit\newblock
  \iftoggle{bbx:isbn}
  {\printfield{isbn}}
  {}%
  \newunit\newblock
\usebibmacro{yjams:finentry}}

\DeclareBibliographyDriver{incollection}{%
  \usebibmacro{bibindex}%
  \usebibmacro{begentry}%
  \usebibmacro{author/translator+others}%
  \setunit{\printdelim{nametitledelim}}\newblock
  \usebibmacro{title}%
  \newunit\newblock
  \usebibmacro{yjams:container}%
  \newunit
  \iffieldundef{maintitle}
  {\printfield{volume}%
  \printfield{part}}
  {}%
  \newunit\newblock
  \usebibmacro{series+number}%
  \newunit\newblock
  \usebibmacro{yjams:publisher+location}%
  \printfield{edition}%
  \newunit
  \usebibmacro{date}%
  \newunit
  \printfield{note}%
  \newunit\newblock
  \usebibmacro{chapter+pages}%
  \usebibmacro{yjams:language}%
  \newunit\newblock
\usebibmacro{yjams:finentry}}

\DeclareBibliographyDriver{inproceedings}{%
  \usebibmacro{bibindex}%
  \usebibmacro{begentry}%
  \usebibmacro{author/translator+others}%
  \setunit{\printdelim{nametitledelim}}\newblock
  \usebibmacro{title}%
  \newunit\newblock
  \usebibmacro{yjams:container}%
  \newunit
  \iffieldundef{maintitle}
  {\printfield{volume}%
  \printfield{part}}
  {}%
  \newunit\newblock
  \usebibmacro{series+number}%
  \newunit\newblock
  \printlist{organization}%
  \newunit
  \usebibmacro{publisher+location+date}%
  \newunit\newblock
  \printfield{note}%
  \newunit\newblock
  \usebibmacro{chapter+pages}%
  \usebibmacro{yjams:language}%
  \newunit\newblock
\usebibmacro{yjams:finentry}}

\DeclareBibliographyDriver{booklet}{%
  \usebibmacro{bibindex}%
  \usebibmacro{begentry}%
  \usebibmacro{author/editor+others/translator+others}%
  \setunit{\printdelim{nametitledelim}}\newblock
  \usebibmacro{title}%
  \newunit
  \printlist{language}%
  \newunit\newblock
  \usebibmacro{byauthor}%
  \newunit\newblock
  \usebibmacro{byeditor+others}%
  \newunit\newblock
  \printfield{howpublished}%
  \newunit\newblock
  \printfield{type}%
  \newunit\newblock
  \usebibmacro{location+date}%
  \newunit
  \printfield{note}%
  \newunit\newblock
  \usebibmacro{chapter+pages}%
  \newunit
  \printfield{pagetotal}%
  \newunit\newblock
  \iftoggle{bbx:isbn}
  {\printfield{isbn}}
  {}%
  \newunit\newblock
\usebibmacro{yjams:finentry}}

\bool_new:N \l__yjams_organization_heading_bool

\DeclareBibliographyDriver{manual}{%
  \usebibmacro{bibindex}%
  \usebibmacro{begentry}%
  \bool_set_false:N \l__yjams_organization_heading_bool
  \ifboolexpr{
    test {\ifnameundef{author}}
    and
    test {\ifnameundef{editor}}
    and
    not test {\iflistundef{organization}}
  }
  {\bool_set_true:N \l__yjams_organization_heading_bool
    \printlist{organization}%
    \setunit*{\addcomma\space}%
  \printlist{location}}
  {\usebibmacro{author/editor}}%
  \setunit{\printdelim{nametitledelim}}\newblock
  \usebibmacro{title}%
  \newunit
  \printlist{language}%
  \newunit\newblock
  \usebibmacro{byauthor}%
  \newunit\newblock
  \usebibmacro{byeditor}%
  \newunit\newblock
  \bool_if:NF \l__yjams_organization_heading_bool
  {\printlist{organization}%
    \newunit
    \printlist{location}%
  \newunit}%
  \printfield{edition}%
  \newunit\newblock
  \usebibmacro{series+number}%
  \newunit\newblock
  \printfield{type}%
  \newunit
  \printfield{version}%
  \newunit\newblock
  \printlist{publisher}%
  \newunit
  \usebibmacro{date}%
  \newunit
  \printfield{note}%
  \newunit\newblock
  \usebibmacro{chapter+pages}%
  \newunit
  \printfield{pagetotal}%
  \newunit\newblock
  \iftoggle{bbx:isbn}
  {\printfield{isbn}}
  {}%
  \newunit\newblock
\usebibmacro{yjams:finentry}}

\DeclareBibliographyDriver{misc}{%
  \usebibmacro{bibindex}%
  \usebibmacro{begentry}%
  \usebibmacro{author/editor+others/translator+others}%
  \setunit{\printdelim{nametitledelim}}\newblock
  \usebibmacro{title}%
  \newunit
  \printlist{language}%
  \newunit\newblock
  \usebibmacro{byauthor}%
  \newunit\newblock
  \usebibmacro{byeditor+others}%
  \newunit\newblock
  \printfield{howpublished}%
  \newunit\newblock
  \printfield{type}%
  \newunit
  \printfield{version}%
  \newunit\newblock
  \usebibmacro{organization+location+date}%
  \newunit
  \printfield{note}%
  \newunit\newblock
\usebibmacro{yjams:finentry}}

\DeclareBibliographyDriver{proceedings}{%
  \usebibmacro{bibindex}%
  \usebibmacro{begentry}%
  \bool_set_false:N \l__yjams_organization_heading_bool
  \ifnameundef{editor}
  {\bool_set_true:N \l__yjams_organization_heading_bool
  \printlist{organization}}
  {\usebibmacro{yjams:editor+others}}%
  \setunit{\printdelim{nametitledelim}}\newblock
  \usebibmacro{maintitle+title}%
  \newunit
  \printlist{language}%
  \newunit\newblock
  \usebibmacro{event+venue+date}%
  \newunit\newblock
  \usebibmacro{byeditor+others}%
  \newunit\newblock
  \iffieldundef{maintitle}
  {\printfield{volume}%
  \printfield{part}}
  {}%
  \newunit
  \printfield{volumes}%
  \newunit\newblock
  \usebibmacro{series+number}%
  \newunit\newblock
  \printlist{location}%
  \newunit
  \bool_if:NF \l__yjams_organization_heading_bool
  {\printlist{organization}%
  \newunit}%
  \printlist{publisher}%
  \newunit
  \usebibmacro{date}%
  \newunit
  \printfield{note}%
  \newunit\newblock
  \usebibmacro{chapter+pages}%
  \newunit
  \printfield{pagetotal}%
  \newunit\newblock
  \iftoggle{bbx:isbn}
  {\printfield{isbn}}
  {}%
  \newunit\newblock
\usebibmacro{yjams:finentry}}

\DeclareBibliographyDriver{report}{%
  \usebibmacro{bibindex}%
  \usebibmacro{begentry}%
  \usebibmacro{author}%
  \setunit{\printdelim{nametitledelim}}\newblock
  \usebibmacro{title}%
  \newunit
  \printlist{language}%
  \newunit\newblock
  \usebibmacro{byauthor}%
  \newunit\newblock
  \printfield{type}%
  \setunit*{\addspace}%
  \printfield{number}%
  \newunit\newblock
  \printfield{version}%
  \newunit\newblock
  \usebibmacro{institution+location+date}%
  \newunit
  \printfield{note}%
  \newunit\newblock
  \usebibmacro{chapter+pages}%
  \newunit
  \printfield{pagetotal}%
  \newunit\newblock
  \iftoggle{bbx:isbn}
  {\printfield{isrn}}
  {}%
  \newunit\newblock
\usebibmacro{yjams:finentry}}

\DeclareBibliographyDriver{thesis}{%
  \usebibmacro{bibindex}%
  \usebibmacro{begentry}%
  \usebibmacro{author}%
  \setunit{\printdelim{nametitledelim}}\newblock
  \usebibmacro{title}%
  \newunit
  \printlist{language}%
  \newunit\newblock
  \usebibmacro{byauthor}%
  \newunit\newblock
  \printfield{type}%
  \newunit
  \usebibmacro{institution+location+date}%
  \newunit
  \printfield{note}%
  \newunit\newblock
  \usebibmacro{chapter+pages}%
  \newunit
  \printfield{pagetotal}%
  \newunit\newblock
  \iftoggle{bbx:isbn}
  {\printfield{isbn}}
  {}%
  \newunit\newblock
\usebibmacro{yjams:finentry}}

\DeclareBibliographyAlias{conference}{inproceedings}
% --- End yj-amsalpha.sty ---

\addbibresource{references.bib}

% Separate compact footnote-sized back-reference page numbers with a period
% while retaining the hyperlinks supplied by BibLaTeX's pageref formatter.
\renewcommand*{\bibpagerefpunct}{\addperiod\space}
\renewbibmacro*{pageref}{%
  \iflistundef{pageref}
  {}
  {\printtext{%
      {\footnotesize
    \printlist[pageref][-\value{listtotal}]{pageref}}}%
\nopunct}}

% Force citation labels upright regardless of surrounding font shape
% (e.g. inside italic theorem bodies from \theoremstyle{plain}).
\renewcommand{\mkbibbrackets}[1]{\textup{[#1]}}
\renewcommand*{\yjamsalpha@bracketbreak}{]\addcomma\space[}

% Force uppercase initials after hyphens in given names.
% Biber preserves source case (e.g., Wu-yen -> W.-y.), but initials
% should be capitalized (W.-Y.).

% --- Chinese bibliography headings ---
\bool_if:NT \g__yjbib_chinese_bool
{
  \AddToHook{begindocument}
  {
    \defbibheading{bibliography}[\textbf{參考資料}]
    {
    \section*{#1} }
    \defbibheading{bibintoc}[\textbf{參考資料}]
    {
    \section*{#1} \addcontentsline{toc}{section}{#1} }
  }
}
% --- End yj-bib.sty ---

\ExecuteBibliographyOptions{mr=false}

% --- Begin yj-preamble.sty ---
% Materialized from yj-preamble.sty
% Classes using yj-bib must load it before this bundle because yj-refs
% transitively loads hyperref.
\PassOptionsToPackage{dvipsnames}{xcolor}

% --- Begin yj-graphics.sty ---
% Materialized from yj-graphics.sty
\usepackage{pgfplots}
\pgfplotsset{compat=1.17}
\usepackage{quiver}
\ExplSyntaxOff
\usetikzlibrary{babel,calc,positioning,arrows.meta,intersections}
\ExplSyntaxOn
% --- End yj-graphics.sty ---

% --- Begin yj-math.sty ---
% Materialized from yj-math.sty
\usepackage{amsmath}
\usepackage{mathtools}
% Use eqnlines only through its explicit equations environment.  Keep the
% established amsmath/mathtools environments, matrices, and \eqref unchanged.
\usepackage[
  defaults=classic,
  env=none,
  ang=false,
  eqref=false
]{eqnlines}
% eqnlines resets the inter-row penalty; apply our display-break policy last.
\allowdisplaybreaks
\usepackage{amssymb}
% Font ownership invariant: do not load math-font packages here.
% yj-fonts loads newtxtext before newtxmath so math operators match the text font.
\usepackage{stmaryrd}
% stmaryrd does not ship bold fonts; keep bold math from requesting U/stmry/b/n.
\SetSymbolFont{stmry}{bold}{U}{stmry}{m}{n}
\usepackage{dsfont}

% --- Operators ---

\DeclareMathOperator{\rmRe}{Re}
\DeclareMathOperator{\rmIm}{Im}

\clist_map_inline:nn
{
  supp, id, Id, rank, End, Aut, im, sgn, Hom, Ob, Pow, nullity, Gal, Irr,
  Coh, Per,
  Alt, Tr, GL, Sym, Spec, coker, SO, Gr, SU, Rad, Map, Res, Ind, Cl, FS,
  Bij, SL, Frac, Li, Ai, Log, ind, sech, csch, disc, Sp, arsinh, artanh,
  lcm, Ann, Ker, Sub, dist, Int,
  colim, vol, Coker, tr,
  dom, ran, cl, sign, Ad, DT
}
{ \exp_args:Nc \DeclareMathOperator {#1} {#1} }

% --- Delimiters ---
% newtxmath's extensible vertical bar can be taller than the surrounding
% delimiters at some automatically selected sizes.  For starred paired
% delimiters, use the font's native bar while it is large enough.  Only when
% TeX requests a larger delimiter do we draw a centered rule of that requested
% extent.  The rule is smashed so it does not affect the outer delimiter
% calculation.
\cs_new_nopar:Npn \__yj_math_middle: { \middle }
\cs_new_protected:Npn \__yj_middle_bar:nn #1#2
{
  \cs_if_eq:NNTF \delimsize \__yj_math_middle:
  { \mathpalette \__yj_auto_middle_bar:nn { {#1}{#2} } }
  { \mathrel { \delimsize\vert } }
}
\cs_new_protected:Npn \__yj_auto_middle_bar:nn #1#2
{ \__yj_auto_middle_bar_aux:nnn {#1} #2 }
\cs_new_protected:Npn \__yj_auto_middle_bar_aux:nnn #1#2#3
{
  \hbox_set:Nn \l_tmpa_box { $\m@th#1#2#3$ }
  \hbox_set:Nn \l_tmpb_box { $\m@th#1\vcenter{\hbox{}}$ }
  \dim_set:Nn \l_tmpa_dim
  { \box_ht:N \l_tmpa_box - \box_ht:N \l_tmpb_box }
  \dim_set:Nn \l_tmpb_dim
  { \box_dp:N \l_tmpa_box + \box_ht:N \l_tmpb_box }
  \dim_compare:nNnT { \l_tmpa_dim } < { \l_tmpb_dim }
  { \dim_set_eq:NN \l_tmpa_dim \l_tmpb_dim }
  \dim_set:Nn \l_tmpb_dim
  { \l_tmpa_dim * \delimiterfactor / 500 }
  \dim_set:Nn \l_tmpa_dim
  { 2 \l_tmpa_dim - \delimitershortfall }
  \dim_compare:nNnT { \l_tmpb_dim } < { \l_tmpa_dim }
  { \dim_set_eq:NN \l_tmpb_dim \l_tmpa_dim }
  \hbox_set:Nn \l_tmpa_box { $\m@th#1\vert$ }
  \dim_set:Nn \l_tmpa_dim
  { \box_ht:N \l_tmpa_box + \box_dp:N \l_tmpa_box }
  \dim_compare:nNnTF { \l_tmpb_dim } > { \l_tmpa_dim }
  {
    \mathrel
    {
      \smash
      {
        \vcenter
        {
          \hbox
          {
            \vrule
            width \fontdimen 8 \textfont 3
            height \l_tmpb_dim
          }
        }
      }
    }
  }
  { \mathrel { \vert } }
}
\DeclarePairedDelimiterX\setc[2]{\{}{\}}{
  \,#1 \__yj_middle_bar:nn {#1}{#2} #2\,
}
\DeclarePairedDelimiterX\vargs[2]{(}{)}{
  #1 \__yj_middle_bar:nn {#1}{#2} #2
}
\DeclarePairedDelimiter{\abs}{\lvert}{\rvert}

\DeclarePairedDelimiter{\angbr}{\langle}{\rangle}
\DeclarePairedDelimiterX{\angvargs}[2]{\langle}{\rangle}{
  #1 \__yj_middle_bar:nn {#1}{#2} #2
}

% The `c' in \cdbr means "composed": it builds [\![...\]\!] from
% two ordinary square brackets, unlike \dbr which uses \lBrack/\rBrack.
\NewDocumentCommand{\cdbr}{m}{\mathopen{[\![}#1\mathclose{]\!]}}

% --- Utilities ---
% Differential with math-atom spacing and an ordinary italic d.
\NewDocumentCommand{\dd}{}{\mathop{}\!d}

\usepackage{centernot}
\NewDocumentCommand{\ndiv}{O{}}{\mathrel{\centernot{#1\vert}}}

\usepackage{mleftright}

\cs_new_protected:Nn \__yj_math_cycle:nn
{
  \seq_set_split:Nnn \l_tmpa_seq { , } {#2}
  (
    \int_zero:N \l_tmpa_int
    \seq_map_inline:Nn \l_tmpa_seq
    {
      \int_incr:N \l_tmpa_int
      \int_compare:nNnF { \l_tmpa_int } = { 1 } {#1}
      ##1
    }
  )
}

\NewDocumentCommand{\cycle}{O{\;} m}
{ \__yj_math_cycle:nn {#1} {#2} }
% --- End yj-math.sty ---

% --- Begin yj-fonts.sty ---
% Materialized from yj-fonts.sty
% --- Text and math fonts ---
% Keep newtxtext before newtxmath. Loading newtxmath earlier makes math
% operators fall back to Latin Modern instead of the newtx/Termes text family.
\usepackage{newtxtext}
\defaultfontfeatures{}

\usepackage[varvw]{newtxmath}
\usepackage{anyfontsize}

% --- Symbol overrides (CM large operators) ---
\DeclareSymbolFont{largesymbolsCM}{OMX}{cmex}{m}{n}

\let\sum\relax
\DeclareMathSymbol{\sum}{\mathop}{largesymbolsCM}{"50}

\let\prod\relax
\DeclareMathSymbol{\prod}{\mathop}{largesymbolsCM}{"51}

\let\coprod\relax
\DeclareMathSymbol{\coprod}{\mathop}{largesymbolsCM}{"60}

\let\rtimes\relax
\DeclareMathSymbol{\rtimes}{\mathbin}{AMSb}{"6F}

% --- Math alphabets ---
\usepackage[cal=cm, bb=ams]{mathalfa}
\DeclareMathAlphabet{\eulerScr}{U}{eus}{m}{n}
% --- End yj-fonts.sty ---

% Standard table environments and column rules are document-type independent.
\usepackage{array}
\usepackage{booktabs}
\usepackage{longtable}
\usepackage{tabularx}
% Hyperref installs its float-anchor patch only when float is already loaded.
\usepackage{float}

% --- Begin yj-refs.sty ---
% Materialized from yj-refs.sty
\IfPackageLoadedF{hyperref}
{ \PassOptionsToPackage{pdfusetitle}{hyperref} }
\usepackage{xr-hyper}
\usepackage{hyperref}
\tl_const:Nn \c__yjrefs_shade_tl { 85 }
\hypersetup{
  linkcolor  = YellowOrange!\c__yjrefs_shade_tl!black,
  citecolor  = Aquamarine!\c__yjrefs_shade_tl!black,
  urlcolor   = violet!\c__yjrefs_shade_tl!black,
  colorlinks = true,
}

% Classes with nonstandard title/author storage register their metadata sources
% here. Explicit pdftitle/pdfauthor values remain authoritative.
\tl_new:N \g__yjrefs_pdf_title_source_tl
\tl_new:N \g__yjrefs_pdf_author_source_tl
\tl_new:N \l__yjrefs_pdf_title_tl
\tl_new:N \l__yjrefs_pdf_author_tl
\bool_new:N \g__yjrefs_custom_pdf_metadata_bool

\cs_new_protected:Npn \yjrefs_set_pdf_metadata_sources:nn #1#2
{
  \tl_gset:Nn \g__yjrefs_pdf_title_source_tl { #1 }
  \tl_gset:Nn \g__yjrefs_pdf_author_source_tl { #2 }
  \bool_gset_true:N \g__yjrefs_custom_pdf_metadata_bool
}

\pdfstringdefDisableCommands
{
  \def\\{~}
  \let\newline\\
  \def\and{;~}
  \def\thanks#1{}
  \def\orcid#1{}
  \let\cjk\@empty
}

\cs_new_protected:Npn \__yjrefs_apply_pdf_metadata:
{
  \bool_if:NT \g__yjrefs_custom_pdf_metadata_bool
  {
    % Materialize the registered source macros before \maketitle can clear
    % them. Bare math shifts are presentation-only and can be dropped from
    % otherwise plain metadata such as ``$K$-type''.
    \exp_args:NNV \tl_set:No
    \l__yjrefs_pdf_title_tl \g__yjrefs_pdf_title_source_tl
    \exp_args:NNV \tl_set:No
    \l__yjrefs_pdf_author_tl \g__yjrefs_pdf_author_source_tl
    \tl_remove_all:Nn \l__yjrefs_pdf_title_tl { $ }
    \tl_remove_all:Nn \l__yjrefs_pdf_author_tl { $ }
    \tl_if_empty:NT \@pdftitle
    {
      \hypersetup
      { pdftitle = { \tl_use:N \l__yjrefs_pdf_title_tl } }
    }
    \tl_if_empty:NT \@pdfauthor
    {
      \hypersetup
      { pdfauthor = { \tl_use:N \l__yjrefs_pdf_author_tl } }
    }
  }
}

\AddToHook{begindocument/before}{\__yjrefs_apply_pdf_metadata:}

\usepackage{zref-clever}
\zcsetup{cap, nameinlink=tsingle, noabbrev}
\usepackage{bookmark}
\usepackage{xurl}

% Funding identifiers retain their literal text and surrounding font while
% using the same automatic breakpoints as other identifier strings.
\DeclareUrlCommand\grantnumber{\urlstyle{same}}
% --- End yj-refs.sty ---

% --- Begin yj-theorems.sty ---
% Materialized from yj-theorems.sty
\bool_new:N \g__yjthm_chinese_bool
\bool_new:N \g__yjthm_noparent_bool

\keys_define:nn { yj-theorems }
{
  chinese  .bool_gset:N = \g__yjthm_chinese_bool,
  chinese  .default:n   = true,
  noparent .bool_gset:N = \g__yjthm_noparent_bool,
  noparent .default:n   = true,
}

\usepackage{zref-clever}
\usepackage{keytheorems}

\bool_if:NTF \g__yjthm_chinese_bool
{
  \theoremstyle{plain}
  \bool_if:NTF \g__yjthm_noparent_bool
  {
    \newkeytheorem{theorem}[
      name=定理,
      refname={定理,定理}, Refname={定理,定理}
    ]
  }
  {
    \newkeytheorem{theorem}[
      name=定理, parent=section,
      refname={定理,定理}, Refname={定理,定理}
    ]
  }
  \newkeytheorem{corollary}[
    name=推論, sibling=theorem,
    refname={推論,推論}, Refname={推論,推論}
  ]
  \newkeytheorem{conjecture}[
    name=猜想, sibling=theorem,
    refname={猜想,猜想}, Refname={猜想,猜想}
  ]
  \newkeytheorem{hypothesis}[
    name=假設, sibling=theorem,
    refname={假設,假設}, Refname={假設,假設}
  ]
  \newkeytheorem{lemma}[
    name=引理, sibling=theorem,
    refname={引理,引理}, Refname={引理,引理}
  ]
  \newkeytheorem{property}[
    name=性質, sibling=theorem,
    refname={性質,性質}, Refname={性質,性質}
  ]
  \newkeytheorem{proposition}[
    name=命題, sibling=theorem,
    refname={命題,命題}, Refname={命題,命題}
  ]
  \newkeytheorem{proposition*}[name=命題, numbered=false]
  \newkeytheorem{claim}[name=斷言, numbered=false]
  \newkeytheorem{fact}[
    name=事實, sibling=theorem,
    refname={事實,事實}, Refname={事實,事實}
  ]

  \theoremstyle{definition}
  \newkeytheorem{exercise}[
    name=習題, sibling=theorem,
    refname={習題,習題}, Refname={習題,習題}
  ]
  \newkeytheorem{exercise*}[name=習題, numbered=false]
  \newkeytheorem{example}[
    name=範例, sibling=theorem,
    refname={範例,範例}, Refname={範例,範例}
  ]
  \newkeytheorem{example*}[name=範例, numbered=false]
  \newkeytheorem{definition}[
    name=定義, sibling=theorem,
    refname={定義,定義}, Refname={定義,定義}
  ]
  \newkeytheorem{problem}[
    name=問題, sibling=theorem,
    refname={問題,問題}, Refname={問題,問題}
  ]

  \theoremstyle{remark}
  \newkeytheorem{remark}[
    name=備註, sibling=theorem,
    refname={備註,備註}, Refname={備註,備註}
  ]
  \newkeytheorem{remark*}[name=備註, numbered=false]
  \newkeytheorem{theorem*}[name=定理, numbered=false]
  \newkeytheorem{fact*}[name=事實, numbered=false]
  \newkeytheorem{definition*}[name=定義, numbered=false]
  \newkeytheorem{lemma*}[name=引理, numbered=false]
}
{
  \theoremstyle{plain}
  \bool_if:NTF \g__yjthm_noparent_bool
  { \newkeytheorem{theorem} }
  { \newkeytheorem{theorem}[parent=section] }
  \newkeytheorem{corollary}[sibling=theorem]
  \newkeytheorem{conjecture}[
    sibling=theorem,
    refname={conjecture,conjectures},
    Refname={Conjecture,Conjectures}
  ]
  \newkeytheorem{hypothesis}[
    sibling=theorem,
    refname={hypothesis,hypotheses},
    Refname={Hypothesis,Hypotheses}
  ]
  \newkeytheorem{lemma}[sibling=theorem]
  \newkeytheorem{property}[
    sibling=theorem,
    refname={property,properties},
    Refname={Property,Properties}
  ]
  \newkeytheorem{proposition}[sibling=theorem]
  \newkeytheorem{proposition*}[name=Proposition, numbered=false]
  \newkeytheorem{claim}[name=Claim, numbered=false]
  \newkeytheorem{fact}[
    sibling=theorem,
    refname={fact,facts},
    Refname={Fact,Facts}
  ]

  \theoremstyle{definition}
  \newkeytheorem{exercise}[
    sibling=theorem,
    refname={exercise,exercises},
    Refname={Exercise,Exercises}
  ]
  \newkeytheorem{exercise*}[name=Exercise, numbered=false]
  \newkeytheorem{example}[sibling=theorem]
  \newkeytheorem{example*}[name=Example, numbered=false]
  \newkeytheorem{definition}[sibling=theorem]
  \newkeytheorem{problem}[
    sibling=theorem,
    refname={problem,problems},
    Refname={Problem,Problems}
  ]

  \theoremstyle{remark}
  \newkeytheorem{remark}[sibling=theorem]
  \newkeytheorem{remark*}[name=Remark, numbered=false]
  \newkeytheorem{theorem*}[name=Theorem, numbered=false]
  \newkeytheorem{fact*}[name=Fact, numbered=false]
  \newkeytheorem{definition*}[name=Definition, numbered=false]
  \newkeytheorem{lemma*}[name=Lemma, numbered=false]
}

% --- Custom instance key: number ---
\keys_define:nn { keytheorems/thmuse }
{
  number .code:n =
  { \setcounter{theorem}{ \int_eval:n { #1 - 1 } } },
}

% --- Shared counters ---
\cs_set_eq:NN \c@equation \c@theorem
\cs_set_eq:NN \c@figure \c@theorem

% --- Chinese cross-reference type names (zref-clever) ---
\bool_if:NT \g__yjthm_chinese_bool
{
  % Section: "第 X 節" format — refbounds puts 節 after the number
  \zcRefTypeSetup{section}{
    Name-sg = {第}, name-sg = {第},
    Name-pl = {第}, name-pl = {第},
    refbounds            = {,,{},},
    +refbounds-rb        = {,,{},},
    +refbounds-re        = {,,{},},
    refbounds-first-sg   = {,,{節},},
    +refbounds-last      = {,,{節},},
    refbounds-last-re    = {,,{節},},
  }
  \zcRefTypeSetup{subsection}{
    Name-sg = 小節, name-sg = 小節,
    Name-pl = 小節, name-pl = 小節,
  }
  \zcRefTypeSetup{equation}{
    Name-sg = 式, name-sg = 式,
    Name-pl = 式, name-pl = 式,
  }
  \zcRefTypeSetup{figure}{
    Name-sg = 圖, name-sg = 圖,
    Name-pl = 圖, name-pl = 圖,
  }
  \zcRefTypeSetup{table}{
    Name-sg = 表, name-sg = 表,
    Name-pl = 表, name-pl = 表,
  }
}

% --- Chinese proof name ---
\bool_if:NT \g__yjthm_chinese_bool
{
  \AddToHook{begindocument}
  {
    \cs_if_exist:NT \proofname
    { \renewcommand*{\proofname}{\textbf{證明：}\nopunct} }
  }
}

% --- Solution environment ---
\bool_if:NTF \g__yjthm_chinese_bool
{
  \NewDocumentEnvironment{solution}{O{\textbf{解答：}\nopunct}}
  {
  \begin{proof}[#1] }
    {
  \end{proof} }
}
{
  \NewDocumentEnvironment{solution}{O{Solution}}
  {
  \begin{proof}[#1] }
    {
  \end{proof} }
}
% --- End yj-theorems.sty ---

% --- General formatting ---
\usepackage{comment}

\usepackage{enumitem}
\setlist[enumerate,1]{
  font=\normalfont,
  label=\textup{(\roman*)},
  listparindent=\parindent
}
% --- End yj-preamble.sty ---

% Keep body text and references fully justified while improving paragraph fit
% through glyph expansion and margin protrusion rather than loose global
% tolerance settings or ragged-right bibliography text.
\usepackage[protrusion=true]{microtype}

% Permit standard English hyphenation with two letters retained on either side.
% This adds genuine word-break opportunities without changing paragraph
% tolerance or warning thresholds.
\AddToHook{begindocument/end}{%
  \lefthyphenmin=2\relax
  \righthyphenmin=2\relax
}

\usepackage[margin=1in]{geometry}

% Retain amsart's zero natural paragraph separation while restoring the
% standard stretch needed to balance flush-bottom pages at paragraph boundaries.
\setlength{\parskip}{0pt plus 1pt}

% The page builder reserves the footnote insertion, not the separator rule.
% Keep that rule vertically neutral when the shared fonts change footnotesep,
% retaining amsart's 5pc width and 11pt reference baseline separation.
\cs_set:Npn \footnoterule
{
  \kern \dimexpr \footnotesep - 11.4pt \relax
  \hrule \@width 5pc \@height .4pt
  \kern \dimexpr 11pt - \footnotesep \relax
}

% BibLaTeX lists replace paragraph glue with their own item spacing. Retain
% compact entries while allowing flush-bottom bibliography pages to balance
% when amsart keeps the final author-address block together.
\setlength{\bibitemsep}{0pt plus .25\baselineskip}
% amsart's address block starts with a no-break request, but a list-ending
% breakpoint can already precede it. Keep the bibliography's final entry
% attached to the author addresses when they are present.
\AtBeginBibliography{
  \ifx\addresses\@empty\else
  \@endparpenalty\@M
  \fi
}

% Shared author blocks may include bare ORCID identifiers for journal wrappers.
% Plain amsart drafts do not print ORCID metadata.
\ProvideDocumentCommand{\orcid}{m}{}

% Hyperref's generic pdfusetitle support sees the optional short title and only
% the last amsart author. Use the full title and accumulated author list instead.
\yjrefs_set_pdf_metadata_sources:nn { \@title } { \authors }
\ExplSyntaxOff
\makeatother
\title[Critical L-values and genus-zero GW invariants]{Critical L-values and genus-zero Gromov--Witten invariants}

\author[Wu-yen Chuang]{Wu-yen Chuang\orcid{0000-0003-3230-3252}}
\address{Department of Mathematics and TIMS, National Taiwan University, Taipei, Taiwan}
\email{wychuang@gmail.com}

\author[Yi-Jing Tseng]{Yi-Jing Tseng\orcid{0009-0008-8349-6772}}
\address{Department of Mathematics, National Taiwan University, Taipei, Taiwan}
\email{yjtseng.math@gmail.com}


\begin{document}

\begin{abstract}
  We establish a period comparison for a two-dimensional rational Hodge summand of type $(3,0)+(0,3)$ in a smooth projective threefold over $\mathbb Q$ with $h^{3,0}=1$. Given an algebraic projector onto the summand and a correspondence whose action on the motive of a rational weight-four newform is nonzero after projection to the summand, restriction to the holomorphic de Rham lines identifies motivic morphisms with maps between these lines.
  %Chosen rational holomorphic generators therefore determine a unique isomorphism of Betti homological motives, whose marked Betti action expresses geometric periods in terms of two critical $L$-values.
  Fixing rational generators of the holomorphic de Rham lines determines a unique isomorphism of homological motives. Its Betti realization then yields an exact relation between geometric periods and two critical modular $L$-values.
  For the Fermat sextic, we construct the correspondence over $\mathbb Q$ and compute its Betti and de Rham actions, determining the scale from algebraic traces and cup-product pairings. Critical-value evaluations and, under mirror comparison and explicit continuation hypotheses, a continued genus-zero Gromov--Witten partition identity follow as corollaries.
\end{abstract}

\maketitle

\section{Introduction}

At a rank-two attractor, the holomorphic three-form spans the $(3,0)$ line of a two-dimensional rational Hodge summand. Numerical relations between Deligne periods and critical values at such points were studied in \cite{candelas-2020-742acd45}\cite[Sections~4--5]{yang-2021-6c28a647}. We study how an algebraic correspondence over $\mathbb Q$ between this summand and the motive of a weight-four newform compares periods of specified algebraic holomorphic forms in rational Betti markings. Its Betti matrix gives the coordinate comparison, and its de Rham action fixes the multiplicative constant.

Let $X/\mathbb Q$ be a smooth projective threefold with $h^{3,0}(X)=1$ whose extreme Hodge lines together are the complexification of a rational plane. An algebraic orthogonal projector $p_X$ onto that plane in degree three, zero elsewhere, defines the homological motive $M_X=(X,p_X,0)$ and determines a unique normalized rational newform $f$ of weight four and trivial nebentypus. If a rational-coefficient correspondence over $\mathbb Q$ from the fixed Scholl motive $M_f$ acts nontrivially after projection by $p_X$, it induces an isomorphism $M_f\simeq M_X$ with an explicit inverse. The Hodge conjecture supplies the projector in \zcref{intrinsic-projector-algebraicity}. Given the projector, the Tate conjecture supplies the comparison cycle in \zcref{conditional-derived-modular-comparison}.

Our main theorem, \zcref{intrinsic-absolute-comparison}, identifies the space of motivic morphisms with the space of maps between the holomorphic de Rham lines by restriction. Nonzero rational generators $\omega_f$ and $\Omega$ therefore select the unique isomorphism $\Phi_\Omega$ satisfying $(\Phi_\Omega)_{\mathrm{dR}}(\omega_f)=\Omega$. Its Betti realization identifies the corresponding rational Hodge structures, and naturality gives the period column \eqref{absolute-geometric-period-column} in any rational marking. Evaluating the inverse image of the geometric period class on the two fixed modular symbols gives the critical-value corollary \eqref{absolute-critical-values-from-periods}.

The weighted Fermat mirror is a complex rank-two attractor \cite[Section~8.3.1]{moore-1998-0cb358b1}. Its extreme $\Omega$-factor is associated with a weight-four form with complex multiplication of level $108$ \cite[Theorem~2(1), Section~8.1]{schimmrigk-2011-1f7180e5}. For this form, $f=\mathtt{108.4.a.c}$, \cite[Table~5]{bastian-2024-7db67c77} records $L(f,1)$ and $L(f,2)$. The weight-four functional equation relates its $L(f,1)$ entry to the untwisted $L(f,3)$ used here.

For the Fermat mirror quotient $X_0$, a common cover relates a Kummer model on the modular side to a degree-two quotient on the sextic side. Composing the graphs gives a rational-coefficient cycle over $\mathbb Q$ and an algebraic projector, so the general theorem applies unconditionally. The form-normalized isomorphism $\Phi_6$ of \eqref{sextic-form-normalized-cycle} sends $\omega_f$ to the residue form $\Omega_6$ fixed in \eqref{sextic-mum-residue}. In the basis dual to the two critical modular symbols and the large-radius period marking continued along the upper negative axis, its Betti matrix is $27C$, with $C$ as in \eqref{sextic-geometric-comparison-matrix}.

The integral in \zcref{sextic-euler-integral-period} identifies the sextic holomorphic line in this marking. By adjointness, the two algebraic traces show that the form-normalized comparison multiplies the normalized modular cup pairing by $6$. The source and target period columns determine its Betti matrix up to a positive scalar, and the cup pairings fix that scalar as $27$. The scale is therefore determined before the endpoint amplitude is evaluated. A beta integral then gives the two critical values in terms of the gamma function in \zcref{sextic-critical-gamma-evaluation}.

The continued Gromov--Witten partition identity is a consequence of the marked comparison. Partition and Bessel expansions for attractor coordinates were developed in \cite[Section~6, eqs.~(6.4)--(6.6)]{candelas-2021-56ab2e45}. The identities proposed in \cite[eqs.~(42)--(44), Section~5]{candelas-2025-ad265528} relate critical-value ratios to partition series of genus-zero invariants. Given a local mirror comparison in a rational period frame, the marked endpoint column expresses the mirror coordinate as a fractional-linear function of the two critical values. Lagrange inversion gives the local partition expansion. Analytic continuation along the inversion path uses the finite period limits and the square-root branch to identify its endpoint value with the mirror coordinate. For the sextic, a bound along the upper negative axis verifies the branch condition and gives \zcref{sextic-partition-comparison}.

\zcref{background} fixes the cohomological and period normalizations. \zcref{intrinsic-absolute-periods} establishes the period comparison and its critical-value and Gromov--Witten corollaries. \zcref{aesz8-sextic-geometric} constructs the algebraic correspondence for the Fermat sextic, proves the motivic comparison and calibrates its marked periods. \zcref{conclusion} summarizes its consequences. \zcref{sextic-analytic-appendix} proves the Taylor divergence and the square-root bound for the sextic.

\section{Preliminaries and normalizations}\label{background}

For a rational vector space $V$ and rational map $a$, write $V_{\mathbb C}=V\otimes_{\mathbb Q}\mathbb C$ and $a_{\mathbb C}=a\otimes_{\mathbb Q}1$. Coordinates on $V$ are specified by an isomorphism $\mathcal B\colon V\xrightarrow{\sim}\mathbb Q^n$. We write $\operatorname{CH}^r(Z)_{\mathbb Q}=\operatorname{CH}^r(Z)\otimes_{\mathbb Z}\mathbb Q$; rational coefficients do not specify the field of definition of a cycle.

\subsection{Hodge structures, realizations and correspondences}

For a rational pure Hodge structure of weight $n$, use $H_{B,\mathbb C}=\bigoplus_{a+b=n}H^{a,b}$, $F^rH_{B,\mathbb C}=\bigoplus_{a\ge r}H^{a,n-a}$ and scalar conjugation $\overline{H^{a,b}}=H^{b,a}$; morphisms preserve Hodge types. If $h^{3,0}=1$ and the rational subspace $H_B\cap(H^{3,0}\oplus H^{0,3})$ has dimension two, its complexification is $H^{3,0}\oplus H^{0,3}$. Our Tate convention is $\mathbb Q(m)$ of type $(-m,-m)$ and $H(m)^{a-m,b-m}=H^{a,b}$; a codimension-$r$ cycle class has type $(0,0)$ after twisting by $(r)$.

A specified scalar polarization $Q_B$ is $(-1)^n$-symmetric, pairs only complementary Hodge types, and satisfies $(-1)^{n(n-1)/2}i^{a-b}Q_{B,\mathbb C}(v,\overline v)>0$ for every nonzero $v\in H^{a,b}$. In weight three, Iritani uses minus the ordinary cup pairing on the residual cohomology \cite[Definition~6.5]{iritani-2011-7b0e4739}.

The restriction of $Q_B$ to a Hodge substructure $H'$ is nondegenerate: its radical is a Hodge substructure, while a nonzero component $v\in (H')^{a,b}$ in the radical would contradict this positivity. Thus the Tate-valued form identifies $(H')^\vee\simeq H'(n)$. In particular, for a weight-three structure of types $(3,0)+(0,3)$, one has $Q_{B,\mathbb C}(\omega,\overline\omega)\ne0$ whenever $0\ne\omega\in (H')^{3,0}$.

For a smooth projective threefold $X$, let $Q_X(x,y)=\int_{X(\mathbb C)}x\smile y$, using the complex orientation. This pairing is nondegenerate and alternating. Its identification with a polarization requires the Lefschetz decomposition and the Hodge--Riemann signs on each primitive summand. If $H^1(X)=0$, then $H^3(X)$ is primitive; otherwise the ample-class multiple of $H^1$ must also be included.

Fix a number field $k\subset\mathbb C$, an algebraic closure $\overline k\subset\mathbb C$, and a prime $\ell$. For smooth projective $X/k$, write $H_B^j(X)=H^j(X(\mathbb C),\mathbb Q)$, $H_{\mathrm{dR}}^j(X/k)=\mathbb H^j(X,\Omega^\bullet_{X/k})$ and $H_\ell^j(X)=H^j_{\mathrm{\acute{e}t}}(X_{\overline k},\mathbb Q_\ell)$. The last realization carries the action of $G_k=\operatorname{Gal}(\overline k/k)$. We use the comparison isomorphisms
\[
  I_X\colon H_{\mathrm{dR}}^j(X/k)\otimes_k\mathbb C
  \xrightarrow{\sim}H_B^j(X)\otimes_{\mathbb Q}\mathbb C,
  \qquad
  I_{B,\ell}\colon H_B^j(X)\otimes_{\mathbb Q}\mathbb Q_\ell
  \xrightarrow{\sim}H_\ell^j(X),
\]
where $I_X$ identifies the Hodge filtrations. Termwise GAGA for coherent sheaf cohomology on the projective variety $X$ \cite[Section~12, Theorem~1]{serre-1956-89fe9571}, the hypercohomology spectral sequence and the holomorphic Poincar\'e lemma give $I_X$, its filtration compatibility and its naturality for pullback, products and proper trace. Finite-coefficient Betti--\'etale comparison, inverse limits and scalar extension give $I_{B,\ell}$ \cite[Theorem~21.5]{milne-1998-6a5fd522}; this comparison also respects pullback, products, cycle classes and proper trace with the corresponding degree and Tate twists \cite[Sections~22--25]{milne-1998-6a5fd522}.

For smooth projective $S,X/k$, the map $H_\ell^j(S)\to H_\ell^j(X)$ obtained from a rational map $H_B^j(S)\to H_B^j(X)$ by Betti--\'etale comparison need not be $G_k$-equivariant.

Our cohomological convention uses $\mathbb Q_\ell(1)=(\varprojlim_n\mu_{\ell^n})\otimes_{\mathbb Z_\ell}\mathbb Q_\ell$ for the cyclotomic representation, where the transition maps on roots of unity are $\zeta\mapsto\zeta^\ell$. At a good rational prime $p\ne\ell$, geometric Frobenius acts on $\mathbb Q_\ell(1)$ by $p^{-1}$. In the convention in which $\mathbb Q_\ell(1)$ has Hodge--Tate weight $-1$, the rank-two weight-three factors below have weights $\{0,3\}$. For a two-dimensional $\mathbb Q_\ell$-representation $M_\ell$ of $G_{\mathbb Q}$, a perfect alternating equivariant pairing $M_\ell\times M_\ell\to\mathbb Q_\ell(-3)$ gives $\det M_\ell\simeq\mathbb Q_\ell(-3)$ and $\det(\operatorname{Frob}_p\mid M_\ell)=p^3$. An arithmetic-Frobenius or dual-representation convention reverses the relevant eigenvalues; every comparison in this paper uses the displayed cohomological convention.

For smooth projective $Z/k$, the cycle-class map $\operatorname{CH}^r(Z)\to H^{2r}_{\mathrm{\acute{e}t}}(Z_{\overline k},\mathbb Q_\ell(r))$ factors through rational equivalence, and a cycle defined over $k$ has $G_k$-invariant class \cite[Section~23]{milne-1998-6a5fd522}.

\begin{definition}\label{degree-three-cycle-action}
  Let $k\subset\mathbb C$ be a number field and let $S,X$ be smooth projective threefolds over $k$. For $\Gamma\in\operatorname{CH}^3(S\times_kX)\otimes_{\mathbb Z}\mathbb Q$, let $\operatorname{cl}_B(\Gamma)\in H_B^6(S\times_kX)(3)$ be its Betti cycle class. Define $\Gamma_{*,B}$ to be the composite
  \[
    H_B^3(S)\xrightarrow{\operatorname{pr}_S^*}H_B^3(S\times_kX)
    \xrightarrow{\smile\operatorname{cl}_B(\Gamma)}H_B^9(S\times_kX)(3)
    \xrightarrow{(\operatorname{pr}_X)_*}H_B^3(X).
  \]
  Thus $\Gamma_{*,B}(x)=(\operatorname{pr}_X)_*(\operatorname{pr}_S^*x\smile\operatorname{cl}_B(\Gamma))$. The pushforward lowers degree by six and removes the twist $(3)$, since $\operatorname{pr}_X$ has relative dimension three. Define $\Gamma_{*,\mathrm{dR}}$ and $\Gamma_{*,\ell}$ by the same graded operations in their realizations; the latter has type $H_\ell^3(S)\to H_\ell^3(X)$.
\end{definition}
If $\Gamma$ is defined over $k$, its twisted Betti class has type $(0,0)$, so $\Gamma_{*,B}$ is a rational Hodge morphism. Functoriality of comparison gives
\[
  I_X\circ(\Gamma_{*,\mathrm{dR}}\otimes1)
  =(\Gamma_{*,B}\otimes1)\circ I_S,
  \qquad
  I_{B,\ell,X}\circ(\Gamma_{*,B}\otimes1)
  =\Gamma_{*,\ell}\circ I_{B,\ell,S},
\]
and $g\Gamma_{*,\ell}=\Gamma_{*,\ell}g$ for every $g\in G_k$. These assertions follow by applying the comparison maps to pullback, multiplication by the cycle class and the twisted trace in \zcref{degree-three-cycle-action} \cite[Sections~23--25]{milne-1998-6a5fd522}. The induced cohomological action depends only on the cycle's cohomology class.

For a rational-coefficient correspondence $e$ on $S$ with $e_B^2=e_B$, the space $e_BH_B^3(S)$ is a rational Hodge substructure. The algebraic cycle class gives corresponding projectors in every comparison realization.

We work in the idempotent-complete category of homological motives over $\mathbb Q$ with rational coefficients, with correspondences modulo Betti homological equivalence. We write an effective summand as $(Z,p,0)$ when the total Betti action of the codimension-$\dim Z$ correspondence $p$ is idempotent. Its identity is the homological class of $p$, and $\operatorname{Hom}(M_1,M_2)$ for motives is taken in this category. Total Betti realization is faithful by the defining equivalence relation. In particular, the construction of an object in this category requires an all-degree idempotent action, but no equality of cycle representatives in the Chow group.

For smooth projective threefolds, duality identifies $\operatorname{Hom}(H_\ell^3(S),H_\ell^3(X))$ with the K\"unneth summand $(H_\ell^3(S)\otimes H_\ell^3(X))(3)\subseteq H_\ell^6(S\times_kX)(3)$. Representing an equivariant map by a cycle is an additional algebraicity statement.

\begin{conjecture}[note={Hodge conjecture \cite[Conjecture~6]{charles-2014-8683f1f1}}]\label{rational-hodge-conjecture}
  Let $Z$ be a smooth projective complex variety and let $r\geq0$. The rational Hodge classes of type $(0,0)$ in $H^{2r}(Z,\mathbb Q(r))$ are generated over $\mathbb Q$ by the cohomology classes of codimension-$r$ algebraic cycles on $Z$.
\end{conjecture}

\begin{conjecture}[note={Tate conjecture \cite[Conjecture~1.1]{totaro-2017-4fdc1de0}}]\label{global-rational-tate-conjecture}
  Let $k$ be a field finitely generated over $\mathbb Q$, let $Z/k$ be smooth and projective, let $r\geq0$, and let $\ell$ be any prime. The cycle-class map $\operatorname{CH}^r(Z)\otimes_{\mathbb Z}\mathbb Q_\ell \longrightarrow H^{2r}_{\mathrm{\acute{e}t}}(Z_{\overline k},\mathbb Q_\ell(r))^{G_k}$ is surjective, with cycles defined over $k$.
\end{conjecture}

The Hodge conjecture supplies the algebraic projector in \zcref{intrinsic-projector-algebraicity}, and the Tate conjecture supplies a nonzero correspondence from $M_f$ after applying $e_f$ in \zcref{conditional-derived-modular-comparison}. Subsequent comparison results take the existence of these objects as hypotheses; for the Fermat sextic they are constructed unconditionally in \zcref{sextic-all-degree-projector,sextic-algebraic-comparison}.

\subsection{Modular forms and critical periods}

For $\tau$ in the upper half-plane $\mathfrak H$, put $q=e^{2\pi i\tau}$. We write $S_4^{\mathrm{new}}(\Gamma_0(N),\varepsilon)$ for the cuspidal weight-four newspace with nebentypus $\varepsilon$. A normalized rational newform is a simultaneous Hecke eigenform $f=\sum_{n\ge1}a_n(f)q^n$ with $a_1(f)=1$ and $a_n(f)\in\mathbb Q$. The modular-symbol result below allows compatible nebentypus; the geometric factor constructed later has trivial nebentypus.

For a normalized rational cuspidal weight-four newform $f$, define initially for $\rmRe s>3$
\[
  L(f,s)=\sum_{n\ge1}a_n(f)n^{-s},
  \qquad
  L(f,\chi_{-4},s)=\sum_{n\ge1}\chi_{-4}(n)a_n(f)n^{-s},
\]
where $\chi_{-4}(2m)=0$ and $\chi_{-4}(2m+1)=(-1)^m$; elsewhere these denote analytic continuations. The second series is the fixed coefficient twist, including the omitted Euler factor at two. Its identification with a primitive newform twist would require the bad-prime factors.

Fix a sufficiently fine modular level and a smooth projective Kuga--Sato threefold $S_f/\mathbb Q$, obtained by resolving the two-fold fibre product of the compactified universal generalized elliptic curve. Scholl's finite-group projector selects the symmetric-square parabolic cohomology and annihilates every degree other than three. The Hecke factor and descent give a rational algebraic correspondence $e_f$ whose total Betti action is the projector onto the $f$-summand \cite[Section~1.2.0, Theorems~1.2.1, 1.2.4, Remark~1.2.6]{scholl-1990-bb6eb16d}. Fix this presentation and put $M_f=(S_f,e_f,0)$, with
\begin{align*}
  M_{f,B}&=(e_f)_{*,B}H_B^3(S_f),\\
  M_{f,\mathrm{dR}}&=(e_f)_{*,\mathrm{dR}}H_{\mathrm{dR}}^3(S_f/\mathbb Q),\\
  M_{f,\ell}&=(e_f)_{*,\ell}H_\ell^3(S_f).
\end{align*}
These are rank-two realizations; $M_{f,B}$ has weight three and types $(3,0)+(0,3)$. Restricting $I_{S_f}$ gives
\[
  I_f\colon M_{f,\mathrm{dR}}\otimes_{\mathbb Q}\mathbb C\xrightarrow{\sim}M_{f,B}\otimes_{\mathbb Q}\mathbb C.
\]
Its holomorphic de Rham line is $\ell_f=F^3M_{f,\mathrm{dR}}$, represented on the analytic uniformization by $f(\tau)\dd\tau\wedge\dd z_1\wedge\dd z_2$ \cite[Section~4.2.0]{scholl-1990-bb6eb16d}.

Use the cuspidal new Hecke quotient of the degree-two modular symbols, with the path and modular relations and $T_n=a_n(f)$. Eichler--Shimura and the Scholl realization identify it with $M_{f,B}^\vee$ \cite[Theorem~3, Proposition~14]{merel-1994-bc469ad4}, \cite[Section~1.2.0, Theorem~1.2.1]{scholl-1990-bb6eb16d}. For a homogeneous quadratic $P\in\mathbb Q[X,Y]$ and $\alpha,\beta\in\mathbb P^1(\mathbb Q)$, write $\Phi_f(P\{\alpha,\beta\})=\int_\alpha^\beta f(z)P(z,1)\dd z$ in this quotient.

For trivial nebentypus, $M_{f,\ell}$ has determinant $\mathbb Q_\ell(-3)$ and good-prime polynomial $1-a_p(f)T+p^3T^2$. In this case fix the rational generator $\omega_f$ of $\ell_f$ by its Tate expression
\[
  \omega_f=f(q)\frac{\dd q}{q}\wedge\frac{\dd u_1}{u_1}\wedge\frac{\dd u_2}{u_2}.
\]
Choose the fibre cycles oriented from $0$ to $1$ and from $0$ to $\tau$ in $\mathbb C/(\mathbb Z+\mathbb Z\tau)$, with the ordered product orientation on the fibre square. Put $q=e^{2\pi i\tau}$ and $u_j=e^{2\pi i\zeta_j}$. Each logarithmic differential contributes $2\pi i$, so the modular-symbol pairing is $s_{\mathbb C}(I_f\omega_f)=(2\pi i)^3\Phi_f(s)$.

Deligne's conjecture relates critical motivic special values to determinant periods defined from Betti--de Rham comparison \cite[Section~1.7, Conjecture~1.8]{deligne-1979-61965bb8}. The modular-symbol calculation identifies the holomorphic modular line directly. Its comparison with a geometric motive is made in the next section.

\begin{proposition}\label{critical-coordinate-basis}\label{uniform-critical-column}
  Let $f\in S_4^{\mathrm{new}}(\Gamma_0(N),\varepsilon)$ be normalized with rational Fourier coefficients and $\dim_{\mathbb Q}M_{f,B}=2$. In its cuspidal new quotient let $s_0,s_\chi$ be the images of $X^2\{0,\infty\}$ and $(X-Y/4)^2\{1/4,\infty\}-(X-3Y/4)^2\{3/4,\infty\}$. Then $(s_0,s_\chi)$ is a rational basis of $M_{f,B}^\vee$. Both critical values are nonzero and $(L(f,3),2iL(f,\chi_{-4},3))^T$ spans the holomorphic line in the dual coordinates. For trivial nebentypus the period column is
  \begin{equation}\label{algebraic-modular-period-column}
    \begin{pmatrix}s_0(I_f\omega_f)\\s_\chi(I_f\omega_f)
    \end{pmatrix}
    =-2
    \begin{pmatrix}L(f,3)\\2iL(f,\chi_{-4},3)
    \end{pmatrix}.
  \end{equation}
\end{proposition}

\begin{proof}
  The character version of the modular-symbol pairing identifies integration of the indicated symbols with evaluation on the holomorphic modular line. Mellin integration of $\tau^2 f(\tau)$ gives $\Phi_f(s_0)=-2(2\pi i)^{-3}L(f,3)$. For $s_\chi$, shifting both paths to the imaginary axis replaces $f(z)$ by $f(z+1/4)-f(z+3/4)$. This is $2i\sum_{n\ge1}a_n\chi_{-4}(n)e^{2\pi inz}$, so $\Phi_f(s_\chi)=-4i(2\pi i)^{-3}L(f,\chi_{-4},3)$. The cusp integrals converge and the Mellin identities follow by continuation from a right half-plane.

  At $s=3$ both Euler products converge absolutely and are nonzero. At good primes the roots have absolute value $r^{3/2}$ by \cite[Theorem~8.2]{deligne-1974-370b49c1}. At bad primes $a_{r^j}=a_r^j$ and $a_n=O_f(n^2)$ give $\abs*{a_r}\le r^2<r^3$. The coefficient twist omits the factor at two and multiplies the odd-prime roots by signs. The real Fourier coefficients make the first pairing purely imaginary and the second real, so the two nonzero pairings prove rational independence. The displayed column follows from the Tate factor.
\end{proof}

These coordinates remain available when $L(f,2)=0$.

\subsection{Frobenius and genus-zero period conventions}\label{frobenius-prepotential-coordinates}

For a smooth projective family $\pi\colon\mathcal X\to S$ of threefolds over a complex curve, use the rational local system $R^3\pi_*\mathbb Q$, its bundle $\mathcal H=(R^3\pi_*\mathbb Q)\otimes\mathcal O_S$, flat alternating pairing to $\mathbb Q(-3)$, and Hodge filtration with $\nabla(F^r\mathcal H)\subseteq F^{r-1}\mathcal H\otimes\Omega_S^1$. A selected rank-four subvariation carries the induced filtration and restricted polarization. For $\Omega\in F^3\mathcal H$, its coordinates are $\int_{\gamma_j}\Omega$ in a rational locally flat cycle basis, continued with that basis. A scalar solution basis alone does not specify this rational marking.

For a rank-four variation on a punctured disc, let $T_0$ be the monodromy of a positive loop about $z=0$ and $N=\log T_0$. \emph{Maximal unipotent monodromy} (MUM) means $(T_0-I)^4=0$, $N^4=0$ and $N^3\ne0$.

A coordinate, a branch of $\log z$, a holomorphic form and a rational basis still have to be fixed. With $\theta=z\dd/\dd z$, cyclicity of the chosen form means that $\Omega,\nabla_\theta\Omega,\nabla_\theta^2\Omega,\nabla_\theta^3\Omega$ are a basis of the selected differential module. Its next derivative gives a fourth-order scalar Picard--Fuchs equation. A general fourth-order operator need not come from such geometric data.

Let $\mathcal L=\sum_{j=0}^sz^j\mathcal L_j(\theta)$, where $\mathcal L_j\in\mathbb Q[\theta]$, $\deg \mathcal L_j\leq4$, $\mathcal L_0(\theta)=\theta^4$, and $\theta=z\dd/\dd z$. Set $\alpha_0(\rho)=1$ and define the higher Frobenius coefficients by
\begin{equation}\label{constructed-frobenius-recurrence}
  \alpha_n(\rho)=-\frac{1}{(n+\rho)^4}\sum_{j=1}^{\min(s,n)}\mathcal L_j(n-j+\rho)\alpha_{n-j}(\rho),\qquad n\geq1.
\end{equation}
On a fixed logarithm branch, let $[\rho^r]$ extract the coefficient of $\rho^r$ and put
\begin{equation}\label{frobenius-solutions}
  w_r(z)=[\rho^r]\mleft(z^\rho\sum_{n\geq0}\alpha_n(\rho)z^n\mright),\qquad0\leq r\leq3.
\end{equation}
Then $w_0(0)=1$ and the leading logarithm of $w_r$ is $(\log z)^r/r!$. We use the regular coefficient series
\[
  h_r(z)=\sum_{n\geq0}[\rho^r]\alpha_n(\rho)z^n,\qquad
  w_r(z)=\sum_{j=0}^r\frac{(\log z)^j}{j!}h_{r-j}(z).
\]
Thus $h_0(0)=1$ and $h_r(0)=0$ for $r>0$. Put
\begin{equation}\label{mirror-quantities}
  \begin{aligned}
    t&=\frac{w_1}{2\pi i w_0},& q&=e^{2\pi it}=z\exp(h_1/h_0),\\
    V(z)&=\frac{w_2}{w_0}-\frac12\left(\frac{w_1}{w_0}\right)^2
    =\frac{h_2}{h_0}-\frac12\left(\frac{h_1}{h_0}\right)^2.
  \end{aligned}
\end{equation}
The quantities in \eqref{mirror-quantities} come from the B-model period equation. The A-model invariants below define the prepotential. The local mirror comparison \eqref{local-full-mirror-frame} identifies its column with the normalized periods in a specified rational flat frame.

We use $Y$ for the A-model and $X$ for the B-model or an endpoint. Here a \emph{Calabi--Yau threefold} $Y$ is smooth and projective, with $K_Y\simeq\mathcal O_Y$ and $H^1(Y,\mathcal O_Y)=H^2(Y,\mathcal O_Y)=0$. For such $Y$, the moduli stack $\overline{\mathcal M}_{0,0}(Y,\beta)$ of genus-zero stable maps in an effective class $\beta$ carries a virtual fundamental class of virtual dimension zero \cite[Section~2.1]{iritani-2011-7b0e4739}. Its degree is a rational invariant, with automorphism weights included. Degree grouping with an ample class and diagonal specialization give the conventions used throughout the paper.

For a primitive ample class $H$, use the classical coefficients of \cite[eqs.~(5)--(6)]{candelas-2025-ad265528}:
\begin{equation}
  Y_{111}=\int_YH^3>0,\qquad
  b=\frac1{24}\int_Yc_2(Y)H,\qquad
  c_0=\frac{\zeta(3)}{2(2\pi i)^3}\int_Yc_3(Y).
\end{equation}
Here $\int_Yc_3(Y)=\chi(Y)$, and the rational coefficient $Y_{011}$ depends on the chosen frame.

Set $p_j=w_j/(2\pi i)^j$. The bases in \cite[eqs.~(2)--(3)]{candelas-2025-ad265528} are $(p_0,p_1,Y_{111}p_2,Y_{111}p_3)$ and $(w_0,w_1,Y_{111}w_2,Y_{111}w_3)$. The unnormalized Frobenius column in the prepotential frame is $\Pi_{\mathrm{raw}}=(\mathcal F_0,\mathcal F_1,z^0,z^1)^T$ \cite[eqs.~(12)--(14)]{candelas-2025-ad265528}, with
\begin{equation}\label{constructed-prepotential-periods}
  \begin{aligned}
    z^0&=p_0,& \qquad z^1&=p_1,\\
    \mathcal F_1&=bp_0-Y_{011}p_1-Y_{111}p_2,&
    \mathcal F_0&=2c_0p_0+bp_1+Y_{111}p_3.
  \end{aligned}
\end{equation}

\begin{definition}\label{genus-zero-period-convention}
  Fix a geometric Calabi--Yau threefold A-model $Y$ and a primitive ample class $H$. Sum the connected, unmarked genus-zero Gromov--Witten invariants over curve classes of $H$-degree $d$ to obtain
  \begin{equation*}
    N_d^{\mathrm{GW}}=\sum_{\int_\beta H=d}
    \deg[\overline{\mathcal M}_{0,0}(Y,\beta)]^{\mathrm{vir}},\qquad d>0,
  \end{equation*}
  where the sum ranges over effective curve classes $\beta\in H_2(Y,\mathbb Z)$. When all degree variables are set equal, sum the multidegree invariants of total degree $d$. Fix a rational prepotential frame with classical coefficients $Y_{111},Y_{011},b,c_0$. In the flat coordinate $t$, define
  $q=e^{2\pi it}$ and $\mathcal I(t)=(2\pi i)^{-3}\sum_{d\geq1}N_d^{\mathrm{GW}}q^d$, initially as a formal series in $q$. Define the prepotential by
  \[
    F^{(0)}(t)=-\frac{Y_{111}}6t^3-\frac{Y_{011}}2t^2+bt+c_0-\mathcal I(t).
  \]
  In the coordinate order of \cite[eq.~(8)]{candelas-2025-ad265528}, put
  \begin{equation}
    \Pi(t)=
    \begin{pmatrix}
      2F^{(0)}(t)-t\partial_tF^{(0)}(t)\\
      \partial_tF^{(0)}(t)\\
      1\\t
    \end{pmatrix}.
  \end{equation}
\end{definition}

$N_d^{\mathrm{GW}}$ are virtual degrees, with no integrality assumption. Formal differentiation in $t$, using $q=e^{2\pi it}$, gives $\partial_t\mathcal I(t)=(2\pi i)^{-2}\sum_{d\ge1}dN_d^{\mathrm{GW}}q^d$. The same equality holds analytically wherever the series converges, and subsequently for their analytic continuations on a common branch.

In a geometric realization of the prepotential, the raw periods satisfy $t=z^1/z^0$, $\mathcal F_1/z^0=\partial_t F^{(0)}$ and $\mathcal F_0/z^0=2F^{(0)}-t\partial_t F^{(0)}$. The normalized column in the definition is $\Pi=\Pi_{\mathrm{raw}}/z^0$, defined only in the chart $z^0\ne0$. Multiplying all raw periods by one nonzero function does not change these ratios. Changing a Betti frame, by contrast, changes the coordinate functions and must be tracked separately.

Mirror comparison identifies the quantum connection of $Y$ with $(\mathcal H,\nabla)$, with the stated periods and normalization. Its geometric scope is specified below.

The class $\widehat\Gamma_Y$ connects this comparison to rational cycles. For ordinary Chern roots $\delta_j$ of $T_Y$ with $c_1(Y)=0$, the convention used below is
\begin{equation}\label{constructed-gamma-class}
  \widehat\Gamma_Y=\prod_{j=1}^{3}\Gamma\mleft(1+\frac{\delta_j}{2\pi i}\mright)
  =1+\frac{c_2(Y)}{24}-\frac{\zeta(3)c_3(Y)}{(2\pi i)^3}.
\end{equation}
Here $\sum_j\delta_j^2=-2c_2(Y)$, $\sum_j\delta_j^3=3c_3(Y)$ and $-\zeta(2)/(2\pi i)^2=1/24$; terms of cohomological degree greater than six vanish. This is the degree-rescaled version of the Gamma class in \cite[Definition~3.6]{iritani-2011-7b0e4739}: its factors $(2\pi i)^{\deg/2}$ on the Chern character and the overall factor $(2\pi i)^{-3}$ give the same top-degree integral as \eqref{constructed-gamma-class} with the ordinary Chern character. The divisor coordinate in that convention is $(2\pi i)tH$.

\subsection{Continuation and endpoint periods}

To continue a local holomorphic function, fix the path $\gamma\colon[0,1)\to S$, logarithm branch and scalar normalization. Its endpoint value is the limit on that branch. A formal identity in $\mathbb C\cdbr{q}$ becomes an equality near $q=0$ after convergence is proved. Equality of continued endpoint limits requires continuation along the chosen path, independently of convergence of the original $q$-series there.

Fix the raw period coordinates of a local holomorphic form. At a singular endpoint, fix a finite base change and a punctured sector in its local coordinate $u$. A nonvanishing holomorphic function $a(u)$ on that sector acts by $\Omega\longmapsto a(u)\Omega$ and $(\mathcal F_0,\mathcal F_1,z^0,z^1)\longmapsto a(u)(\mathcal F_0,\mathcal F_1,z^0,z^1)$.

The gauge may depend on $u$. Its effect on the scalar differential equation must then be included; a nonconstant gauge is not a constant change of solution basis. The common factors cancel in normalized ratios, whereas the limit of $a(u)z^0$ depends on the prescribed gauge. At a finite-monodromy point, a finite base change can remove monodromy; at a unipotent degeneration the limiting object is generally mixed, so an ordinary pure Hodge comparison requires an identified pure realization rather than an unqualified substitution into smooth-fibre cohomology.

A square-root branch is determined by its value at a point where the radicand is nonzero and the continuation path.

For the raw periods of \zcref{genus-zero-period-convention}, the normalized column $\Pi=\Pi_{\mathrm{raw}}/z^0$ has third coordinate one and fourth coordinate $t$. Along a fixed continuation path on which $z^0\ne0$, a finite limit $\Pi_*$ therefore gives $t_*=\Pi_{*,4}$.

\section{Period comparisons}\label{intrinsic-absolute-periods}

\subsection{Algebraic projectors and modular realizations}

For a smooth projective threefold $X/\mathbb Q$ with $h^{3,0}(X)=1$, put
\[
  H_X^{\mathrm{att}}=H_B^3(X,\mathbb Q)\cap(H^{3,0}(X)\oplus H^{0,3}(X)),\qquad
  \ell_X=F^3H_{\mathrm{dR}}^3(X/\mathbb Q).
\]
If $\dim_{\mathbb Q}H_X^{\mathrm{att}}=2$, comparison identifies $\ell_X\otimes_{\mathbb Q}\mathbb C$ with $H^{3,0}(X)\subset H_{X,\mathbb C}^{\mathrm{att}}$. A correspondence $e$ acting as the cup-product orthogonal projector onto $H_X^{\mathrm{att}}$ in degree three and acting by zero in every other degree defines $M_X=(X,e,0)$.

\begin{theorem}\label{intrinsic-projector-algebraicity}
  Assume the Hodge conjecture of \zcref{rational-hodge-conjecture}. Let $X/\mathbb Q$ be a smooth projective threefold, with $h^{3,0}(X)=1$ and $\dim_{\mathbb Q}H_X^{\mathrm{att}}=2$. Let $e_B\in\operatorname{End}_{\mathbb Q}(H_B^*(X,\mathbb Q))$ be the $Q_X$-orthogonal projector onto $H_X^{\mathrm{att}}$ in degree three, extended by zero in every other degree: $e_B|_{H_X^{\mathrm{att}}}=I_{H_X^{\mathrm{att}}}$ and $e_B|_{(H_X^{\mathrm{att}})^\perp}=0$. There is $\Gamma\in\operatorname{CH}^3(X\times_{\mathbb Q}X)_{\mathbb Q}$ with $\Gamma_{*,B}=e_B$.
\end{theorem}

\begin{proof}
  The plane $H_{X,\mathbb C}^{\mathrm{att}}=H^{3,0}(X)\oplus H^{0,3}(X)$ is primitive. The Hodge--Riemann relations make its cup pairing nondegenerate, giving
  \[
    H_B^3(X,\mathbb Q)
    =H_X^{\mathrm{att}}\overset{\perp,Q_X}{\oplus}(H_X^{\mathrm{att}})^\perp.
  \]
  Thus $e_B$ is a rational Hodge endomorphism. Poincar\'e duality and K\"unneth identify its total action with a rational Hodge class in $H_B^6(X\times X,\mathbb Q(3))$, represented by an algebraic cycle under the Hodge conjecture.

  The Hilbert-scheme argument of \cite[Remark~46]{charles-2014-8683f1f1}, applied componentwise, replaces this cycle by an algebraically equivalent cycle over $\overline{\mathbb Q}$. Choose a finite Galois extension $K/\mathbb Q$ over which a representative $\Gamma_K$ is defined. Algebraic equivalence preserves the cycle class, so Betti--de Rham comparison shows that its de Rham action has rank two in degree three, is the identity on $F^3H^3_{\mathrm{dR}}(X/\mathbb Q)\otimes_{\mathbb Q}K$, and vanishes in every other degree. These properties persist under $\operatorname{Gal}(K/\mathbb Q)$, since the cohomology spaces and filtration are defined over $\mathbb Q$. For each $\sigma$, the Betti action of $\Gamma_K^\sigma$ is a rational Hodge morphism fixing $H^{3,0}(X)$ and hence, by scalar conjugation, $H^{0,3}(X)$. These lines exhaust its rank, so preservation of Hodge types makes it zero on $H^{2,1}(X)\oplus H^{1,2}(X)$. Every conjugate therefore has total action $e_B$, and $\Gamma=[K:\mathbb Q]^{-1}\operatorname{cor}_{K/\mathbb Q}(\Gamma_K)$ has the required action.
\end{proof}

An already constructed nonzero modular correspondence supplies the projector without an additional algebraicity hypothesis.

\begin{lemma}\label{correspondence-induced-projector}
  Let $S,X/\mathbb Q$ be smooth projective threefolds with $h^{3,0}(X)=1$. Let $e\in\operatorname{CH}^3(S\times_{\mathbb Q}S)_{\mathbb Q}$ act on total Betti cohomology as a projector whose image is a two-dimensional Hodge structure $M\subset H_B^3(S,\mathbb Q)$ of types $(3,0)+(0,3)$. Let $\Gamma\in\operatorname{CH}^3(S\times_{\mathbb Q}X)_{\mathbb Q}$ act nontrivially on $M^{3,0}$. Put $Z_\Gamma=\Gamma\circ e\circ{}^te\circ{}^t\Gamma$, and let $\Delta_X\subset X\times_{\mathbb Q}X$ be the diagonal. Then $\dim_{\mathbb Q}H_X^{\mathrm{att}}=2$, and
  \begin{equation}
    d_\Gamma=-\frac12\deg\mleft(Z_\Gamma\cdot[\Delta_X]\mright)\in\mathbb Q_{>0},
    \qquad
    p_X=d_\Gamma^{-1}Z_\Gamma
  \end{equation}
  is a self-adjoint algebraic correspondence over $\mathbb Q$ whose total cohomological action is the orthogonal projector onto $H_X^{\mathrm{att}}$ in degree three and zero in every other degree. In this category, the isomorphism $[\Gamma\circ e]\colon(S,e,0)\to(X,p_X,0)$ has inverse
  \[
    [\Gamma\circ e]^{-1}=d_\Gamma^{-1}[e\circ{}^te\circ{}^t\Gamma\circ p_X].
  \]
\end{lemma}
\begin{proof}
  In degree three, let $a$ be the Betti action of $\Gamma\circ e$ and $a^\dagger$ its cup-product adjoint. The image of $a$ is $H_X^{\mathrm{att}}$ by rationality and the nonzero holomorphic action. This plane is nondegenerate by Hodge--Riemann, so $aa^\dagger$ kills its orthogonal complement and restricts to a scalar on the plane. Indeed, self-adjointness on a symplectic plane forces equal diagonal entries and zero off-diagonal entries. The support of $e$ and duality make $Z_\Gamma$ act by zero outside degree three. Its intersection with the diagonal is therefore minus the trace on this plane \cite[Section~25]{milne-1998-6a5fd522}, identifying the scalar with $d_\Gamma$ and the trace with $2d_\Gamma$. For $0\ne\omega\in H^{3,0}(X)$, adjointness gives
  \[
    d_\Gamma=\frac{iQ_S(a^\dagger\omega,\overline{a^\dagger\omega})}
    {iQ_X(\omega,\overline{\omega})}>0,
  \]
  by Hodge--Riemann, since $a^\dagger$ is injective on the image plane and both holomorphic classes are primitive. Thus $p_X$ has the asserted total action. It is self-adjoint since $Z_\Gamma={}^tZ_\Gamma$.

  In homological motives, $[\Gamma\circ e]$ has source projector $e$ and target projector $p_X$. Composing it with the proposed inverse gives $d_\Gamma^{-1}[Z_\Gamma p_X]=[p_X]$. The reverse composite acts as the identity on the source summand because $a$ restricts to an isomorphism from $M$ onto $H_X^{\mathrm{att}}$. Faithfulness of Betti realization gives the reverse identity $[e]$.
\end{proof}

The homological class of $p_X$ is independent of the correspondence $\Gamma$: its total Betti action is the uniquely determined orthogonal projector onto $H_X^{\mathrm{att}}$ in degree three and zero otherwise.

\begin{theorem}\label{intrinsic-modular-realization}
  Let $X/\mathbb Q$ be a smooth projective threefold, with $h^{3,0}(X)=1$ and $\dim_{\mathbb Q}H_X^{\mathrm{att}}=2$. Suppose that a correspondence $\Gamma\in\operatorname{CH}^3(X\times_{\mathbb Q}X)_{\mathbb Q}$ acts on Betti cohomology as the $Q_X$-orthogonal projector $e_B$ onto $H_X^{\mathrm{att}}$ in degree three and as zero in every other degree. For each prime $\ell$, set $H_{X,\ell}^{\mathrm{att}}=\Gamma_{*,\ell} H^3_{\mathrm{\acute{e}t}}(X_{\overline{\mathbb Q}},\mathbb Q_\ell)$. Then there is a unique normalized weight-four newform $f\in S_4^{\mathrm{new}}(\Gamma_0(N))$ with trivial nebentypus and rational Fourier coefficients such that $H_{X,\ell}^{\mathrm{att}}\simeq M_{f,\ell}$ as $\mathbb Q_\ell[G_{\mathbb Q}]$-modules for every prime $\ell$.
\end{theorem}

\begin{proof}
  Compatibility with Betti--\'etale comparison makes $\Gamma_{*,\ell}$ idempotent, of rank two in degree three and zero in every other degree. Since $\Gamma$ is defined over $\mathbb Q$, its image $H_{X,\ell}^{\mathrm{att}}$ is a $G_{\mathbb Q}$-stable direct summand. Smooth projective cohomology makes it finitely ramified. Filtered comparison for smooth projective cohomology and passage to the direct summand give the de Rham property and Hodge--Tate weights $\{0,3\}$ \cite[Theorem~5.1]{bhatt-2018-4c44c329}\cite[Proposition~6.3.3]{brinon-2009-6811eb78}. It is therefore potentially semistable \cite[Theorem~9.3.5]{brinon-2009-6811eb78}. At a good prime $p\ne\ell$, its Frobenius eigenvalues have complex absolute value $p^{3/2}$ by purity \cite[Theorem~1.6]{deligne-1974-370b49c1}.

  For any $\ell$, an invariant line over a finite coefficient extension would be a de Rham character, hence a finite-order character times the character of $\mathbb Q_\ell(m)$ for some $m\in\mathbb Z$ \cite[p.~2]{dieulefait-2003-b01f8879}. At common good primes its geometric-Frobenius eigenvalues would have absolute value $p^{-m}$, contradicting purity of weight three. Thus $H_{X,\ell}^{\mathrm{att}}$ is absolutely irreducible for every $\ell$.

  The plane $H_X^{\mathrm{att}}$ is primitive, and the Hodge--Riemann relations make its restricted cup-product pairing nondegenerate. Comparison therefore gives a perfect alternating $G_{\mathbb Q}$-equivariant pairing with values in $\mathbb Q_\ell(-3)$. In dimension two this identifies $\bigwedge^2H_{X,\ell}^{\mathrm{att}}\simeq\mathbb Q_\ell(-3)$. Complex conjugation acts by $-1$ on this determinant, so the representation is odd.

  Choose one auxiliary prime $\ell_0\ge5$. The representation $H_{X,\ell_0}^{\mathrm{att}}$ is continuous, absolutely irreducible, odd, finitely ramified and potentially semistable at $\ell_0$, with distinct Hodge--Tate weights. The modularity theorem \cite[Theorem~1.0.4]{pan-2022-73f5d316} therefore produces a cuspidal eigenform up to twist. The choice $\ell_0\ge5$ avoids the additional residual restriction at three. Absorb the finite-order part of the geometric twist into the primitive newform twist and normalize the resulting form $f$. For a place $\lambda$ above $\ell_0$ of the coefficient field of $f$, the resulting isomorphism, after extending coefficients and using our cohomological convention, identifies $H_{X,\ell_0}^{\mathrm{att}}$ with $M_{f,\lambda}(m)$ for an integer $m$. If $k$ is the weight of $f$, comparison of Hodge--Tate weights gives $\{-m,k-1-m\}=\{0,3\}$, so $m=0$ and $k=4$. Since the determinant of the modular representation is the nebentypus times the character of $\mathbb Q_\ell(-3)$ in weight four, the restricted cup pairing forces trivial nebentypus. At this stage the Fourier coefficients of $f$ may lie in a number field larger than $\mathbb Q$.

  We determine that coefficient field from the all-degree action of $\Gamma$. Outside a finite set of primes, $X$ has smooth projective reduction and $\Gamma$ specializes. For such $p$ and every $\ell\ne p$, the specialized correspondence satisfies
  \[
    \sum_{j=0}^6(-1)^j
    \operatorname{Tr}\mleft(\Gamma_{*,\ell}\circ\operatorname{Frob}_p
    \mid H^j_{\mathrm{\acute{e}t}}(X_{\overline{\mathbb F}_p},\mathbb Q_\ell)\mright)
    =-\operatorname{Tr}(\operatorname{Frob}_p\mid H_{X,\ell}^{\mathrm{att}}).
  \]
  The trace formula for correspondences expresses the left-hand side as a rational intersection number involving the specialized cycle and Frobenius, independently of $\ell$ \cite[Section~25]{milne-1998-6a5fd522}. Denote the resulting common rational trace by $a_p$. The eigenvalues on the summand are among those on smooth proper $H^3$, so they are algebraic integers. Consequently $a_p\in\mathbb Z$, and $\det(1-T\operatorname{Frob}_p\mid H_{X,\ell}^{\mathrm{att}}) =1-a_pT+p^3T^2\in\mathbb Z[T]$ for every $\ell\ne p$.

  The modular isomorphism at $\ell_0$ gives $a_p(f)=a_p$ at all common good primes outside a finite set. Every coefficient-field conjugate $f^\sigma$ has the same coefficients at these primes. Strong multiplicity one identifies $f^\sigma$ with $f$, so all Fourier coefficients are rational.

  For any prime $\ell$, the representations $H_{X,\ell}^{\mathrm{att}}$ and $M_{f,\ell}$ have the same Frobenius polynomials at all common good primes outside a finite set. The former is absolutely irreducible by the purity argument, and the latter is the semisimple cuspidal modular realization. Chebotarev density and Brauer--Nesbitt give $H_{X,\ell}^{\mathrm{att}}\simeq M_{f,\ell}$ \cite[Proposition~9.4.5]{kedlaya-2020-4f9884b9}. Any other normalized newform realizing this system has the same Hecke eigenvalues at almost every prime and is equal to $f$ by strong multiplicity one.
\end{proof}

An explicit projector cycle or an independently identified zeta factor supplies the Frobenius traces needed to recognize $f$; a finite coefficient match alone does not establish the almost-everywhere equality in the modularity argument.

\begin{theorem}\label{conditional-derived-modular-comparison}
  Let $X/\mathbb Q$ be a smooth projective threefold with $h^{3,0}(X)=1$ and $\dim_{\mathbb Q}H_X^{\mathrm{att}}=2$. Suppose that a correspondence $e\in\operatorname{CH}^3(X\times_{\mathbb Q}X)_{\mathbb Q}$ acts as the cup-product orthogonal projector onto $H_X^{\mathrm{att}}$ in degree three and as zero in every other degree. Let $f$ be the resulting rational weight-four newform, and fix its Scholl realization $(S_f,e_f)/\mathbb Q$. Assume the Tate conjecture of \zcref{global-rational-tate-conjecture}. Then there is a cycle $\Gamma\in\operatorname{CH}^3(S_f\times_{\mathbb Q}X)_{\mathbb Q}$ with $e_{*,B}\Gamma_{*,B}|_{M_{f,B}}\ne0$.
\end{theorem}

\begin{proof}
  Fix a prime $\ell$ and an isomorphism $M_{f,\ell}\to H_{X,\ell}^{\mathrm{att}}$ from \zcref{intrinsic-modular-realization}. Extend it to degree-three cohomology by the source projector and target inclusion. Poincar\'e duality and K\"unneth identify the equivariant map with an invariant class in $H_\ell^6(S_f\times X)(3)$. By Tate it is a $\mathbb Q_\ell$-linear combination of classes of rational-coefficient cycles over $\mathbb Q$. At least one of these cycles acts nontrivially on $M_{f,\ell}$ after projection to $H_{X,\ell}^{\mathrm{att}}$, since their combination induces the chosen isomorphism. Betti--\'etale comparison gives the required Betti action.
\end{proof}

Given a comparison cycle, the chosen holomorphic forms determine the normalized isomorphism below.

\subsection{The canonical comparison and critical \texorpdfstring{$L$}{L}-values}\label{presentation-period-calibration}

\begin{theorem}\label{intrinsic-absolute-comparison}
  Let $X/\mathbb Q$ be a smooth projective threefold with $h^{3,0}(X)=1$ and $\dim_{\mathbb Q}H_X^{\mathrm{att}}=2$. Suppose there is a correspondence $e\in\operatorname{CH}^3(X\times_{\mathbb Q}X)_{\mathbb Q}$ acting as the cup-product orthogonal projector onto $H_X^{\mathrm{att}}$ in degree three and as zero in every other degree. Let $f$ be the unique normalized rational cuspidal newform of weight four and trivial nebentypus supplied by \zcref{intrinsic-modular-realization}, so that $M_{f,\ell}\simeq H_{X,\ell}^{\mathrm{att}}$ for every prime $\ell$.

  Fix a Scholl realization $M_f=(S_f,e_f,0)$ over $\mathbb Q$, put $M_X=(X,e,0)$, and suppose there is a cycle $\Gamma\in\operatorname{CH}^3(S_f\times_{\mathbb Q}X)_{\mathbb Q}$ such that $e_{*,B}\Gamma_{*,B}|_{M_{f,B}}\ne0$. Restriction to the holomorphic de Rham lines is an isomorphism
  \begin{equation}\label{absolute-holomorphic-restriction}
    \operatorname{Hom}(M_f,M_X)
    \xrightarrow{\sim}\operatorname{Hom}_{\mathbb Q}(\ell_f,\ell_X),
    \qquad u\longmapsto u_{\mathrm{dR}}|_{\ell_f}.
  \end{equation}
  Every nonzero $u$ in the source is an isomorphism.

  The projected cycle $[e\circ\Gamma\circ e_f]$ induces an isomorphism with inverse supplied by \zcref{correspondence-induced-projector}. Its realizations identify $M_{f,B}$ with $H_X^{\mathrm{att}}$ and $M_{f,\ell}$ with $e_{*,\ell}H_\ell^3(X)$ equivariantly for every $\ell$.

  Fix a field $\mathbb Q\subseteq k\subseteq\mathbb C$ and $0\ne\Omega\in\ell_X\otimes_{\mathbb Q}k$. There is a unique form-normalized isomorphism $\Phi_\Omega\in\operatorname{Hom}(M_f,M_X)\otimes_{\mathbb Q}k$ satisfying $(\Phi_\Omega)_{\mathrm{dR},k}(\omega_f)=\Omega$. It is rational if and only if $\Omega\in\ell_X$. Identify $M_{f,B}$ with $\mathbb Q^2$ in the basis dual to $(s_0,s_\chi)$. For any rational linear map $\mathcal B\colon H_X^{\mathrm{att}}\to\mathbb Q^m$, $m\ge1$, one has
  \begin{equation}\label{absolute-geometric-period-column}
    \mathcal B_{\mathbb C}I_X(\Omega)=-2\mathcal B_{\mathbb C}(\Phi_\Omega)_{B,\mathbb C}
    \begin{pmatrix}L(f,3)\\2iL(f,\chi_{-4},3)
    \end{pmatrix}.
  \end{equation}
\end{theorem}

\begin{proof}
  By rationality and conjugation, $e\circ\Gamma$ acts nontrivially on $\ell_f$. Apply \zcref{correspondence-induced-projector} to $e_f$ and $e\circ\Gamma$. The resulting projector agrees with $e$ on total Betti cohomology, so the target is $M_X$ in homological motives. Betti--\'etale comparison gives the equivariant realizations.

  If a rational $u$ vanishes on $\ell_f$, its Betti action vanishes on both conjugate Hodge lines, which exhaust $M_{f,\mathbb C}$. Faithfulness gives $u=0$. The projected cycle acts nontrivially on $\ell_f$ and the space on the right of \eqref{absolute-holomorphic-restriction} is one-dimensional, so restriction is an isomorphism. Evaluating at $\omega_f$ identifies the rational Hom space with $\ell_X$. Extending coefficients to $k$ gives $\Phi_\Omega$ and the rationality criterion. Naturality gives
  \begin{equation}\label{absolute-chosen-form-equality}
    I_X(\Omega)=(\Phi_\Omega)_{B,\mathbb C}I_f(\omega_f).
  \end{equation}
  Apply $\mathcal B_{\mathbb C}$ and use \eqref{algebraic-modular-period-column} to obtain \eqref{absolute-geometric-period-column}.
\end{proof}

\begin{corollary}
  For the form-normalized comparison $\Phi_\Omega$ of \zcref{intrinsic-absolute-comparison},
  \begin{equation}\label{absolute-critical-values-from-periods}
    L(f,3)=-\frac12s_{0,\mathbb C}
    \mleft((\Phi_\Omega)_{B,\mathbb C}^{-1}I_X(\Omega)\mright),
    \qquad
    L(f,\chi_{-4},3)=-\frac1{4i}s_{\chi,\mathbb C}
    \mleft((\Phi_\Omega)_{B,\mathbb C}^{-1}I_X(\Omega)\mright).
  \end{equation}
\end{corollary}

\begin{proof}
  Apply $(\Phi_\Omega)_{B,\mathbb C}^{-1}$ to \eqref{absolute-chosen-form-equality} and evaluate the two modular symbols using \eqref{algebraic-modular-period-column}. The inverse-image class is unchanged under rescaling $\Omega$.
\end{proof}

\subsection{Mirror comparison and rational period frames}

A rational period satisfying one genus-zero charge identity will determine the full rational prepotential frame. We construct it near the MUM point and compute its cup pairing for the toric comparison. Endpoint specialization enters in the next subsection.

Let $Y$ be a smooth projective Calabi--Yau threefold with primitive ample class $H$. Use the invariants $N_d^{\mathrm{GW}}$ of \zcref{genus-zero-period-convention}, with degrees measured by $H$. Put $Y_{111}=\int_YH^3$ and use its classical prepotential frame.

Fix a projective family $\pi\colon\mathcal X\to S$ over $\mathbb Q$, smooth over $S^\circ$, a rational coordinate $z$ at a MUM point and a rational rank-four subvariation. Write $\mathcal H_B$ for its rational local system and $\mathcal H_{\mathrm{dR}}$ for its algebraic de Rham bundle. Fix an algebraic section $\Omega_{\mathrm{alg}}\in H^0(S^\circ,F^3\mathcal H_{\mathrm{dR}})$ cyclic under $\nabla_\theta$, where $\theta=z\dd/\dd z$. Require that the scalar relation for $\nabla_\theta^4\Omega_{\mathrm{alg}}$, in the chosen coordinate and form gauge, has the normalized operator
\begin{equation}\label{constructed-picard-fuchs-operator}
  \mathcal L=\sum_{j=0}^s z^j\mathcal L_j(\theta),\qquad
  \mathcal L_j\in\mathbb Q[\theta],\qquad \deg \mathcal L_j\leq4,\qquad \mathcal L_0(\theta)=\theta^4.
\end{equation}
On every disc $\abs*{\rho}\leq r_0<1/2$, the quotients $\mathcal L_j(n-j+\rho)/(n+\rho)^4$ in \eqref{constructed-frobenius-recurrence} are uniformly bounded. The finite recurrence gives $\sup\abs*{\alpha_n(\rho)}\leq A^n$ for some $A>0$, so $h_0,\ldots,h_3$ converge near zero. In \eqref{mirror-quantities}, $h_0\in1+z\mathbb Q\cdbr{z}$, $h_1/h_0,V\in z\mathbb Q\cdbr{z}$ and $q\in z+z^2\mathbb Q\cdbr{z}$. In particular $q(z)$ has a local holomorphic inverse at zero.

Use the Frobenius column \eqref{constructed-prepotential-periods} and $\mathcal I(t),\mathcal I'(t)$ from \zcref{genus-zero-period-convention}, regarded as formal series in $q=e^{2\pi it}$. On a simply connected large-radius branch $U$, fix a rational flat functional $\Gamma_{\mathcal O}\in H^0(U,\mathcal H_B^\vee)$ and a fixed nonzero complex constant multiple $\Omega$ of $\Omega_{\mathrm{alg}}$. Here $\Gamma_{\mathcal O}(\Omega)(z)=\Gamma_{\mathcal O,z,\mathbb C}(I_z(\Omega(z)))$, where $I_z$ is the fibrewise de Rham--Betti comparison. On the fixed logarithm branch, suppose the period has the Frobenius expansion
\begin{equation}\label{single-charge-formal-identity}
  \frac{\Gamma_{\mathcal O}(\Omega)}{(2\pi i)^3}
  =w_0\mleft(-\frac{Y_{111}}6t^3-bt-2c_0
  +2\mathcal I(t(z))-t\mathcal I'(t(z))\mright)
  \quad\text{in }\mathbb C\cdbr{z}[\log z].
\end{equation}

Equality in $\mathbb C\cdbr{z}[\log z]$ is coefficientwise in both $z$ and $\log z$. This constant rescaling fixes the common period scale without changing the operator or rational functional, and a new coordinate normalized by $q=z+O(z^2)$ requires a new comparison.

\begin{proposition}\label{constructed-line-bundle-extraction}\label{constructed-rational-period-frame}
  Fix the cyclic rational rank-four subsystem, MUM coordinate, constant form gauge and frame \eqref{constructed-prepotential-periods} used in \eqref{single-charge-formal-identity}. Suppose that a rational flat functional $\Gamma_{\mathcal O}$ has the complete formal expansion \eqref{single-charge-formal-identity} for the virtual degrees $N_d^{\mathrm{GW}}$, with $Y_{111}\in\mathbb Q^\times$ and $Y_{011},b\in\mathbb Q$. Then the $q$-series defining $\mathcal I(t)$ and $\mathcal I'(t)$ converge near $q=0$, and
  \begin{equation}\label{single-charge-gw-pf}
    \mathcal I'(t(z))=Y_{111}\left(\frac{p_2}{p_0}-\frac{t^2}{2}\right)
    =\frac{Y_{111}V(z)}{(2\pi i)^2}.
  \end{equation}
  The MUM logarithm applied to $\Gamma_{\mathcal O}$ determines the full rational prepotential frame.
\end{proposition}
\begin{proof}
  Comparing the $z^0$ logarithmic polynomial in \eqref{single-charge-formal-identity} with the Frobenius basis fixes all four coefficients:
  \begin{equation}
    \frac{\Gamma_{\mathcal O}(\Omega)}{(2\pi i)^3}=-Y_{111}p_3-bp_1-2c_0p_0=-\mathcal F_0.
  \end{equation}
  Let $N=\log T^\vee$ for positive MUM monodromy on the rational flat functionals. On periods $Np_0=0$ and $Np_j=p_{j-1}$ for $1\le j\le3$. In \eqref{single-charge-formal-identity}, $N$ differentiates the explicit polynomial in $t$ because the translation $t\mapsto t+1$ fixes $q$ and $w_0$. Applying $N$ to the two expressions for the charge gives
  \[
    -Y_{111}p_2-bp_0=w_0\left(-\frac{Y_{111}}2t^2-b-\mathcal I'\right).
  \]
  Dividing by $p_0=w_0$ gives \eqref{single-charge-gw-pf}. Since $V(0)=0$ and $q'(0)=1$, the right-hand side is holomorphic in $q$ and vanishes at zero. Integrating $\mathcal I'/(2\pi iq)$ with zero constant gives convergence of $\mathcal I$.

  The powers of $N$ applied to $\Gamma_{\mathcal O}$ give the rational functionals
  \begin{align*}
    \delta_0&=-Y_{111}^{-1}N^3\Gamma_{\mathcal O},&
    \delta_1&=-Y_{111}^{-1}N^2\Gamma_{\mathcal O},\\
    \delta_{\mathcal F_1}&=N\Gamma_{\mathcal O}+2b\delta_0-Y_{011}\delta_1,&
    \delta_{\mathcal F_0}&=-\Gamma_{\mathcal O}.
  \end{align*}
  Their periods satisfy
  \begin{equation}\label{single-charge-prepotential-frame}
    \frac1{(2\pi i)^3}
    \begin{pmatrix}
      \delta_{\mathcal F_0}(\Omega)\\\delta_{\mathcal F_1}(\Omega)\\\delta_0(\Omega)\\\delta_1(\Omega)
    \end{pmatrix}
    =(\mathcal F_0,\mathcal F_1,z^0,z^1)^T=w_0\Pi(t).
  \end{equation}
  Their leading logarithmic degrees $3,2,0,1$ give a rational basis of $\mathcal H_B^\vee$.
  The rational flat coordinate map is
  \begin{equation}\label{constructed-betti-coordinates}
    \mathcal B_z\colon\mathcal H_{B,z}\xrightarrow{\sim}\mathbb Q^4,
    \qquad
    \mathcal B_z(x)=\mleft(\delta_{\mathcal F_0}(x),\delta_{\mathcal F_1}(x),\delta_0(x),\delta_1(x)\mright)^T.
  \end{equation}
\end{proof}

For a one-parameter Batyrev pair satisfying the toric hypotheses of \cite[Section~6.4, Theorems~6.9--6.10]{iritani-2011-7b0e4739}, we verify \eqref{single-charge-formal-identity} in the untwisted rank-four ambient sector. In particular, the lattice points of the Newton polytope must generate the character lattice of the mirror torus. These hypotheses concern this verification by toric mirror symmetry; the preceding reconstruction applies whenever \eqref{single-charge-formal-identity} holds.

For this toric comparison, the component of $\prod_j\Gamma(1+\delta_j)$ in $H^{2k}(Y)$ is $(2\pi i)^k$ times the corresponding component of $\widehat\Gamma_Y$ in \eqref{constructed-gamma-class}. For $d>0$, let $\angbr*{\cdots}_{0,n,d}$ denote the sum over curve classes of $H$-degree $d$ and put $\tau_k(\alpha)=\psi^k\operatorname{ev}^*\alpha$. With the fundamental solution $L_Y$ of \cite[eq.~(2.8)]{iritani-2011-7b0e4739}, auxiliary parameter $1$, divisor coordinate $2\pi itH$, and all other coordinates zero, the source formula is
\[
  \int_YL_Y(2\pi itH,1)\prod_j\Gamma(1+\delta_j)
  =\int_Ye^{-2\pi itH}\prod_j\Gamma(1+\delta_j)
  +\sum_{d=1}^\infty q^d
  \angbr*{\frac{e^{-2\pi itH}\prod_j\Gamma(1+\delta_j)}{-1-\psi},1}_{0,2,d}.
\]
The two-point virtual dimension is two and the Gamma class has no $H^2$ component. Dimension, string, dilaton and divisor equations therefore give
\begin{align*}
  \angbr*{\frac{e^{-2\pi itH}\prod_j\Gamma(1+\delta_j)}{-1-\psi},1}_{0,2,d}
  &=-\angbr*{\tau_2(1),1}_{0,2,d}
  -2\pi it\angbr*{\tau_1(H),1}_{0,2,d}\\
  &=-\angbr*{\tau_1(1)}_{0,1,d}-2\pi it\angbr*{H}_{0,1,d}\\
  &=(2-2\pi idt)N_d^{\mathrm{GW}}.
\end{align*}
The unit-section comparison of \cite[Theorem~6.9]{iritani-2011-7b0e4739} uses the fixed residue form. On the one-parameter untwisted rank-four sector, the period of the oriented compact point-class cycle of \cite[Theorem~6.10]{iritani-2011-7b0e4739}, divided by $(2\pi i)^3$, is holomorphic, solves $\mathcal L$ and has value $1$ at $z=0$. The MUM recurrence makes it equal to $w_0$. The K-group framing of \cite[Definition~3.6]{iritani-2011-7b0e4739} then gives the complete coefficientwise charge
% Published locators checked in Theorems 6.9--6.10, pp. 2947--2950;
% the orientation correction is footnote 16.
\begin{equation}\label{toric-structure-sheaf-charge}
  \begin{split}
    \frac{\Gamma_{\mathcal O}(\Omega)}{(2\pi i)^3w_0}
    &=\int_Ye^{-tH}\widehat\Gamma_Y
    +\frac1{(2\pi i)^3}\sum_{d=1}^\infty(2-2\pi idt)N_d^{\mathrm{GW}}q^d\\
    &=-\frac{Y_{111}}6t^3-bt-2c_0
    +2\mathcal I(t)-t\mathcal I'(t).
  \end{split}
\end{equation}
Thus the structure-sheaf period verifies \eqref{single-charge-formal-identity} coefficientwise in the formal degree series.

\begin{proposition}\label{toric-period-pairing}
  For the rank-four Batyrev comparison used in \eqref{toric-structure-sheaf-charge}, take $\Gamma_{\mathcal O}$ to be the cycle corresponding to $\mathcal O_Y$ under its $K$-group framing. Then the ordinary cup matrix in the rational coordinates \eqref{constructed-betti-coordinates} is
  \[
    J_4=
    \begin{pmatrix}0&I_2\\-I_2&0
    \end{pmatrix}.
  \]
\end{proposition}

\begin{proof}
  On a mirror fibre $\mathcal X_z$, use
  \[
    a\cdot b=\int_{\mathcal X_z}\operatorname{PD}(a)\cup\operatorname{PD}(b),
    \qquad
    \int_a\omega=\int_{\mathcal X_z}\omega\cup\operatorname{PD}(a).
  \]
  With this convention, the mirror comparison identifies homological intersection with minus the Euler pairing \cite[Definition~6.5, Theorem~6.9, eqs.~(6.8)--(6.9)]{iritani-2011-7b0e4739}. Positive MUM continuation corresponds to tensoring by $\mathcal O_Y(-H)$ \cite[proof of Theorem~5.7]{iritani-2011-7b0e4739}. For its logarithm $N$ on homology, the inverse images of $N^j\Gamma_{\mathcal O}$ have Chern characters $(-H)^j$. Riemann--Roch with $\operatorname{Td}(Y)=1+c_2(Y)/12$ gives, for $0\le i,j\le3$,
  \[
    (N^i\Gamma_{\mathcal O})\cdot(N^j\Gamma_{\mathcal O})
    =-\int_Y H^i(-H)^j\left(1+\frac{c_2(Y)}{12}\right).
  \]
  Substituting the coordinate functionals of \zcref{constructed-rational-period-frame} gives
  \[
    \delta_{\mathcal F_0}\cdot\delta_0
    =\delta_{\mathcal F_1}\cdot\delta_1=1,\qquad
    \delta_{\mathcal F_0}\cdot\delta_{\mathcal F_1}=-2b+2b=0,
  \]
  and all other upper-triangular pairings vanish. The intersection matrix is $J_4$, so the cup matrix in the dual coordinates is $-J_4^{-1}=J_4$.
\end{proof}

In applications, the full local comparison is the rational flat coordinate map \eqref{constructed-betti-coordinates} together with the two identities near the large-radius point
\begin{equation}\label{local-full-mirror-frame}
  \frac{\mathcal F_1}{z^0}=(F^{(0)})'(t),\qquad
  \frac{\mathcal F_0}{z^0}=2F^{(0)}(t)-t(F^{(0)})'(t),\qquad
  t=\frac{z^1}{z^0}.
\end{equation}
The potential $F^{(0)}$ here is convergent in the selected classical frame. \zcref{constructed-line-bundle-extraction} supplies this comparison from \eqref{single-charge-formal-identity}. A direct construction of the rational flat map $\mathcal B$ and \eqref{local-full-mirror-frame} is sufficient in its place. Equality of scalar Gromov--Witten and Picard--Fuchs coefficients does not by itself supply the rational functionals.

\subsection{The continued Gromov--Witten partition identity}\label{geometric-genus-zero-data}

Fix the mirror data of the preceding subsection, with $Y_{111}>0$, a nonzero constant multiple $\Omega$ of $\Omega_{\mathrm{alg}}$, and a rational flat coordinate map $\mathcal B$ satisfying \eqref{local-full-mirror-frame} and $\mathcal B_{z,\mathbb C}(\Omega(z))=(2\pi i)^3w_0\Pi(t)$ near the large-radius point. The first identity in \eqref{local-full-mirror-frame} and the raw frame give $\sum_{d\ge1}dN_d^{\mathrm{GW}}q(z)^d=Y_{111}V(z)$. Holomorphy of $V$ at zero and $q'(0)=1$ give convergence of this degree series near $q=0$. For $x_0\in\mathbb R$ and $y_0>0$, the holomorphic implicit-function theorem gives the local solution
\begin{equation}\label{scalar-normalized-inversion}
  (T-x_0)^2+y_0^2=\frac1{2\pi^2Y_{111}}
  \sum_{d\ge1}dN_d^{\mathrm{GW}}s^de^{2\pi idT},
  \qquad T(0)=x_0+iy_0,
\end{equation}
since the derivative in $T$ at the initial point is $2iy_0$.

We prove the expansion of \cite[Section~6.2]{candelas-2021-56ab2e45} and \cite[eqs.~(42)--(44)]{candelas-2025-ad265528} by formal inversion.

\begin{lemma}\label{partition-expansion}
  Fix the genus-zero coefficients $N_d^{\mathrm{GW}}$ of \zcref{genus-zero-period-convention}, $Y_{111}>0$, $x_0\in\mathbb R$ and $y_0>0$. For a partition $\mathfrak p=(1^{\mu_1}2^{\mu_2}\cdots)$, put $l(\mathfrak p)=\sum_k\mu_k$ and $\widetilde N_{\mathfrak p}=\prod_k(kN_k^{\mathrm{GW}})^{\mu_k}/\mu_k!$. The formal solution of \eqref{scalar-normalized-inversion} with $T(0)=x_0+iy_0$ is
  \begin{equation}\label{normalized-scalar-partition-expansion}
    T(s)=x_0+iy_0-\frac{i}{\sqrt{2\pi^3Y_{111}}}
    \sum_{j\ge1}s^je^{2\pi ijx_0}
    \sum_{\mathfrak p\vdash j}\widetilde N_{\mathfrak p}
    \left(\frac{j}{2\pi y_0Y_{111}}\right)^{l(\mathfrak p)-1/2}
    K_{l(\mathfrak p)-1/2}(2\pi jy_0),
  \end{equation}
  where $K_\nu$ is the modified Bessel function of the second kind.
\end{lemma}
\begin{proof}
  Put $G(q)=\sum_{d\ge1}dN_d^{\mathrm{GW}}q^d$ and choose the formal square root
  \[
    \mathfrak{w}(q)=i\sqrt{y_0^2-\frac{G(q)}{2\pi^2Y_{111}}},\qquad \mathfrak{w}(0)=iy_0.
  \]
  The series $s(q)=q e^{-2\pi i(x_0+\mathfrak{w}(q))}$ has a nonzero linear coefficient, so its formal inverse $q(s)$ gives the unique solution $T(s)=x_0+\mathfrak{w}(q(s))$ of \eqref{scalar-normalized-inversion} with the specified initial value. For $j\ge1$, Lagrange inversion followed by differentiation of the exponential gives
  \[
    [s^j]T=\frac{e^{2\pi ijx_0}}{2\pi ij}[q^j]e^{2\pi ij\mathfrak w(q)}.
  \]
  For $a,b>0$ and $l\ge0$, on the branch $\sqrt{b^2-u}=b$ at zero, the half-integer value and the derivative identity for $K_\nu$ \cite[eqs.~(10.29.4), (10.39.2)]{dlmf-2010-fa11482d} give
  \begin{align*}
    [u^l]e^{-a\sqrt{b^2-u}}
    &=\frac{(-1)^l}{2^ll!}(b^{-1}\partial_b)^le^{-ab}\\
    &=\frac{a}{\sqrt\pi l!}
    \left(\frac{a}{2b}\right)^{l-1/2}K_{l-1/2}(ab).
  \end{align*}
  Substitute $a=2\pi j$, $b=y_0$ and $u=G(q)/(2\pi^2Y_{111})$. Since $G(0)=0$, only $1\leq l\leq j$ contributes. The identity $[q^j]G(q)^l=l!\sum_{\mathfrak p\vdash j,l(\mathfrak p)=l}\widetilde N_{\mathfrak p}$ cancels $1/l!$ and gives \eqref{normalized-scalar-partition-expansion}.
\end{proof}

When $\sum_{d\ge1}dN_d^{\mathrm{GW}}q^d$ converges near zero this is an equality of holomorphic functions near $s=0$. In the geometric identities, $T_\gamma(1)$ is the value obtained by continuing the function defined by the full series along the inversion path induced by the selected period branch. Each inner partition sum is finite; ordinary convergence of the Taylor series at $s=1$ is a separate assertion.

Let $\gamma\colon[0,1)\to S^\circ$ run from the large-radius base point $z_b$ towards an endpoint $z_*$. Use the selected Frobenius determination of $t$ and assume the normalized column continues holomorphically along the path with finite limit $\Pi_*$. At a smooth endpoint take $X=\mathcal X_{z_*}$ and let $\sigma_\gamma$ be the polarized projection $\operatorname{proj}\colon H_B^3(X,\mathbb Q)\to\mathcal H_{B,z_*}$ followed by rational parallel transport from $z_*$ to $z_b$ along $\gamma$. At a singular endpoint, choose a smooth projective pure model $X/\mathbb Q$ and a rational map $\sigma_\gamma\colon H_B^3(X,\mathbb Q)\to\mathcal H_{B,z_b}$. Put $\mathcal B_\gamma=\mathcal B_{z_b}\sigma_\gamma$. At a smooth endpoint the rank-four subsystem contains the holomorphic line and its conjugate, so $\operatorname{proj}$ fixes $H_X^{\mathrm{att}}$. The normalized column lies in the transported holomorphic line by mirror comparison, including at removable quotient zeros by holomorphic extension. Continuity at the endpoint therefore gives $\Pi_*\in\mathcal B_{\gamma,\mathbb C}(H^{3,0}(X))$. In the singular case require
\begin{equation}\label{marked-pure-specialization}
  \sigma_\gamma|_{H_X^{\mathrm{att}}}\ \text{injective},\qquad
  \Pi_*\in\mathcal B_{\gamma,\mathbb C}(H^{3,0}(X)).
\end{equation}

\begin{corollary}\label{geometric-partition-identity}
  Use the local mirror comparison \eqref{local-full-mirror-frame} with the rational flat frame and constant form gauge specified in \zcref{geometric-genus-zero-data}. Along $\gamma$, take the marking $\mathcal B_\gamma$, subject to \eqref{marked-pure-specialization} at a singular endpoint. Assume that the normalized period column continues holomorphically along $\gamma$ with finite limit $\Pi_*$.

  Let $X/\mathbb Q$ have $h^{3,0}(X)=1$, $\dim_{\mathbb Q}H_X^{\mathrm{att}}=2$, and an algebraic orthogonal projector onto $H_X^{\mathrm{att}}$ in degree three, zero in other degrees. Let $f$ be the resulting newform and take a nonzero rational algebraic comparison $u_B\colon M_{f,B}\to H_X^{\mathrm{att}}$. In the source basis dual to $(s_0,s_\chi)$ put $(c_{ij})=\mathcal B_\gamma u_B$. Suppose $\Pi_{*,3}=1$ and $t_*=\Pi_{*,4}$ has positive imaginary part. The last two rows of the matrix are independent and determine $(\Lambda,\Upsilon)\in\mathbb Q^2$ by $(0,1,-\Lambda,-\Upsilon)(c_{ij})=0$. Put $x_0=-(Y_{011}+\Upsilon)/Y_{111}$ and $y_0^2=2(\Lambda-b)/Y_{111}-x_0^2$.
  Suppose $y_0^2>0$, put $y_0=\sqrt{y_0^2}$, and assume the branch $\mathfrak w^2=-y_0^2+V/(2\pi^2)$ starting at $iy_0$ continues holomorphically along the initial MUM segment and $\gamma$ with a finite limit of positive imaginary part.

  Then $T$ has expansion \eqref{normalized-scalar-partition-expansion} near zero and its continued value on the induced inversion path satisfies
  \begin{equation}\label{geometric-partition-l-values}
    T_\gamma(1)=\frac{c_{41}L(f,3)+2ic_{42}L(f,\chi_{-4},3)}
    {c_{31}L(f,3)+2ic_{32}L(f,\chi_{-4},3)}.
  \end{equation}
  The denominator is nonzero. Rescaling the comparison leaves the ratio unchanged.
\end{corollary}

\begin{proof}
  The nonzero rational Hodge map $u_B$ is an isomorphism because its complexified kernel is a conjugation-stable sum of the two Hodge lines. The image under $(c_{ij})$ of the critical column in \zcref{uniform-critical-column} spans the marked holomorphic line. Since $\Pi_{*,3}=1$, the denominator is nonzero and the ratio equals $t_*$. A rational dependence of the last two rows would force a nontrivial rational linear relation between $1$ and $t_*$, contrary to $\rmIm t_*>0$. Hence the two rows form a basis of the row space and give the asserted pair.

  The raw formula $\mathcal F_1/z^0=b-Y_{011}t-Y_{111}t^2/2+Y_{111}V/(4\pi^2)$ has endpoint $\Lambda+\Upsilon t_*$. Consequently
  \[
    \frac{V_*}{2\pi^2}
    =t_*^2+\frac{2(Y_{011}+\Upsilon)}{Y_{111}}t_*+\frac{2(\Lambda-b)}{Y_{111}}
    =(t_*-x_0)^2+y_0^2.
  \]
  Positive imaginary parts give $\mathfrak w_*=t_*-x_0$. On the surface of the continued branches with the MUM point adjoined, put $s=qe^{-2\pi i(x_0+\mathfrak w)}$. At zero, $s(0)=0$ and $s'(0)=e^{-2\pi i(x_0+iy_0)}\ne0$. The mirror identity and implicit uniqueness identify $x_0+\mathfrak w\circ s^{-1}$ with the local solution of \eqref{scalar-normalized-inversion}.

  The critical points of $s$ are discrete. On successive compact portions of the lifted period path, take small detours avoiding them, with the oscillations of $s,\mathfrak w$ on the detours tending to zero along the tail. The local inverses of $s$ continue $x_0+\mathfrak w\circ s^{-1}$ along the image of the modified path under $s$, and the detours preserve both limits. Thus $s_*=e^{2\pi i(t_*-x_0-\mathfrak w_*)}=1$ and $T_\gamma(1)=x_0+\mathfrak w_*=t_*$. The local expansion follows from \zcref{partition-expansion}.
\end{proof}

\section{The Fermat sextic}\label{aesz8-sextic-geometric}

For the Fermat sextic the Hecke character determines the modular source and its base equation. We use them to construct the motivic isomorphism $\Phi_6$ sending $\omega_f$ to $\Omega_6$. We then determine its marked Betti matrix from the two traces and the cup forms in the fixed coordinates.

\subsection{Fermat geometry and the modular form}

Let $Y\subset\mathbb P(2,1,1,1,1)$ be a smooth sextic missing the stacky point, and put $H=c_1(\mathcal O_Y(1))$. The weighted Euler and normal sequences give $c(T_Y)=(1+2H)(1+H)^4/(1+6H)$ and $\int_YH^3=3$. Choosing the prepotential frame with $Y_{011}=1/2$ gives the classical constants
\begin{equation}\label{sextic-classical-constants}
  (Y_{111},c_2(Y)\mathbin{\cdot}H,\chi(Y),Y_{011}) =(3,42,-204,1/2), \qquad b=\frac74, \qquad c_0=-\frac{102\zeta(3)}{(2\pi i)^3}.
\end{equation}
The mirror family is obtained from
\begin{equation}\label{sextic-covering-family}
  \widetilde X_\psi\colon u^3+x_1^6+x_2^6+x_3^6+x_4^6 -6\psi ux_1x_2x_3x_4=0 \subset\mathbb P(2,1,1,1,1)
\end{equation}
by the diagonal group
\begin{equation}\label{sextic-diagonal-group}
  G= \frac{\setc*{(\lambda_0,\ldots,\lambda_4)\in \mu_3\times\mu_6^4}{\prod_{j=0}^4\lambda_j=1}} {\setc*{(a^2,a,a,a,a)}{a\in\mu_6}}.
\end{equation}
Thus $X_\psi$ denotes a fixed maximal projective crepant partial (MPCP) desingularization of $\widetilde X_\psi/G$, and $X_0$ is the resolved quotient at $\psi=0$.  We fix
\begin{equation}\label{sextic-endpoint-branch}
  z=(6\psi)^{-6}=-R,\qquad 6\psi=R^{-1/6}e^{-i\pi/6},\qquad R\to\infty,
\end{equation}
with the upper negative-axis determination of the continued periods.

On the quotient torus put
\[
  t_0=\frac{u^2}{x_1x_2x_3x_4},
  \qquad
  t_i=\frac{x_i^6}{ux_1x_2x_3x_4}\quad(1\leq i\leq4).
\]
The equations are
\begin{equation}\label{sextic-quotient-torus}
  t_0^2t_1t_2t_3t_4=1, \qquad t_0+t_1+t_2+t_3+t_4=6\psi.
\end{equation}
The exponent matrix has diagonal $(1,3,6,6)$, so the quotient has degree $108$.  It is the Batyrev Laurent mirror for the exponent relation $(2,1,1,1,1)$ and is defined over $\mathbb Q$ \cite[Theorem~4.2.2 and Corollary~4.2.3]{batyrev-1994-5b5891d2}. Its Newton simplex $\operatorname{conv}\{e_0,e_1,e_2,e_3,(-2,-1,-1,-1)\}$ has only its five vertices and the origin as lattice points, while its proper faces have no interior correction terms.  \cite[Corollary~4.5.1]{batyrev-1994-5b5891d2} gives $h^{2,1}(X_0)=6-5=1$. The torus critical equations give $4\psi^6=1$, and the proper faces have no critical points. Hence $\psi=0$ is a smooth fibre.  The degree-$108$ quotient map $\widetilde X_0\dashrightarrow X_0$ is rational over $\mathbb Q$; a resolution of its graph defines pullback and pushforward on every realization.

Choose the primitive invariant logarithmic form in the MUM orientation and put
\begin{equation}\label{sextic-mum-residue}
  \Omega_X=\operatorname{Res}
  \frac{\bigwedge_{j=0}^3\dd t_j/t_j}{\sum_i t_i-6\psi},
  \qquad
  \Omega_{\mathrm{MUM}}=-6\psi\Omega_X.
\end{equation}
The anticanonical hypersurface and its crepant MPCP resolution have trivial canonical bundle by adjunction. Thus the smooth projective connected fibre $X_0$ has $h^{3,0}(X_0)=1$, and the residue extends to a nonzero regular generator $\Omega_6=\Omega_X(0)\in F^3H^3_{\mathrm{dR}}(X_0/\mathbb Q)$ in this fixed orientation.

Put $\omega=e^{2\pi i/3}$, $K=\mathbb Q(\omega)$ and $\mathcal O_K=\mathbb Z[\omega]$. In this fixed embedding in $\mathbb C$, one has $1+2\omega=i\sqrt3$. Every nonzero integral ideal prime to $6\mathcal O_K$ has a unique generator $\alpha\equiv1\pmod{3\mathcal O_K}$: the ring $\mathcal O_K$ is norm-Euclidean, and its six units represent $(\mathcal O_K/3\mathcal O_K)^\times$. Define
\begin{equation}\label{sextic-geometric-hecke-character}
  \varepsilon(\alpha)\in\mu_3,\qquad
  \varepsilon(\alpha)\equiv\alpha\pmod{2\mathcal O_K},\qquad
  \Psi((\alpha))=\alpha^3\varepsilon(\alpha)^{-1}.
\end{equation}
Products of primary generators are primary, so $\Psi$ is multiplicative. Moreover, $\Psi((\alpha))=\alpha^3$ when $\alpha\equiv1\pmod{6\mathcal O_K}$. Define $f$ by the Hecke theta series
\begin{equation}\label{sextic-hecke-theta-series}
  f(\tau)=\sum_{\mathfrak a}\Psi(\mathfrak a)q^{N\mathfrak a},
  \qquad q=e^{2\pi i\tau},\qquad \tau\in\mathfrak H,
\end{equation}
where the sum runs over nonzero integral ideals prime to $6\mathcal O_K$.

For $x\in\mathcal O_K$ prime to $6$, write $x=u\alpha$ with $\alpha$ primary and $u\in\mu_6$, and let $\varepsilon(x)\in\mu_3$ be its representative modulo $2$. Then $\Psi((x))/x^3=u\varepsilon(x)^{-1}$ defines a finite character, extended by zero on nonunit classes modulo $6$. Its factor at $2$ is the nontrivial cubic character. Its factor at $3$ is nontrivial on the class of $\omega$ in the kernel of reduction modulo $i\sqrt3$, so the local conductors are $(2)$ and $(3)$ and the full conductor is $6\mathcal O_K$. For integers $n$ prime to $6$, one has $\Psi((n))/n^3=\chi_{-3}(n)$, where $\chi_{-3}$ is the quadratic character modulo three. Since $K$ has discriminant $-3$, the Hecke theta theorem makes $f$ a weight-$4$ newform of exact level $3N(6\mathcal O_K)=108$ and trivial nebentypus $\chi_{-3}^2$ \cite[Theorem~1.4]{schutt-2009-5629f886}. Conjugation permutes the ideals of a fixed norm and satisfies $\Psi(\overline{\mathfrak a})=\overline{\Psi(\mathfrak a)}$, so each coefficient lies in $\mathcal O_K\cap\mathbb R=\mathbb Z$. The norm-one ideal gives the coefficient $1$ at $q$, hence $f$ is normalized. Its LMFDB label is $\mathtt{108.4.a.c}$ \cite{lmfdb-2016-9e768d50}.

Put $g(\xi)=f(\xi/6)$, $Q=e^{2\pi i\xi}$, and $\zeta=e^{\pi i/3}$, with $Q^{1/6}=e^{\pi i\xi/3}$. For a primary $\alpha$ prime to $6$, its norm is odd and congruent to $1$ modulo $3$. Thus $g(\xi)=Q^{1/6}(1+O(Q))$ has a convergent expansion in powers of $Q$ after its leading factor, and $g(\xi+1)=\zeta g(\xi)$.

Put $T\xi=\xi+1$ and $W_3\xi=-1/(3\xi)$. Reduction modulo $3$ identifies the cosets of $\Gamma_0(3)$ with $\mathbb P^1(\mathbb F_3)$, so its index is four. Bezout's identity identifies a reduced rational cusp with $\infty$ when its denominator is divisible by $3$, and with $0$ otherwise; their widths are one and three. Put
\[
  U_3=W_3T^{-1}W_3^{-1}
  =
  \begin{pmatrix}1&0\\3&1
  \end{pmatrix}.
\]
For the first column $(a,3c)^T$ of a matrix in $\Gamma_0(3)$, one has $3\nmid a$. Left multiplication by a power of $U_3$ makes $\abs*{c}\leq\abs*{a}/2$. If $c\ne0$, a power of $T$ then makes the new $\abs*{a}$ at most $3\abs*{c}/2\leq3\abs*{a}/4$. Since $a$ remains nonzero modulo $3$, the descent ends at $c=0$, where the matrix is $\pm T^m$. Thus $T,U_3$ generate $\Gamma_0(3)$ modulo its centre.

The newform transformation law \cite[Theorem~4.6.15(2)]{miyake-2006-05f1a867} and the reality of the coefficients give, for some $\epsilon_f\in\{1,-1\}$,
\[
  f\mleft(-\frac1{108\tau}\mright)=\epsilon_f108^2\tau^4f(\tau),
  \qquad g(W_3\xi)=9\epsilon_f\xi^4g(\xi).
\]
Let $\eta$ denote the Dedekind eta function. Put $A_\eta(\xi)=\eta(\xi)^{10}/\eta(3\xi)^2$, $B_\eta(\xi)=\eta(3\xi)^{10}/\eta(\xi)^2$ and $\mathfrak h=729B_\eta/A_\eta$. The eta transformations give
\begin{align*}
  A_\eta(T\xi)&=\zeta A_\eta(\xi),& B_\eta(T\xi)&=\zeta B_\eta(\xi),\\
  A_\eta(W_3\xi)&=243\xi^4B_\eta(\xi),& B_\eta(W_3\xi)&=\frac{\xi^4}{3}A_\eta(\xi),
\end{align*}
and $\mathfrak h(W_3\xi)=729/\mathfrak h(\xi)$. By \cite[Tables~2 and~7]{maier-2009-f70ffeb6}, $\mathfrak h$ is a Hauptmodul of $X_0(3)$. The quotient $g/A_\eta$ is invariant under $\Gamma_0(3)$ and satisfies $g(W_3\xi)/A_\eta(W_3\xi)=27\epsilon_f g(\xi)/(\mathfrak h(\xi)A_\eta(\xi))$.
Since $\mathfrak h=729Q+O(Q^2)$ and $g/A_\eta=1+O(Q)$, this quotient is holomorphic away from a simple pole at $0$, hence affine in $\mathfrak h$. The Fricke law gives $g=A_\eta+27\epsilon_f B_\eta$.

For $\alpha=-2-3\omega$, one has $N(\alpha)=7$, $\alpha\equiv1\pmod3$, $\varepsilon(\alpha)=\omega$, and $\alpha^3=19+18\omega$. The two ideals of norm seven are $(\alpha)$ and $(\overline\alpha)$: they are distinct because a common prime divisor would divide $\alpha-\overline\alpha=-3i\sqrt3$, whose norm is prime to seven. Consequently
\[
  \operatorname{Tr}_{K/\mathbb Q}\mleft(\alpha^3\omega^{-1}\mright)
  =\operatorname{Tr}_{K/\mathbb Q}(-1-19\omega)=17.
\]
The products give $A_\eta=Q^{1/6}(1-10Q+O(Q^2))$ and $B_\eta=Q^{7/6}(1+O(Q))$, so $17=-10+27\epsilon_f$ and $\epsilon_f=1$. Thus
\begin{equation}\label{sextic-geometric-modular-form}
  f(\tau)=\frac{\eta(6\tau)^{10}}{\eta(18\tau)^2}
  +27\frac{\eta(18\tau)^{10}}{\eta(6\tau)^2},
\end{equation}
and
\begin{equation}\label{sextic-hecke-transformations}
  g(\xi+1)=\zeta g(\xi),\qquad
  g\mleft(-\frac1{3\xi}\mright)=9\xi^4g(\xi),
  \qquad \xi\in\mathfrak H.
\end{equation}

Put
\begin{equation}\label{sextic-weight-two-base-form}
  \mathscr A(\xi)=1+12\sum_{n=1}^\infty Q^n\sum_{d\mid n,3\nmid d}d.
\end{equation}
Use the normalized Eisenstein series $E_4=1+240Q+O(Q^2)$ and $E_6=1-504Q+O(Q^2)$.
The same tables give $j=(\mathfrak h+27)(\mathfrak h+3)^3/\mathfrak h$. Put $\lambda=\mathfrak h/(\mathfrak h+27)$. Differentiating the eta product gives $\mathscr A=\dd\log\mathfrak h/\dd\log Q$. Logarithmic differentiation of $j$, together with $\dd\log j/\dd\log Q=-E_6/E_4$ and $E_6^2=E_4^3(1-1728/j)$, yields
\[
  E_4=(1+8\lambda)\mathscr A^2,\qquad
  E_6=(1-20\lambda-8\lambda^2)\mathscr A^3.
\]
Differentiating $\mathfrak h(W_3\xi)=729/\mathfrak h(\xi)$ gives $\mathscr A(W_3\xi)=-3\xi^2\mathscr A(\xi)$.
Put $\sigma(\tau)=f(\tau)/\mathscr A(6\tau)^2$ and $v=1-2\lambda(6\tau)$.

\begin{proposition}\label{sextic-hecke-base}
  The function $\lambda$ induces an isomorphism $X_0(3)\to\mathbb P^1$, taking the cusps $\infty,0$ and the cubic elliptic point to $0,1,\infty$, respectively. It satisfies $\lambda(W_3\xi)=1-\lambda(\xi)$. For $\xi=6\tau$, the functions $\sigma$ and $v$ satisfy
  \[
    27\sigma^6=\lambda(1-\lambda),\qquad v^2=1-108\sigma^6.
  \]
  The first equation defines an extension of degree six of $\mathbb C(\lambda)$.
\end{proposition}

\begin{proof}
  The Hauptmodul $\mathfrak h$ takes the values $0,\infty,-27$ at the two cusps and the cubic elliptic point, so the formula for $\lambda$ and the Fricke law give the stated isomorphism and values. Use $\eta^{24}=E_4^3/j$, $A_\eta^6=729\eta^{48}/\mathfrak h$ and $g=A_\eta(1+\mathfrak h/27)$ to obtain
  \[
    \frac{g^6}{\mathscr A^{12}}
    =\frac{\mathfrak h}{(\mathfrak h+27)^2}
    =\frac{\lambda(1-\lambda)}{27}.
  \]
  Substitution of $\xi=6\tau$ proves both equations; the first is Eisenstein at $\lambda=0$.
\end{proof}

The differential of the coordinate has the normalization
\begin{equation}\label{sextic-cubic-theta-differential}
  \dd\lambda=6\lambda(1-\lambda)\mathscr A\frac{\dd q}{q}
  =162\sigma^6\mathscr A\frac{\dd q}{q},
\end{equation}
where $q=e^{2\pi i\tau}$ and $\lambda,\mathscr A$ are evaluated at $\xi=6\tau$.

Along $\tau=iy$, $y>0$, the eta product is positive, so \eqref{sextic-geometric-modular-form} gives $f(iy)>0$. The coefficients in \eqref{sextic-weight-two-base-form} are positive, so $\mathscr A(6iy)>0$ and $\sigma(iy)>0$. The real Fourier expansions make $\lambda(6iy)$ real, and $\lambda(6iy)(1-\lambda(6iy))=27\sigma(iy)^6>0$ gives
\begin{equation}\label{sextic-hecke-source-positivity}
  0<\lambda(6iy)<1,\qquad \sigma(iy)>0,\qquad f(iy)>0.
\end{equation}

Let $(S_f,e_f)/\mathbb Q$ be the fixed smooth projective Scholl realization of $f$, with $M_f=(S_f,e_f,0)$ and $\omega_f$ normalized as in \eqref{algebraic-modular-period-column}.

\subsection{Algebraic correspondence and motivic comparison}\label{sextic-explicit-common-covers}\label{sextic-explicit-kummer-descent}

A common cover supplies the correspondence required by \zcref{intrinsic-absolute-comparison}. On the modular side, the product of two elliptic quotients of the universal curve modulo simultaneous negation gives a Kummer model $\mathscr M$ over a genus-two base. On the sextic side, the residue form on $X_0$ descends to a model $\mathcal Q$ along a degree-two rational map. Both models receive maps from a threefold $\mathscr H$ over that base. The reciprocal quotients give the two elliptic cubics of $\mathscr M$. A quadric secant construction using the cubed fibre coordinates, together with the original roots, gives the map to $\mathcal Q$.

The ratio of the pulled-back forms is a rational function on $\mathscr H$. The reciprocal involutions and an order-three symmetry give the traces of this ratio and its reciprocal. The resulting graph correspondence acts nontrivially on the modular holomorphic form. Adjointness turns the two trace identities into the factor relating the cup pairings in \zcref{sextic-absolute-de-rham-calibration}.

The modular source uses the functions $\sigma,\lambda$ of \zcref{sextic-hecke-base} at $\xi=6\tau$. With $v=1-2\lambda$, their identity gives $v^2=1-108\sigma^6$.
Thus the base is the smooth projective model of
\begin{equation}\label{sextic-genus-two-base}
  \mathscr C\colon v^2=1-108\sigma^6,\qquad \lambda=\frac{1-v}{2},\qquad \beta=108\sigma^6.
\end{equation}
Define the two cubics by
\begin{align*}
  P_{6,\lambda}(x)&=x^3-\frac{1+8\lambda}{48}x+\frac{1-20\lambda-8\lambda^2}{864},\\
  P_{18,\lambda}(x)&=-\frac1{27}P_{6,1-\lambda}(-3x).
\end{align*}
Let $\mathscr M$ be a smooth projective model of the threefold over $\mathscr C$ given by
\[
  Z^2=P_{6,\lambda}(X_1)P_{18,\lambda}(X_2),\qquad
  \Omega_{\mathscr M}=\frac{\dd\sigma\wedge\dd X_1\wedge\dd X_2}{432vZ}.
\]
Use graph correspondences $[\Gamma_h]$ in source times target, so that their actions are $h_*$ and those of their transposes are $h^*$ \cite[Tags~0FFZ and 0FH2]{stacks-2018-7627557c}. For a rational map, these actions are computed on a resolution of its graph. They preserve holomorphic top forms, and on the dense finite \'etale locus the pushforward of such a form is its trace.

\begin{lemma}\label{sextic-modular-map}
  There is a smooth projective three-dimensional $\mathbb Q$-scheme $\overline W_{\mathbb Q}$ with a morphism $P_{\mathbb Q}\colon\overline W_{\mathbb Q}\to S_f$ and a dominant rational map $\alpha\colon\overline W_{\mathbb Q}\dashrightarrow\mathscr M$ over $\mathbb Q$. Its complex components compactify universal elliptic fibre squares over fine modular curves $U_j$. Put $d_j=\deg(U_j\to X_0(108))$ and $d=\sum_jd_j$, with degrees computed in projective modular groups. The form $\Omega_{\mathscr M}$ extends to a regular form and satisfies
  \[
    \alpha^*\Omega_{\mathscr M}=P_{\mathbb Q}^*\omega_f\ne0,\qquad
    \alpha_*P_{\mathbb Q}^*\omega_f=1944d\Omega_{\mathscr M}.
  \]
\end{lemma}

\begin{proof}
  Let $N$ be a congruence level for the fixed presentation of $S_f$, with defining subgroup containing $\Gamma(N)$. Set $n=108N$ and $L_n=\mathbb Q(\zeta_n)$, and take a smooth open $U$ of the geometrically connected component of $X(n)_{L_n}$ carrying the ordered basis $(1/n,\tau/n)$. Its universal elliptic curve $\mathcal E\to U$ gives $W=\mathcal E\times_U\mathcal E$ and the forgetful map $P\colon W\dashrightarrow(S_f)_{L_n}$ into the component represented by $U$ \cite[Sections~1.0--1.2]{scholl-1990-bb6eb16d}.

  Choose a nonzero rational invariant differential $\delta$ on $\mathcal E/U$ and write the fixed cusp section as $f_\delta\delta^4$. For $j\in\{6,18\}$ let $\phi_j$ be the quotient by $\angbr*{1/j}$ and use short coordinates $x_j,y_j$ with $\delta_j=\dd_{/U}x_j/(2y_j)$ and $\phi_j^*\delta_j=j\delta$. Put $A_\delta=-4x_6(\phi_6(1/18))$, using the full level structure. At the Tate cusp with its canonical differential $\delta_{\mathrm{can}}$, this order-three point on $\operatorname{Tate}(q^6)$ has parameter $\omega=e^{2\pi i/3}$. The short Tate coordinate \cite[eqs.~(A1.2.1), (A1.2.7)]{katz-1973-b9891d54} there has constant term $-1/4$ and coefficient of $q^{6m}$ equal to $\sum_{k\mid m}k(\omega^k+\omega^{-k}-2)=-3\sum_{k\mid m,3\nmid k}k$. Thus $A_{\delta_{\mathrm{can}}}=\mathscr A(6\tau)$ by \eqref{sextic-weight-two-base-form}. Set
  \[
    X_j=\frac{x_j}{A_\delta},\qquad Z=\frac{y_6y_{18}}{A_\delta^3},\qquad
    \sigma=\frac{f_\delta}{A_\delta^2}.
  \]
  Define $\lambda$ by the coefficient $-(1+8\lambda)/48$ of $X_6$ in the normalized short cubic, and put $v=1-2\lambda$. Under $\delta\mapsto c\delta$, $A_\delta$ scales like $x_6$, so these functions are independent of the differential frame.

  On this Tate chart the $6$-quotient has parameter $6\tau$ and normalized short coefficients $-E_4/(48\mathscr A^2)$ and $E_6/(864\mathscr A^3)$ by \cite[eq.~(A1.2.1)]{katz-1973-b9891d54}. The Eisenstein identities above identify its normalized cubic with $P_{6,\lambda}$ and give the analytic functions $\lambda,\sigma$ used in \zcref{sextic-hecke-base}.
  On the Tate chart, $(3\xi)(\mathbb Z+\mathbb ZW_3\xi)=\mathbb Z+\mathbb Z(3\xi)$. The short coefficients scale with weights four and six. Since $\mathscr A(W_3\xi)=-3\xi^2\mathscr A(\xi)$ and $\lambda(W_3\xi)=1-\lambda(\xi)$, the normalized coordinates on the $18$-quotient are those on the $6$-quotient at $W_3\xi$ scaled by $-1/3$ in $x$ and by $-1/27$ in $y^2$. Pulling back along the two isogenies gives
  \begin{equation}\label{sextic-cubic-theta-eisenstein}
    \frac{y_6^2}{A_\delta^3}=P_{6,\lambda}(X_6),\qquad
    \frac{y_{18}^2}{A_\delta^3}=P_{18,\lambda}(X_{18}).
  \end{equation}
  Their product defines $\alpha\colon W\dashrightarrow\mathscr M_{L_n}$. Over the common base, $\alpha^*X_1=X_6$, $\alpha^*X_2=X_{18}$ and $\alpha^*Z=Z$.
  Since $\lambda=27q^6+O(q^{12})$, the map $U\dashrightarrow\mathscr C$ is dominant and generically finite. The curve $\mathscr C$ is geometrically integral because $1-108\sigma^6$ has simple zeros. The radicand has a pole of order three at $X_1=\infty$, so $\mathscr M$ is geometrically integral. The functions $X_6,X_{18}$ are algebraically independent over $L_n(U)$, since each elliptic isogeny followed by its $x$-coordinate is finite. Thus $\alpha$ is dominant.

  Differentiating $27\sigma^6=\lambda(1-\lambda)$ and using \eqref{sextic-cubic-theta-differential} gives $\dd\sigma/v=\sigma\mathscr A(6\tau)\dd q/q$. The relative identity $\dd_{/U}X_6\wedge\dd_{/U}X_{18}/(4Z)=A_\delta\delta_6\wedge\delta_{18}$ on the two quotient curves pulls back under $\phi_6\times\phi_{18}$ with factor $108$. Wedging with $\dd\sigma$ removes the dependence on relative lifts. Thus
  \begin{equation}\label{sextic-direct-incidence-scholl-scale}
    \alpha^*\Omega_{\mathscr M}
    =\frac{\dd\sigma}{v}\wedge A_\delta\operatorname{pr}_1^*\delta\wedge\operatorname{pr}_2^*\delta
    =P^*\omega_f,
  \end{equation}
  where $\operatorname{pr}_i$ are the projections from $W$. The last equality uses the Tate expression for $\omega_f$ preceding \eqref{algebraic-modular-period-column}, with $A_{\delta_{\mathrm{can}}}=\mathscr A(6\tau)$ and $f=\sigma\mathscr A(6\tau)^2$. The product $A_\delta\operatorname{pr}_1^*\delta\wedge\operatorname{pr}_2^*\delta$ is independent of the frame.

  Resolve $P$ on a smooth projective compactification $\overline W/L_n$ and regard it as a $\mathbb Q$-scheme $\overline W_{\mathbb Q}$ with the resulting morphism $P_{\mathbb Q}$ to $S_f$. Smoothness and dimension are preserved under this restriction of the structure map because $L_n/\mathbb Q$ is finite separable. The map $\alpha$ is defined over $\mathbb Q$ on $\overline W_{\mathbb Q}$ and includes all its complex components. Its pullback of $\Omega_{\mathscr M}$ equals $P_{\mathbb Q}^*\omega_f$ by \eqref{sextic-direct-incidence-scholl-scale}, so $P_{\mathbb Q}^*\omega_f\ne0$ by dominance.

  On each $U_j$, the Tate chart identifies $\lambda$ with the degree-one coordinate of $X_0(3)$ in \zcref{sextic-hecke-base}. Conjugation by $\operatorname{diag}(6,1)$ takes $\Gamma_0(108)$ into $\Gamma_0(3)$ and preserves covolume, so the ratio $216/4$ of their projective indices gives $\deg(U_j\to X_0(3))=54d_j$. The right side of $27\sigma^6=\lambda(1-\lambda)$ has valuation $1$ at $\lambda=0$, giving a degree-six extension and $\deg(U_j\to\mathscr C)=9d_j$. Let $W_j$ be the universal fibre square over $U_j$. The two isogenies and the simultaneous-sign quotient give $\deg(\alpha|_{W_j})=6\cdot18\cdot2\cdot9d_j=1944d_j$.

  The trace along $\alpha$ gives the stated pushforward and shows that $\Omega_{\mathscr M}$ extends to a regular form.
\end{proof}

For $\Sigma=a+b+1$, let $\mathcal Q$ be a smooth projective model of
\[
  y^2=b\Sigma^2u^3-4a^4b^3,\qquad
  \Omega_{\mathcal Q}=\frac{\dd a\wedge\dd b\wedge\dd u}{2y}.
\]
The right-hand side of the defining equation has odd pole order at $u=\infty$, so $\mathcal Q$ is geometrically integral.

\begin{lemma}\label{sextic-quotient-map}
  There is a dominant rational map $\varpi\colon X_0\dashrightarrow\mathcal Q$ of degree two over $\mathbb Q$. The form $\Omega_{\mathcal Q}$ extends to a regular form and satisfies
  \[
    \varpi^*\Omega_{\mathcal Q}=\Omega_6,\qquad
    \varpi_*\Omega_6=2\Omega_{\mathcal Q}.
  \]
\end{lemma}

\begin{proof}
  On the quotient torus at $\psi=0$ take $a=t_0/t_1$, $b=t_2/t_1$, $u=a^2bt_1^2$ and $y=a^3b^2t_1^2(t_3-t_4)$. In \eqref{sextic-quotient-torus} we have $t_3+t_4=-t_1\Sigma$ and $t_3t_4=(a^2bt_1^4)^{-1}$, so the four functions satisfy the equation of $\mathcal Q$. Conversely, $t_0=at_1$, $t_2=bt_1$ and
  \[
    t_1^2=\frac{u}{a^2b},\qquad
    t_{3,4}=\frac{-t_1\Sigma\pm y/(a^3b^2t_1^2)}2.
  \]
  Thus the function field of $X_0$ is generated by $t_1$ over the field of the four coordinates. This subfield has transcendence degree three, so the coordinates define a dominant map $\varpi\colon X_0\dashrightarrow\mathcal Q$. The involution $(t_0,t_1,t_2,t_3,t_4)\mapsto(-t_0,-t_1,-t_2,-t_4,-t_3)$ fixes that field and negates $t_1$, so the extension is quadratic.

  In ambient torus coordinates $(a,b,t_1,t_3)$ with $t_4=(a^2bt_1^4t_3)^{-1}$, the logarithmic form and defining equation have
  \[
    \bigwedge_{j=0}^3\frac{\dd t_j}{t_j}
    =-\frac{\dd a\wedge\dd b\wedge\dd t_1\wedge\dd t_3}{abt_1t_3},
    \qquad
    \frac{\partial(t_1\Sigma+t_3+t_4)}{\partial t_3}
    =\frac{t_3-t_4}{t_3}.
  \]
  Moving $\dd t_3$ past three differentials gives the residue with the differential of the defining equation first. Since $\dd u\equiv2a^2bt_1\dd t_1$ modulo $\dd a,\dd b$,
  \[
    \Omega_6=\frac{\dd a\wedge\dd b\wedge\dd t_1}{abt_1(t_3-t_4)}
    =\varpi^*\Omega_{\mathcal Q}.
  \]
  Since $\varpi$ has degree two, its graph pushforward gives the stated trace and extends $\Omega_{\mathcal Q}$ to a regular form.
\end{proof}

Put $\kappa=\mathbb Q(\mathscr C)$ and $G_{2v}(r)=r^6+2vr^3+1$. Let $\mathscr H$ be the threefold over $\mathscr C$ with generic fibre
\[
  w^2=-3G_{2v}(r_1)G_{2v}(r_2).
\]
The radicand has valuation one at a simple zero of $G_{2v}(r_1)$ even over $\overline\kappa$, since $v^2\ne1$ makes $G_{2v}$ squarefree. Thus $\mathscr H$ is geometrically integral. Use a smooth projective model for $\mathscr H$.

\begin{lemma}\label{sextic-cover-maps}
  There are dominant rational maps $\pi_{\mathscr M}\colon\mathscr H\dashrightarrow\mathscr M$ and $\pi_{\mathcal Q}\colon\mathscr H\dashrightarrow\mathcal Q$ over $\mathbb Q$, of degrees $4$ and $72$, respectively. For $\mu=(r_1+1)(r_2-1)/r_2$, their pullbacks satisfy
  \begin{equation}\label{sextic-reverse-cover-form-ratio}
    \pi_{\mathcal Q}^*\Omega_{\mathcal Q}
    =\frac{648}{\mu}\pi_{\mathscr M}^*\Omega_{\mathscr M}.
  \end{equation}
\end{lemma}

\begin{proof}
  For the map to $\mathscr M$, assign the target coordinates by
  \begin{equation}
    \begin{aligned}
      X_1&=-\frac14-\frac{(1+v)r_1}{6(r_1-1)^2},&
      X_2&=\frac1{12}-\frac{(1-v)r_2}{18(r_2+1)^2},\\
      Z&=\frac{(1-v^2)w}{3888(r_1-1)^3(r_2+1)^3}.
    \end{aligned}
  \end{equation}
  Expansion of the two cubics gives
  \begin{equation}\label{sextic-direct-incidence-cubic-identities}
    P_{6,\lambda}(X_1)=-\frac{(1+v)^2G_{2v}(r_1)}{432(r_1-1)^6},\qquad
    P_{18,\lambda}(X_2)=\frac{(1-v)^2G_{2v}(r_2)}{11664(r_2+1)^6}.
  \end{equation}
  Their product gives $Z^2=P_{6,\lambda}(X_1)P_{18,\lambda}(X_2)$. The independent $r_i$ give independent $X_i$, so these coordinates define a dominant map $\pi_{\mathscr M}\colon\mathscr H\dashrightarrow\mathscr M$. Over a general target point, each $X_i$ gives a pair of distinct reciprocal roots $r_i,r_i^{-1}$. Each pair of choices recovers $w$ from $Z$, and the two cubic identities give the defining equation for $w$. Thus $\pi_{\mathscr M}$ has degree four. Relatively over $\mathscr C$,
  \[
    \dd X_1=\frac{(1+v)(r_1+1)}{6(r_1-1)^3}\dd r_1,\qquad
    \dd X_2=\frac{(1-v)(r_2-1)}{18(r_2+1)^3}\dd r_2.
  \]
  Substitution in $\Omega_{\mathscr M}$ gives
  \begin{equation}\label{sextic-direct-incidence-cover-scale}
    \pi_{\mathscr M}^*\Omega_{\mathscr M}
    =\frac{(r_1+1)(r_2-1)}{12}\frac{\dd\sigma\wedge\dd r_1\wedge\dd r_2}{vw}.
  \end{equation}

  For the map from $\mathscr H$ to the sextic quotient, set $R_i=r_i^3$. In coordinates $(E,F,x)$ on the quadric $x^2=EF-\beta$, the line $(E,F,x)=(1-hR_1,1-hR_2,v+h)$ passes through $(1,1,v)$. Substitution gives $h(hm-R_1-R_2-2v)=0$, where $m=R_1R_2-1$. On the common open where the denominators are nonzero, take the second intersection, $h=(R_1+R_2+2v)/m$. Then $mE=-G_{2v}(r_1)$ and $mF=-G_{2v}(r_2)$, so
  \begin{equation}\label{sextic-quotient-cover-identities}
    w^2=-3m^2EF,\qquad
    EF-\beta=(v+h)^2,\qquad
    \left(\frac{r_1}{r_2}\right)^3=\frac{1-E}{1-F}.
  \end{equation}
  The inverse formulas $R_1=(1-E)/h$ and $R_2=(1-F)/h$ identify $\kappa(R_1,R_2)$ with $\kappa(E,F,h)$, so $E,F$ are algebraically independent over $\kappa$. Put
  \begin{equation}\label{sextic-explicit-quotient-cover-map}
    \begin{aligned}
      a&=\frac{1-E}{E-\beta},&
      b&=-\frac{v^2F}{E-\beta F},\\
      u&=-\frac{abr_1}{3\sigma^2\Sigma r_2},&
      y&=-\frac{a^2bv(v+h)w}{9\sigma^3mE(1-F)},
    \end{aligned}
  \end{equation}
  Put $L=\beta(b+1)-1$. Thus $u$ and $y$ are multiples over $\kappa(E,F)$ of the roots $r_1/r_2$ and $(v+h)w/m$, respectively. From the formulas for $a,b$,
  \[
    \frac{1-E}{1-F}=-\frac{aL}{\Sigma},\qquad
    L=-\frac{v^2E}{E-\beta F},\qquad
    \Sigma=-\frac{L(1-F)}{E-\beta}.
  \]
  Cubing $u$ and squaring $y$ using \eqref{sextic-quotient-cover-identities} gives
  \begin{equation}\label{sextic-common-cover-target-identity}
    u^3=\frac{4a^4b^3L}{\beta\Sigma^4},\qquad
    y^2=\frac{4a^4b^3(bL-\beta\Sigma^2)}{\beta\Sigma^2}
    =b\Sigma^2u^3-4a^4b^3.
  \end{equation}
  Conversely, a general point of $\mathcal Q$ gives
  \[
    \beta=\frac{4a^4b^3}{4a^4b^3(b+1)-\Sigma^4u^3},\qquad
    E=\frac{1+\beta a}{1+a},\qquad F=\frac{bE}{L}.
  \]
  These birational changes identify $\mathbb Q(a,b,u^3)$ with $\mathbb Q(\beta,E,F)$. Choose
  \[
    \sigma^6=\beta/108,\qquad v^2=1-\beta,\qquad
    x^2=EF-\beta,\qquad r_2^3=\frac{1-F}{x-v}.
  \]
  Set $h=x-v$, $r_1=-3\sigma^2\Sigma ur_2/(ab)$, $m=r_1^3r_2^3-1$ and $w=-9\sigma^3mE(1-F)y/(a^2bvx)$.
  The cubic and square identities in \eqref{sextic-common-cover-target-identity} give $r_1^3=(1-E)/h$ and $w^2=-3m^2EF$; the quadric gives $hm=r_1^3+r_2^3+2v$, hence $mE=-G_{2v}(r_1)$ and $mF=-G_{2v}(r_2)$. Thus the recovered point lies on $\mathscr H$ and maps to the chosen target. The exceptional equations $h=0$ and $m=0$ would force $EF=1$ and $(E-F)^2+4\beta(1-E)(1-F)=0$, respectively, and do not hold generically. All $6\cdot2\cdot2\cdot3=72$ choices are separable and recovered from the source point. This proves dominance and degree $72$.

  At fixed $\sigma$, differentiation gives
  \[
    \dd a\wedge\dd b=-\frac{a\Sigma L}{E(1-E)(1-F)}\dd E\wedge\dd F,\qquad
    \dd E\wedge\dd F=-\frac{2h(v+h)}m\dd R_1\wedge\dd R_2.
  \]
  At fixed $a,b$, the cubic gives $\dd u/u=2\dd\sigma/(\sigma L)$. Substituting these identities and \eqref{sextic-explicit-quotient-cover-map} gives
  \begin{equation}\label{sextic-quotient-cover-pullback}
    \begin{aligned}
      \pi_{\mathcal Q}^*\Omega_{\mathcal Q}
      &=\dd\sigma\wedge\frac{u}{\sigma Ly}\dd a\wedge\dd b\\
      &=\frac{6}{vwr_1^2r_2}\dd\sigma\wedge\dd R_1\wedge\dd R_2\\
      &=54r_2\frac{\dd\sigma\wedge\dd r_1\wedge\dd r_2}{vw}.
    \end{aligned}
  \end{equation}
  Taking the ratio of \eqref{sextic-quotient-cover-pullback} and \eqref{sextic-direct-incidence-cover-scale} gives \eqref{sextic-reverse-cover-form-ratio}.
\end{proof}

For $P_{\mathbb Q}$ and $d$ from \zcref{sextic-modular-map}, normalize the cup pairing on $M_{f,B}$ by
\[
  Q_{\mathrm{mod}}(x,y)=\frac1d Q_{\overline W_{\mathbb Q}}(P_{\mathbb Q}^*x,P_{\mathbb Q}^*y).
\]

\begin{theorem}\label{sextic-normalized-motivic-comparison}\label{sextic-all-degree-projector}
  The plane $H_{X_0}^{\mathrm{att}}$ has rational dimension two. There is a self-adjoint algebraic correspondence $p_X$ over $\mathbb Q$ with rational coefficients whose Betti action is the cup-product orthogonal projector onto $H_{X_0}^{\mathrm{att}}$ in degree three and zero in every other degree. This action determines its homological class. For $M_{X_0}=(X_0,p_X,0)$ there is a unique isomorphism of homological motives
  \[
    \Phi_6\colon M_f\xrightarrow{\sim}M_{X_0}
    \qquad\text{with}\qquad
    (\Phi_6)_{\mathrm{dR}}(\omega_f)=\Omega_6.
  \]
  Its Betti realization satisfies
  \begin{equation}\label{sextic-cover-polarization}
    \Phi_6^*Q_{X_0}=6Q_{\mathrm{mod}}.
  \end{equation}
\end{theorem}

\begin{proof}
  Use the maps in \zcref{sextic-modular-map,sextic-quotient-map,sextic-cover-maps}. In $\mu^{-1}=r_2/((r_1+1)(r_2-1))$, the two values of each factor under its reciprocal involution sum to one. Hence $\operatorname{Tr}_{\pi_{\mathscr M}}(\mu^{-1})=1$. After extending constants to contain $\omega$, the map
  \[
    (\sigma,v,r_1,r_2,w)\longmapsto(\omega\sigma,v,\omega r_1,\omega^2r_2,w)
  \]
  fixes $r_1^3,r_2^3,v,w,\sigma^3$ and multiplies $r_1/r_2$ and $\sigma^2$ by $\omega^2$. It is therefore an automorphism over $\mathcal Q$ by \eqref{sextic-explicit-quotient-cover-map}. It multiplies each of the three nonconstant terms of $\mu=1+r_1-r_1/r_2-r_2^{-1}$ by a nontrivial cube root of unity. Their traces vanish, so $\operatorname{Tr}_{\pi_{\mathcal Q}}(\mu)=72$, also over $\mathbb Q$ by base change for trace. The form ratio \eqref{sextic-reverse-cover-form-ratio} gives
  \[
    (\pi_{\mathscr M})_*\pi_{\mathcal Q}^*\Omega_{\mathcal Q}=648\Omega_{\mathscr M},
    \qquad
    (\pi_{\mathcal Q})_*\pi_{\mathscr M}^*\Omega_{\mathscr M}=\Omega_{\mathcal Q}/9.
  \]
  Compose the graphs over $\mathbb Q$ to obtain
  \begin{equation}
    \mathscr G
    ={}^t[\Gamma_{\varpi}]\circ[\Gamma_{\pi_{\mathcal Q}}]
    \circ{}^t[\Gamma_{\pi_{\mathscr M}}]\circ[\Gamma_\alpha]\circ{}^t[\Gamma_{P_{\mathbb Q}}]
    \in\operatorname{CH}^3(S_f\times_{\mathbb Q}X_0)_{\mathbb Q}.
  \end{equation}
  The pushforwards for $\alpha$ and $\varpi$ give
  \begin{equation}
    (\mathscr G)_*\omega_f=216d\Omega_6,\qquad
    ({}^t\mathscr G)_*\Omega_6=1296(P_{\mathbb Q})_*P_{\mathbb Q}^*\omega_f.
  \end{equation}
  The first equality lets us apply \zcref{correspondence-induced-projector} with $e=e_f$ and $\Gamma=\mathscr G$. It gives $p_X$ and the isomorphism represented by $p_X\circ\mathscr G\circ e_f$. Since the projectors fix the chosen forms, the same equality gives
  \begin{equation}\label{sextic-form-normalized-cycle}
    \Phi_6=\frac1{216d}[p_X\circ\mathscr G\circ e_f]_{\mathrm{hom}}.
  \end{equation}
  Uniqueness follows from \eqref{absolute-holomorphic-restriction}. Identify the chosen forms with their Betti images. The cycles and $\Phi_6$ are rational, so the forward and reverse traces also hold on the conjugate forms. Adjointness gives
  \[
    \begin{aligned}
      216dQ_{X_0}(\Omega_6,\overline{\Omega}_6)
      &=Q_{S_f}(\omega_f,({}^t\mathscr G)_*\overline{\Omega}_6)\\
      &=1296dQ_{\mathrm{mod}}(\omega_f,\overline{\omega}_f).
    \end{aligned}
  \]
  Since $\omega_f,\overline{\omega}_f$ form a basis, this equality determines the alternating pairing and proves \eqref{sextic-cover-polarization}. In particular, $Q_{\mathrm{mod}}$ is nondegenerate.
\end{proof}

\subsection{Marked periods at the Fermat point}

Identify $M_{f,B}$ with $\mathbb Q^2$ in the basis dual to the critical modular symbols $(s_0,s_\chi)$. The coordinate map $\mathcal B_\gamma\colon H_B^3(X_0,\mathbb Q)\to\mathbb Q^4$ is fixed by the ambient Gamma classes of $\mathcal O_Y(kH)$, $0\leq k\leq3$, and continued along \eqref{sextic-endpoint-branch}. In the marked target coordinates $(\mathcal F_0,\mathcal F_1,z^0,z^1)$ the comparison matrix used below is
\begin{equation}\label{sextic-geometric-comparison-matrix}
  C=
  \begin{pmatrix} 2&0\\ -1&-3\\ 1&-1\\ 2&0
  \end{pmatrix}.
\end{equation}

Put $x=11664z$ and $\theta_x=x\dd/\dd x$. The periods of $\Omega_{\mathrm{MUM}}=-6\psi\Omega_X$ satisfy the Picard--Fuchs equation with operator
\begin{equation}\label{sextic-picard-fuchs-operator}
  \mathcal L=\theta_x^4-x
  \mleft(\theta_x+\frac16\mright)
  \mleft(\theta_x+\frac13\mright)
  \mleft(\theta_x+\frac23\mright)
  \mleft(\theta_x+\frac56\mright).
\end{equation}

The stack $\mathbb P(2,1,1,1,1)$ is Gorenstein Fano with trivial generic stabilizer, and $\mathcal O(6)$ is nef and convex. In the Laurent coordinates $(t_1,t_2,t_3,t_0)$, put
\[
  W=\frac{t_1+t_2+t_3+t_0+(t_1t_2t_3t_0^2)^{-1}}{6\psi}.
\]
The four coordinate monomials in $W$ have exponent vectors forming a basis of the torus character lattice, so the lattice-generation condition of the toric comparison holds. Its coefficient coordinate is $z=(6\psi)^{-6}$, and the odd coordinate permutation gives
\[
  \operatorname{Res}\frac{\dd t_1/t_1\wedge\dd t_2/t_2\wedge\dd t_3/t_3\wedge\dd t_0/t_0}{W-1}
  =-6\psi\Omega_X=\Omega_{\mathrm{MUM}}.
\]
The comparison of \cite[Theorems~5.7, 6.9--6.10]{iritani-2011-7b0e4739} identifies $K^0_{\mathrm{amb}}(Y)_{\mathbb Q}$ with $H_3(X_\psi,\mathbb Q)$ as rational flat spaces. Here the classes $[\mathcal O_Y(jH)]$, $0\leq j\leq3$, span the ambient $K$-space: their Chern characters span $1,H,H^2,H^3$. The structure-sheaf charge \eqref{toric-structure-sheaf-charge} therefore supplies the input to \zcref{constructed-line-bundle-extraction}. With the toric-stack mirror theorem \cite[Theorem~25 and Corollary~28]{coates-2019-fa6d0720}, it gives
\begin{equation}\label{sextic-local-mirror-identity}
  \sum_{d\geq1}\frac{dN_d^{\mathrm{GW}}}{3}q(z)^d=V(z).
\end{equation}
Specialize the coordinate functionals of \zcref{constructed-rational-period-frame} using \eqref{sextic-classical-constants}. Integration in this order defines $\mathcal B_\gamma$, continued along the chosen branch, and \eqref{single-charge-prepotential-frame} gives
\begin{equation}\label{sextic-iritani-normalization}
  \mathcal B_\gamma I_X(\Omega_{\mathrm{MUM}})
  =(2\pi i)^3\Pi_{\mathrm{raw}}.
\end{equation}

The cup matrix in this marking is $J_4$ by \zcref{toric-period-pairing}.

Put $\mathcal A=\{1/6,1/3,2/3,5/6\}$. For $a\in\mathcal A$, put $k_a(x)=(e^x-1)^{-a}/\mathrm B(a,1-a)$ for $x>0$ and $k_a(x)=0$ for $x\leq0$. Let $k=k_{1/6}*k_{1/3}*k_{2/3}*k_{5/6}$.

\begin{lemma}\label{sextic-euler-integral-period}
  Let $W(\rho,z)=z^\rho\sum_{n\ge0}\alpha_n(\rho)z^n$ be the generating function in \eqref{frobenius-solutions}. On the upper negative axis, with $z=-R+i0$ and $r=11664R$ for $R>0$, its continuation satisfies
  \begin{equation}\label{sextic-euler-period}
    W(\rho,z)=z^\rho\frac{I_r(\rho)}{I_0(\rho)},\qquad
    I_r(\rho)=\int_0^\infty\frac{e^{-\rho x}k(x)}{1+re^{-x}}\dd x,
    \qquad \abs*{\rho}<1/12,
  \end{equation}
  with $\log z=\log R+i\pi$. In particular $w_0(z)=I_r(0)>0$, and
  \[
    A\coloneqq\lim_{R\to\infty}R^{1/6}w_0(-R+i0)
    =11664^{-1/6}\prod_{b\in\mathcal A\setminus\{1/6\}}
    \int_0^\infty e^{x/6}k_b(x)\dd x>0.
  \]
  Locally uniformly for $\abs*{\rho}<1/12$,
  \begin{equation}\label{sextic-euler-endpoint-limit}
    R^{1/6}W(\rho,-R+i0)\longrightarrow
    A\frac{e^{i\pi\rho}\sin(\pi/6)}{\sin\pi(\rho+1/6)}
    \frac{\Gamma(1+2\rho)\Gamma(1+\rho)^4}{\Gamma(1+6\rho)}.
  \end{equation}
\end{lemma}
\begin{proof}
  With $u=e^{-x}$, the beta integral \cite[eq.~(5.12.1)]{dlmf-2010-fa11482d} gives
  \[
    \begin{aligned}
      \int_0^\infty e^{-\rho x}k_a(x)\dd x
      &=\frac1{\mathrm B(a,1-a)}\int_0^1u^{a+\rho-1}(1-u)^{-a}\dd u\\
      &=\frac{\Gamma(a+\rho)}{\Gamma(a)\Gamma(1+\rho)},\qquad \rmRe\rho>-a.
    \end{aligned}
  \]
  Thus $\int k=1$ and $\int_0^\infty e^{cx}k(x)\dd x<\infty$ for every real $c<1/6$. The product of these transforms gives $I_0(\rho)\ne0$ and
  \begin{align*}
    \frac{I_0(n+\rho)}{I_0(\rho)}
    &=\frac{\prod_{a\in\mathcal A}(a+\rho)_n}{(1+\rho)_n^4}\\
    &=11664^{-n}\alpha_n(\rho)
  \end{align*}
  by \eqref{constructed-frobenius-recurrence}, with the same integral for $I_0$ on $\rmRe\rho>-1/6$. Expanding the kernel geometrically for $r<1$ gives \eqref{sextic-euler-period}, which continues to $r>0$ by holomorphy on $\rmRe r>0$. These integrability bounds justify parameter differentiation near zero.

  To take the limit, put $a=1/6$ and fix $0<\epsilon<1/12$. For $t>0$, use the real logarithm for positive powers and substitute $v=te^{-x}$. The beta and reflection formulas \cite[eqs.~(5.12.3), (5.5.3)]{dlmf-2010-fa11482d} give
  \[
    \begin{aligned}
      t^{a+\rho}\int_0^\infty\frac{e^{-\rho x}k_a(x)}{1+te^{-x}}\dd x
      &=\frac1{\mathrm B(a,1-a)}\int_0^t\frac{v^{a+\rho-1}(1-v/t)^{-a}}{1+v}\dd v\\
      &\longrightarrow\frac{\mathrm B(a+\rho,1-a-\rho)}{\mathrm B(a,1-a)}
      =\frac{\sin\pi a}{\sin\pi(a+\rho)}
      \qquad(t\to\infty).
    \end{aligned}
  \]
  The limit is uniform for $\abs*{\rho}\le\epsilon$. On $0<v\le t/2$, the factor $(1-v/t)^{-a}$ is at most $2^a$, and the integrand is bounded in modulus by $O(v^{a-\epsilon-1})$ for $v\le1$ and $O(v^{a+\epsilon-2})$ for $v\ge1$. These common integrable bounds and convergence on compact intervals give uniform convergence there. For $t\ge2$, scaling $v$ by $t$ bounds the terminal part by $O(t^{a+\epsilon-1})$, since $a<1$, and this tends to zero uniformly. For $0<t\le1$, the scaled integral is bounded by $\int_0^\infty e^{\epsilon x}k_a(x)\dd x<\infty$, since $\epsilon<a$ and $t^{a-\epsilon}\le1$. Together with the direct bound for $1\le t\le2$, this gives a uniform bound for every $t>0$.

  Write the convolution as a fourfold integral and put $U=\exp(-\sum_{b\in\mathcal A\setminus\{a\}}x_b)$, so $t=rU$ in the $x_a$ integral. Since $r^{a+\rho}U^\rho t^{-a-\rho}=U^{-a}$, the remaining three integrals are dominated, up to a constant independent of $r$ and $\rho$, by $U^{-a}\prod_{b\ne a}k_b(x_b)$, whose integral is
  \[
    \int_{(0,\infty)^3}U^{-a}\prod_{b\ne a}k_b(x_b)\prod_{b\ne a}\dd x_b
    =\prod_{b\ne a}\int_0^\infty e^{ax}k_b(x)\dd x<\infty,
  \]
  since every remaining $b$ is greater than $a$. Dominated convergence, applied to the supremum of the difference over the parameter disc, gives the formula for $A$ and the corresponding locally uniform limit. Using $z^\rho=R^\rho e^{i\pi\rho}$ and the multiplication formula \cite[eq.~(5.5.6)]{dlmf-2010-fa11482d}, we have
  \[
    I_0(\rho)=11664^{-\rho}
    \frac{\Gamma(1+6\rho)}{\Gamma(1+2\rho)\Gamma(1+\rho)^4},
  \]
  proving \eqref{sextic-euler-endpoint-limit}.
\end{proof}

On a sufficiently small disc $\Delta$ about the smooth fibre $\psi=0$, the fixed toric resolution gives a smooth projective family and $\Omega_X(\psi)$ extends as a holomorphic relative form with nonzero value $\Omega_6$. Identify its rational Betti local system with the fibre at $0$ by parallel transport on $\Delta$, agreeing with the upper continuation in \eqref{sextic-endpoint-branch}, and write $\Pi_\Omega(\psi)=\mathcal B_\gamma I_{X_\psi}(\Omega_X(\psi))$ in this trivialization.

\begin{proposition}\label{sextic-marked-plane}
  The rational plane in \zcref{sextic-normalized-motivic-comparison} has complexification $H^{3,0}(X_0)\oplus H^{0,3}(X_0)$. The marking continued along the upper negative axis \eqref{sextic-endpoint-branch} satisfies $\mathcal B_\gamma(H_{X_0}^{\mathrm{att}})=\operatorname{im}C$. Then
  \begin{equation}\label{sextic-residue-principal-period}
    \Pi_\Omega(0)=-4\pi^3A
    \begin{pmatrix}
      2\\-1-3i\sqrt3\\1-i\sqrt3\\2
    \end{pmatrix}.
  \end{equation}
  The coordinate $(\Pi_\Omega)_3(0)$ is nonzero, and the normalized column extends holomorphically to $\psi=0$, with value
  \begin{equation}\label{sextic-spectral-endpoint-column}
    \left.\frac{\Pi_{\mathrm{raw}}(z(\psi))}{z^0(z(\psi))}\right|_{\psi=0}
    =\frac{\Pi_\Omega(0)}{(\Pi_\Omega)_3(0)}
    =\frac{C(1,i\sqrt3)^T}{1-i\sqrt3}
    =
    \begin{pmatrix}
      t_*\\3-2t_*\\1\\t_*
    \end{pmatrix},
    \qquad t_*=\frac{1+i\sqrt3}{2}.
  \end{equation}
\end{proposition}

\begin{proof}
  Local uniformity in \eqref{sextic-euler-endpoint-limit} permits coefficient extraction through $\rho^3$. The log-Gamma expansion \cite[eq.~(5.7.3)]{dlmf-2010-fa11482d}, with $p_j=w_j/(2\pi i)^j$, gives
  \begin{align*}
    \frac{p_1}{p_0}&\longrightarrow\frac{1+i\sqrt3}{2},\\
    \frac{p_2}{p_0}&\longrightarrow-\frac16+\frac{i\sqrt3}{4},\\
    \frac{p_3}{p_0}&\longrightarrow-\frac{1+i\sqrt3}{8}
    +\frac{17i\zeta(3)}{2\pi^3}.
  \end{align*}
  Substitution in \eqref{constructed-prepotential-periods}, with the constants \eqref{sextic-classical-constants}, gives the ray limit of the normalized column in \eqref{sextic-spectral-endpoint-column}.

  The marking \eqref{sextic-iritani-normalization} and residue gauge \eqref{sextic-mum-residue} give $\Pi_{\mathrm{raw}}(z(\psi))=-6\psi\Pi_\Omega(\psi)/(2\pi i)^3$. On the chosen branch,
  \begin{equation}\label{sextic-residue-period-phase}
    (\Pi_\Omega)_3(\psi)=-\frac{(2\pi i)^3}{6\psi}w_0(-R+i0)
    =-4\pi^3(1-i\sqrt3)R^{1/6}w_0(-R+i0).
  \end{equation}
  Hence $(\Pi_\Omega)_3(0)=-4\pi^3(1-i\sqrt3)A\ne0$. The holomorphic extension of $\Pi_\Omega$ makes $\Pi_\Omega/(\Pi_\Omega)_3$ holomorphic at zero, proving \eqref{sextic-spectral-endpoint-column}. Multiplying by the third coordinate gives \eqref{sextic-residue-principal-period}. Its vector and conjugate span $(\operatorname{im}C)_{\mathbb C}$, which is therefore the marked extreme Hodge plane.
\end{proof}

\subsection{Geometric calibration and special-value identities}

\begin{proposition}
  In the basis dual to $(s_0,s_\chi)$, the form $\omega_f$ and the normalized cup pairing satisfy
  \begin{gather}
    I_f\omega_f=s_0(I_f\omega_f)(1,i\sqrt3)^T,\qquad s_0(I_f\omega_f)<0,\label{sextic-modular-cm-coordinates}\\
    [Q_{\mathrm{mod}}]_{(s_0,s_\chi)^\vee}=486J_2,\qquad
    J_2=
    \begin{pmatrix}0&1\\-1&0
    \end{pmatrix}.
    \label{sextic-modular-cup-matrix}
  \end{gather}
\end{proposition}

\begin{proof}
  Use $g(\xi)=f(\xi/6)$ and $T\xi=\xi+1$ as above. Identify quadratics with coefficient columns in order $1,X,X^2$, and define
  \[
    \rho(T)Q(X)=Q(X-1),\qquad
    \rho(W_3)Q(X)=3X^2Q\mleft(-\frac1{3X}\mright),\qquad
    \mathsf E=\zeta\rho(W_3)\rho(T).
  \]
  The right side defining $\rho(W_3)$ expands to a polynomial. Set $\eta=g(\xi)(\xi-X)^2\dd \xi$, $F(\xi)=\int_{i\infty}^\xi\eta$ and $P=\int_0^{i\infty}\eta$. Cuspidality gives convergence at rational cusps, and holomorphy makes the finite parts path independent. Substitution using \eqref{sextic-hecke-transformations} gives
  \[
    T^*\eta=\zeta\rho(T)\eta,\qquad
    W_3^*\eta=\frac13g(\xi)(1+3\xi X)^2\dd \xi=\rho(W_3)\eta.
  \]
  Since $T\infty=\infty$ and $W_3\infty=0$, the primitives satisfy
  \begin{align*}
    F(T\xi)&=\zeta\rho(T)F(\xi),\\
    F(W_3\xi)&=\rho(W_3)F(\xi)-P,\\
    F(W_3T\xi)&=\mathsf E F(\xi)-P.
  \end{align*}
  The relation $W_3^2=1$ gives $\rho(W_3)P=-P$. Put $\xi_L=(-3+i\sqrt3)/6$. Since $W_3T\xi_L=\xi_L$ and $3\zeta \xi_L^2=1$, one has $(\mathsf E Q)(\xi_L)=Q(\xi_L)$ for every quadratic $Q$. Evaluating the fixed-point primitive equation therefore gives $P(\xi_L)=0$. Together these two relations give
  \[
    P=P(0)(1,i\sqrt3,-3)^T,\qquad
    P(0)=\int_0^{i\infty}\xi^2g(\xi)\dd \xi.
  \]
  Three iterations of the primitive equation, with $(W_3T)^3\infty=-1/2$ and $F(i\infty)=0$, give $F(-1/2)=-(\mathrm{id}+\mathsf E+\mathsf E^2)P$. For $t\in\mathbb Q$, put $I(t)=\int_t^{i\infty}g(\xi)(\xi-t)^2\dd \xi=-F(t)|_{X=t}$. Evaluation and the translation law $I(t+1)=\zeta I(t)$ give
  \begin{align*}
    I(-1/2)&=\mleft((\mathrm{id}+\mathsf E+\mathsf E^2)P\mright)(-1/2)\\
    &=\frac{3-i\sqrt3}{4}P(0),\\
    I(3/2)-I(9/2)&=2\zeta^2I(-1/2)\\
    &=i\sqrt3P(0).
  \end{align*}
  Substituting $\xi=6\tau$ in the symbols gives $\Phi_f(s_0)=P(0)/216$ and $\Phi_f(s_\chi)=i\sqrt3P(0)/216$. The Tate normalization gives $s_0(I_f\omega_f)=(2\pi i)^3P(0)/216=-8\pi^3\int_0^\infty y^2f(iy)\dd y<0$ by \eqref{sextic-hecke-source-positivity}.

  For quadratic columns use $\angbr*{a,b}=a_0b_2-a_1b_1/2+a_2b_0$, preserved by $\rho(T)$ and $\rho(W_3)$. Set $\Theta=\angbr*{F,\overline{\eta}}$. For $\xi=x+iy$, the primitive identities give
  \[
    \dd\Theta=8i\abs*{g(\xi)}^2y^2\dd x\wedge\dd y,\quad
    T^*\Theta=\Theta,\quad
    W_3^*\Theta=\Theta+\angbr*{P,\overline{\eta}}.
  \]
  The second identity uses $\zeta\overline{\zeta}=1$, and the third uses $\rho(W_3)P=-P$.

  The group $\angbr*{T,W_3}$ is the index-two Fricke extension of $\Gamma_0(3)$, since $W_3$ normalizes $\Gamma_0(3)$ and $T,W_3T^{-1}W_3^{-1}$ generate it modulo its centre. Its fundamental domain is
  \[
    D=\setc*{\xi\in\mathfrak H}{\abs*{\rmRe\xi}\le1/2,\ \abs*{\xi}\ge1/\sqrt3}.
  \]
  Indeed, the heights in an orbit have denominators $\abs*{c\xi+d}^2$ or $3\abs*{a\xi+b}^2$ for integer pairs, so a maximum exists. Translate a maximum-height point into the strip. It lies on or outside the circle, since otherwise $W_3$ increases its height. The area of $D$ is $2\pi/3$, the coarea of the group, proving the claim.

  Use $T\xi_L=W_3\xi_L$ and put $B=P(0)(-1,-\zeta,0)^T$. The two polynomial actions give $(\mathsf E-1)B=P$ and $(\rho(W_3)-1)B-P=P(0)(0,1,0)^T$. Thus $h=F(\xi_L)-B$ is fixed by $\mathsf E$. Invariance of the Hermitian pairing gives $\angbr*{P,\overline{h}}=0$, and $\rho(W_3)P=-P$ gives $\angbr*{P,\overline{(\rho(W_3)-1)h}}=0$. The vertical sides cancel by translation and the top tends to zero by cuspidality. The Fricke involution reverses the lower arc and negates $\angbr*{P,\overline{\eta}}$. Pairing its two halves in Stokes' formula therefore gives
  \[
    16i\int_D\abs*{g(\xi)}^2y^2\dd x\dd y
    =-\angbr*{P,\overline{F(T\xi_L)-F(\xi_L)}}
    =\frac{i\sqrt3}{2}\abs*{P(0)}^2.
  \]
  The conjugate of $\Gamma_0(108)$ by $\operatorname{diag}(6,1)$ lies in the Fricke group and has index $108$ there by covolume. The change $\xi=6\tau$ contributes $6^{-4}$ to the Petersson density in the $\xi$-integral. On the fibre square, $\omega_f=(2\pi i)^3f(\tau)\dd\tau\wedge\dd z_1\wedge\dd z_2$. The fibre integrals $\int\dd z_j\wedge\dd\overline{z}_j=-2i\rmIm\tau$ and the sign from reordering the six forms into conjugate pairs give, after division by the fine-level degree,
  \[
    Q_{\mathrm{mod}}(\omega_f,\overline{\omega}_f)
    =-512i\pi^6\frac{108}{6^4}\int_D\abs*{g}^2y^2\dd x\dd y
    =-972i\sqrt3\abs*{s_0(I_f\omega_f)}^2.
  \]
  If the cup matrix is $cJ_2$, the column in \eqref{sextic-modular-cm-coordinates} gives $-2ci\sqrt3\abs*{s_0(I_f\omega_f)}^2$ for the same pairing. Hence $c=486$.
\end{proof}

\begin{proposition}\label{sextic-absolute-de-rham-calibration}\label{sextic-algebraic-comparison}\label{sextic-algebraic-comparison-action}
  For the upper continuation \eqref{sextic-endpoint-branch}, the marked Betti matrix of the form-normalized isomorphism \eqref{sextic-form-normalized-cycle} is $\mathcal B_\gamma(\Phi_6)_B=27C$. Its period column is
  \begin{equation}\label{sextic-regular-endpoint-pairing}
    \Pi_\Omega(0)=-54C
    \begin{pmatrix}L(f,3)\\2iL(f,\chi_{-4},3)
    \end{pmatrix}.
  \end{equation}
\end{proposition}

\begin{proof}
  The two holomorphic columns in \eqref{sextic-modular-cm-coordinates} and \eqref{sextic-residue-principal-period} give
  \[
    \mathcal B_\gamma(\Phi_6)_B(1,i\sqrt3)^T
    =-\frac{4\pi^3A}{s_0(I_f\omega_f)}C(1,i\sqrt3)^T.
  \]
  Comparing real and imaginary parts shows that the marked Betti matrix of $\Phi_6$ is a positive multiple of $C$. The target cup matrix is $J_4$ by \zcref{toric-period-pairing}. \zcref{sextic-cover-polarization,sextic-modular-cup-matrix}, together with direct multiplication, give
  \[
    \begin{aligned}
      \mleft(\mathcal B_\gamma(\Phi_6)_B\mright)^T J_4\mathcal B_\gamma(\Phi_6)_B
      &=6\cdot486J_2,\\
      C^TJ_4C&=4J_2.
    \end{aligned}
  \]
  Positivity of the multiple therefore gives
  \[
    \mathcal B_\gamma(\Phi_6)_B
    =\sqrt{\frac{6\cdot486}{4}}C=27C.
  \]
  Apply \eqref{absolute-chosen-form-equality} to $\Phi_6$ and use \eqref{algebraic-modular-period-column} to obtain \eqref{sextic-regular-endpoint-pairing}.
\end{proof}

\begin{corollary}\label{sextic-critical-ratio}\label{sextic-critical-gamma-evaluation}
  For the fixed modular form $f$ of \eqref{sextic-hecke-theta-series},
  \begin{equation}\label{sextic-absolute-gamma-values}
    L(f,3)=\frac{\sqrt3\Gamma(1/3)^9}{288\cdot2^{1/3}\pi^3},\qquad
    L(f,\chi_{-4},3)=\frac{\Gamma(1/3)^9}{3\cdot2^{19/3}\pi^3}.
  \end{equation}
  In particular, $L(f,\chi_{-4},3)/L(f,3)=\sqrt3/2$.
\end{corollary}
\begin{proof}
  \zcref{sextic-regular-endpoint-pairing,sextic-residue-principal-period}, together with the column rank of $C$, give $L(f,3)=2\pi^3 A/27$ and $L(f,\chi_{-4},3)=(\sqrt3/2)L(f,3)$. The beta integrals in \zcref{sextic-euler-integral-period} and the reflection and duplication formulas \cite[eqs.~(5.5.3), (5.5.5)]{dlmf-2010-fa11482d} give
  \[
    A=11664^{-1/6}\frac{\Gamma(1/6)\sqrt\pi}{\Gamma(1/3)\Gamma(5/6)^4}
    =\frac{3\sqrt3\Gamma(1/3)^9}{64\cdot2^{1/3}\pi^6}.
  \]
  Substitution gives both values.
\end{proof}

Substituting the second entry of \eqref{sextic-spectral-endpoint-column} into \eqref{constructed-prepotential-periods} gives $3-2t_*=7/4-t_*/2-3t_*^2/2+3V_*/(4\pi^2)$. Since $t_*^2=t_*-1$, this gives $V_*=-\pi^2/3$. The ratio of the fourth and third entries in \eqref{sextic-regular-endpoint-pairing}, together with \zcref{sextic-critical-ratio}, gives
\begin{equation}\label{sextic-genus-zero-endpoint}
  t_*=\frac12+i\frac{L(f,\chi_{-4},3)}{L(f,3)},\qquad
  q_*=-e^{-\pi\sqrt3}.
\end{equation}

For \zcref{geometric-partition-identity}, take the projector of \zcref{sextic-all-degree-projector} and $\Phi_6$ of \eqref{sextic-form-normalized-cycle}, whose marked Betti matrix is $27C$. The modular form is the fixed $f$ by the motivic comparison and uniqueness in \zcref{intrinsic-modular-realization}. The row relation $(0,1,-3,2)C=0$ and \eqref{sextic-classical-constants} give $(\Lambda,\Upsilon)=(3,-2)$ and $(x_0,y_0^2)=(1/2,7/12)$. The local mirror identities and the marked endpoint supply the period data. By \eqref{sextic-radicand-bound}, $\mathfrak w^2=-7/12+V/(2\pi^2)$ stays negative on the upper negative axis, and $w_0$ stays positive by \zcref{sextic-euler-integral-period}. Hence the branch starting at $i\sqrt{7/12}$ continues with endpoint $i\sqrt3/2$. Specializing \eqref{normalized-scalar-partition-expansion} gives
\begin{equation}\label{sextic-partition-sum}
  T(s)=\frac12+\frac{i\sqrt{21}}6-\frac{i}{\sqrt{6\pi^3}}
  \sum_{m\geq1}(-1)^ms^m
  \sum_{\mathfrak p\vdash m}\widetilde N_{\mathfrak p}
  \left(\frac{m}{\pi\sqrt{21}}\right)^{l(\mathfrak p)-1/2}
  K_{l(\mathfrak p)-1/2}\mleft(\frac{\pi m\sqrt{21}}3\mright).
\end{equation}
The Taylor series diverges at $s=1$ by \zcref{sextic-taylor-divergence}.

\begin{theorem}\label{sextic-partition-comparison}
  Let $Y\subset\mathbb P(2,1,1,1,1)$ be a smooth sextic missing the stacky point, with degrees measured by $H=c_1(\mathcal O_Y(1))$, and let $T$ be its local solution given by \eqref{sextic-partition-sum}. For the newform $f$ with label $\mathtt{108.4.a.c}$, continuation along the inversion path induced by the upper negative-axis period branch $z=-R+i0$, $R\to\infty$, gives
  \begin{equation}\label{sextic-continued-partition-identity}
    6\frac{L(f,\chi_{-4},3)}{L(f,3)}
    =-6i\left(T_\gamma(1)-\frac12\right).
  \end{equation}
  The right side is the finite value obtained by continuing the function defined near $s=0$ by the partition series.
\end{theorem}

\begin{proof}
  By \zcref{geometric-partition-identity}, $T_\gamma(1)=t_*$. Substituting \eqref{sextic-genus-zero-endpoint} gives the identity.
\end{proof}

\section{Conclusion}\label{conclusion}

A correspondence with nonzero action between the two rank-two realizations, together with their rational holomorphic generators, fixes the period comparison. For the Fermat sextic, the two traces and the cup pairings determine the marked Betti matrix and hence the two critical-value formulas. Dividing by the principal period gives the mirror coordinate in the continued Gromov--Witten identity.

Constructing explicit correspondences for further cases of the identities in \cite[Section~5]{candelas-2025-ad265528} remains open. For the bicubic (AESZ4), the task is to construct, without assuming the Tate conjecture, a cycle over $\mathbb Q$ with nonzero action from the Scholl realization of its weight-four newform to the rank-two summand of the mirror at the attractor point, and to determine its action in the marked period frame.

\section*{Declarations}

\subsection*{Funding}
The authors were partially supported by Taiwan NSTC grant 115-2115-M-002-010-MY2 and NTU Core Consortiums grant 115L890503 (TIMS).

\subsection*{Author contributions}
Both authors contributed to and approved the manuscript.

\subsection*{Data availability}
This study generated and analysed no datasets.

\subsection*{Competing interests}
The authors have no relevant financial or non-financial interests to disclose.

\subsection*{Statement on AI use}
During the preparation of this manuscript, the authors used GPT-6 to discuss proof strategies, perform certain computations, and improve the LaTeX editing. The original ideas, main results, and overall mathematical framework of this work were developed by the authors. All AI-assisted suggestions and outputs were evaluated by the authors, who take full responsibility for the correctness of the manuscript.

\appendix

\section{Analytic estimates for the Fermat sextic}\label{sextic-analytic-appendix}

\subsection{Divergence of the sextic Taylor series}

\begin{proposition}\label{sextic-taylor-divergence}
  Let $T$ be the sextic solution in \eqref{sextic-partition-sum}. At $s=1$ its Taylor terms do not tend to zero. The same holds for the terms of the partition expansion grouped by total degree, so both series diverge there.
\end{proposition}

\begin{proof}
  Work towards the positive real conifold in $x=11664z$, writing $h_j$ at $z=x/11664$. Let $c_n=11664^{-n}\alpha_n(0)$, $d_n=(\log\alpha_n)'(0)$ and $e_n=(\log\alpha_n)''(0)$. Then $h_0=\sum c_nx^n$, $h_1=\sum c_nd_nx^n$ and $2h_2=\sum c_n(d_n^2+e_n)x^n$. The Frobenius product in the proof of \zcref{sextic-euler-integral-period} gives
  \[
    d_n=\sum_{a\in\mathcal A}\sum_{j=0}^{n-1}
    \mleft(\frac1{j+a}-\frac1{j+1}\mright),\qquad
    -e_n=\sum_{a\in\mathcal A}\sum_{j=0}^{n-1}
    \mleft(\frac1{(j+a)^2}-\frac1{(j+1)^2}\mright).
  \]
  Since $j+a\ge a(j+1)$, the first difference is at most $(a^{-1}-1)/(j+1)^2$ and the second at most $(a^{-2}-1)/(j+1)^2$. Summing over $a$ gives the bounds $0\le d_n\le(77/10)\zeta(2)<13$ and $0\le-e_n\le(4469/100)\zeta(2)<75$. For $n\ge1$, H\"older in the Gamma integral gives $\Gamma(n+a)\le\Gamma(n+1)n^{a-1}$. The sum of the four exponents is $2$ and reflection gives $\prod_{a\in\mathcal A}\Gamma(a)=4\pi^2/\sqrt3$ \cite[eqs.~(5.2.1), (5.5.3)]{dlmf-2010-fa11482d}. Hence $c_n\le\sqrt3/(4\pi^2n^2)$ and $\sum_{n\ge1}c_n<1/12$.

  The series converge uniformly on $[0,1]$, with $h_0\ge1$, $0\le h_1/h_0<13/12$ and $-25/8<h_2/h_0<169/24$. Thus $-4<V(x/11664)<8$ by $V=h_2/h_0-(h_1/h_0)^2/2$. Using $\pi>3$, the functions
  \begin{equation}
    y(x)=\sqrt{\frac7{12}-\frac{V(x/11664)}{2\pi^2}},
    \qquad s(x)=-q(x)e^{2\pi y(x)}
  \end{equation}
  satisfy $0<y(x)<1$. Also $q(x)=(x/11664)e^{h_1/h_0}\le e^{13/12}/11664$. Since $13/12+2\pi<8$ and $e^8<11664$, the path satisfies $\abs*{s(x)}<1$ on $[0,1]$.

  In the divided-power Frobenius basis, \cite[Theorem~1]{chen-2008-c6669bf5} gives the conifold monodromy. Conjugating by the raw-frame change of basis in \eqref{constructed-prepotential-periods} yields
  \begin{equation}\label{sextic-conifold-monodromy}
    (\mathcal F_0,\mathcal F_1,z^0,z^1)^T
    \longmapsto(\mathcal F_0,\mathcal F_1,z^0-\mathcal F_0,z^1)^T.
  \end{equation}

  Since $0<q(1)<1$ and $h_0$ has positive limit, the limit of $z^1=z^0t$ is nonzero. The operator has a regular singularity at $x=1$ \cite[eqs.~(16.8.3), (16.8.5)]{dlmf-2010-fa11482d}. Substituting $(x-1)^r$ gives indicial polynomial $r(r-1)^2(r-2)$ up to a nonzero scalar. The solutions $z^1,\mathcal F_1$ fixed by \eqref{sextic-conifold-monodromy} are therefore holomorphic at $1$: a single-valued regular-singular solution has no logarithmic terms and its leading exponent is nonnegative. Thus $h=\mathcal F_1/z^1$ is holomorphic there and nonconstant by independence of the two solutions.

  Suppose that $[s^m]T\to0$. Boundedness of the coefficients makes $T$ holomorphic on the unit disc, which contains the path $s(x)$. Initial implicit uniqueness and the real-analytic identity theorem give $T(s(x))=1/2+iy(x)$ on $0<x<1$. Near $s(1)$ define
  \begin{align*}
    \widehat t(s)&=T(s)-\frac12+\frac{\log(-s)}{2\pi i},\\
    \widehat h(s)&=-\frac12+\frac3{2\widehat t(s)}
    \left((T(s)-1/2)^2+\frac74-\widehat t(s)^2\right),
  \end{align*}
  choosing $\log(-s(1))$ real. Since $\widehat t(s(1))=\log q(1)/(2\pi i)\ne0$, both functions are holomorphic there and the period frame \eqref{constructed-prepotential-periods} gives $\widehat t(s(x))=t(x)$ and $\widehat h(s(x))=h(x)$. Since $h$ is nonconstant, $\widehat h(s)-h(1)=\varphi(s)^e$ with $e\ge1$ and a local coordinate $\varphi$ at $s(1)$. On the branch from $0<x<1$ this gives $s(x)=\varphi^{-1}((h(x)-h(1))^{1/e})$. After $e$ loops about $x=1$ this branch and $t=\widehat t(s)$ return to their original values. But \eqref{sextic-conifold-monodromy} gives
  \[
    t\longmapsto\frac{z^1}{z^0-e\mathcal F_0}\ne\frac{z^1}{z^0}=t,
  \]
  since neither $z^1$ nor $\mathcal F_0$ vanishes identically. This contradiction proves that the Taylor terms at $s=1$ do not tend to zero. The total-degree $m$ term of \eqref{sextic-partition-sum} at $s=1$ equals $[s^m]T$, so the grouped partition series diverges as well.
\end{proof}

\subsection{The sextic radicand estimate}

\begin{proposition}\label{sextic-global-radicand-bound}
  For the sextic operator \eqref{sextic-picard-fuchs-operator}, on the upper negative-axis branch and for every $R>0$,
  \begin{equation}\label{sextic-radicand-bound}
    -\frac7{12}+\frac{V(-R+i0)}{2\pi^2}<-\frac1{12}.
  \end{equation}
\end{proposition}
\begin{proof}
  Write $z=-R+i0$ and $r=11664R$. It suffices to prove $V(z)<\pi^2$. In \eqref{sextic-euler-period}, $I_r(0)>0$ and the exponential moments justify the first two parameter derivatives. Since $w_j=[\rho^j]W$ and $z^\rho$ contributes only a linear term to $\log W$,
  \begin{equation}\label{sextic-log-derivative-radicand}
    2V(z)=(\log I_r)''(0)-(\log I_0)''(0).
  \end{equation}

  For $a\in\{1/6,1/3\}$, put $\phi_a=k_a*k_{1-a}$. When $v>0$, the substitutions $u=(e^x-1)/(e^v-1)$ and $u=(1+e^y)^{-1}$ give
  \begin{align*}
    \phi_a(v)
    &=\mathrm B(a,1-a)^{-2}\int_0^v
    (e^x-1)^{-a}(e^{v-x}-1)^{a-1}\dd x\\
    &=\mathrm B(a,1-a)^{-2}\int_{\mathbb R}
    (1+e^y)^{-(1-a)}(1+e^{v-y})^{-a}\dd y.
  \end{align*}
  Extend the last integrand by zero for $v\leq0$. Its support is convex, its logarithm is concave, and its two $y$-tails decay exponentially. Pr\'ekopa's integration and convolution theorems make $\phi_a$ and $k=\phi_{1/6}*\phi_{1/3}$ log-concave \cite[Theorems~6--7]{prekopa-1973-58032ec3}. They also give log-concavity of
  \[
    S(t)=\int_{\max(t,0)}^\infty k(x)\dd x
    =\int_0^\infty k(t+y)\dd y.
  \]
  The inequality $1+e^{v+s-y}\leq e^s(1+e^{v-y})$, followed by convolution and integration, gives
  \begin{equation}\label{sextic-truncated-integral-shift}
    S(t+s)\geq e^{-s/6}S(t),\qquad t\in\mathbb R,\ s\geq0.
  \end{equation}

  On $(0,\infty)$ the function $S$ decreases strictly from $1$ to $0$. For $t\in\mathbb R$ and $0<u<1$, set $x(t,u)=S^{-1}(uS(t))$ and $m(t)=\int_0^1x(t,u)\dd u$. Log-concavity makes $k/S$ nondecreasing, so for $t>0$
  \[
    \frac{\partial x}{\partial t}(t,u)
    =\frac{(k/S)(t)}{(k/S)(x(t,u))}\in[0,1].
  \]
  For $t\leq0$ the function $x(t,u)$ is constant in $t$. Thus $m$ is $1$-Lipschitz. If $u<v$, then $x(t,u)>x(t,v)$ and the derivative formula makes their gap nonincreasing in $t$. Changing variables by $S$ and expanding the square therefore gives, for real $s,t$,
  \begin{equation}\label{sextic-truncated-integral-bounds}
    \begin{aligned}
      \int_0^1(x(t,u)-m(s))^2\dd u
      &=\frac12\int_0^1\int_0^1
      (x(t,u)-x(t,v))^2\dd u\dd v+(m(t)-m(s))^2\\
      &\leq(\log I_0)''(0)+(t-s)^2.
    \end{aligned}
  \end{equation}
  The exponential integrability of $k$ justifies these integrals.

  Put $\ell(v)=e^{-v}(1+e^{-v})^{-2}$ and $s=\log r$. Since $\int_{-\infty}^{x-s}\ell(v)\dd v=(1+re^{-x})^{-1}$, Tonelli gives $I_r(0)=\int_{\mathbb R}\ell(v)S(s+v)\dd v$. Completing the square and applying \eqref{sextic-truncated-integral-bounds} at $t=s+v$ gives
  \begin{align*}
    (\log I_r)''(0)+\mleft(\frac{I_r'(0)}{I_r(0)}+m(s)\mright)^2
    &=\frac{\int_{\mathbb R}\ell(v)S(s+v)\int_0^1(x(s+v,u)-m(s))^2\dd u\dd v}{I_r(0)}\\
    &\leq(\log I_0)''(0)
    +\frac{\int_{\mathbb R}v^2\ell(v)S(s+v)\dd v}{I_r(0)}.
  \end{align*}
  Subtracting $(\log I_0)''(0)$ and dropping the displayed square gives
  \[
    2V(z)\leq\frac{\int_{\mathbb R}v^2\ell(v)S(\log r+v)\dd v}{I_r(0)}.
  \]
  Since $\ell$ is even and $S$ decreases, $I_r(0)\geq S(\log r)/2$. \zcref{sextic-truncated-integral-shift} gives $\ell(v)S(\log r+v)/I_r(0)\leq2e^{-5\abs*{v}/6}$, hence
  \[
    2V(z)\leq2\int_{\mathbb R}v^2e^{-5\abs*{v}/6}\dd v
    =\frac{1728}{125}<2\pi^2.
  \]
  This proves \eqref{sextic-radicand-bound}.
\end{proof}

\printbibliography

\end{document}